\documentclass{article}

\usepackage[T1]{fontenc}
\usepackage[utf8]{inputenc}
\usepackage[margin=1in,top=1in,bottom=3cm]{geometry}
\usepackage{perpage}
\MakePerPage{footnote}
\usepackage{epsfig}
\usepackage{latexsym}
\usepackage{amsmath}
\usepackage{mathrsfs}
\usepackage{yfonts}
\usepackage{makeidx}
\usepackage{graphicx}
\usepackage{slashed}
\usepackage{hyperref}
\usepackage[titletoc]{appendix}
\usepackage{chngcntr}
\usepackage{bm}
\usepackage{multicol}
\usepackage{nameref}
\usepackage{textcomp}
\usepackage{changepage}
\usepackage{lettrine}
\usepackage{enumitem}
\usepackage[sortcites,backend=bibtex,style=numeric,sorting=anyvt,doi=false,url=false,giveninits=true,isbn=false]{biblatex}
\AtEveryBibitem{\clearfield{month}}
\AtEveryBibitem{\clearfield{day}}
\usepackage{tikzsymbols}
\hypersetup{citecolor=blue,colorlinks=true,linkcolor=blue,pdfstartview={FitH},linktoc=page,pdftitle={Maxwell-Decay-on-DSRN-BH.pdf}pdfauthor={Mokdad Mokdad}}
\usepackage{amsthm,amssymb}
\usepackage{subcaption}
\newtheorem{thm1}{Theorem}
\newtheorem{prop1}[thm1]{Proposition}
\newtheorem{lem1}[thm1]{Lemma}
\newtheorem{cor1}[thm1]{Corollary}

\newtheorem*{thm}{Theorem}

 \newtheoremstyle{TheoremNum}
        {\topsep}{\topsep}              
        {\itshape}                      
        {}                              
        {\bfseries}                     
        {.}                             
        { }                             
        {\thmname{#1}\thmnote{ \bfseries #3}}
    \theoremstyle{TheoremNum}

     \newtheoremstyle{TheoremNum}
        {\topsep}{\topsep}              
        {\itshape}                      
        {}                              
        {\bfseries}                     
        {.}                             
        { }                             
        {\thmname{#1}\thmnote{ \bfseries #3}\thmnumber{}}
    \theoremstyle{TheoremNum}

\renewcommand{\d}{\mathrm{d}}
\newcommand{\scri}{{\mathscr I}}

\newcommand{\R}{\mathbb{R}}
\newcommand{\dl}{\partial}
\newcommand{\slnb}{\slashed{\nabla}}
\newcommand{\sldlt}{\slashed{\Delta}}
\newcommand{\hf}{\frac{1}{2}}
\newcommand{\dlr}{\partial_{r_*}}
\newcommand{\lie}{\mathcal{L}}
\newcommand{\hook}{{\setlength{\unitlength}{11pt}   
                   \begin{picture}(.833,.8)
                   \put(.15,.08){\line(1,0){.35}}
                   \put(.5,.08){\line(0,1){.5}}
                   \end{picture}}}
\makeindex

\title{Decay of Maxwell Fields on Reissner-Nordstrøm-de Sitter Black Holes}
\author{Mokdad Mokdad\thanks{email: mokdad.al.mokdad@gmail.com - Mokdad.Mokdad@univ-brest.fr} \\LMBA -- Université de Bretagne Occidentale}

\begin{document}

\maketitle

\begin{abstract}
	\emph{In this paper we use Morawetz and 
		geometric energy estimates -the so-called vector field method- to prove decay results for the Maxwell field in the static exterior region of the Reissner-Nordstr{\o}m-de Sitter black hole. We prove two types of decay: The first is a uniform decay of the energy of the Maxwell field on achronal hypersurfaces as the hypersurfaces approach timelike infinities. 
	The second decay result is a pointwise decay in time with a rate of $t^{-1}$ which follows from local energy decay by Sobolev estimates. 
	Both results are consequences of bounds on the conformal energy defined by the Morawetz conformal vector field. These bounds are obtained through wave analysis on the middle spin component of the field.
	The results hold for a more general class of spherically symmetric spacetimes with the same arguments used in this paper.}
\end{abstract}

\tableofcontents

\section{Introduction}

In this paper, we discuss the topic of decay, in particular, the decay of Maxwell fields in the exterior static region of Reissner-Nordstr{\o}m-de Sitter black holes (sometimes abbreviated as ``\textbf{RNdS}''). Our motivations are twofold: On the one hand, we use part of the decay results (uniform decay) obtained here to construct a complete conformal scattering theory in a separate paper \cite{mokdad_conformal_2016}. On the other hand, the subject of energy bounds and decay using Morawetz estimates in general relativity has gained attention in the last decades largely due to its fundamental role in the analysis of the nonlinear stability of spacetimes. Let us in short discuss these motivations.  

Linked to the decay problem is whether information about a test field and its history can always be retrieved from  its remnants or from its trace on the boundary of the observable space far in the future. In other words, we would like to know if the field can be completely characterized along with its entire evolution by observing its asymptotic profile in the distant future. The question is of course valid if we reverse the time orientation replacing the above experiments by similar ones about the past. In a complete scattering theory, the scattering operator is an isomorphism between two energy spaces on the future and the past boundaries. To obtain this complete characterization of test fields by their asymptotic profiles (their traces on the conformal boundary), we need to make sure that no information is lost at time-like infinities ($i^\pm$). Information is lost at time-like infinities means that part of the field is confined to a bounded space in  the exterior region for all future or past times. If no part of the energy is concentrated at $i^\pm$, which is to say that the energy flux across an achronal hypersurface decays to zero as the hypersurface approaches $i^+$ or $i^-$, then this guarantees that there is no loss. This is the uniform decay we prove here and we use in \cite{mokdad_conformal_2016}. For more on conformal scattering, we refer to the our paper \cite{mokdad_conformal_2016} and the references within.

On the other hand, global stability problems in the general theory of relativity require specific information about the asymptotic behaviour of the solutions to Einstein's equations. Often, a control provided by precise decay estimates for test fields on the background spacetime is crucial to access these informations. 
A basic problem in general relativity is the question of stability of Minkowski spacetime, that is, whether any asymptotically flat initial data set which is sufficiently close to the trivial one gives rise to a global (i.e. geodesically complete) solution of the Einstein vacuum equations that remains globally close to Minkowski spacetime. The local existence of solutions of the initial value problem was proven by Y. Choquet-Bruhat \cite{choquet-bruhat_theoreme_1952} in 1952. In 1983 partial results were obtained by H. Friedrich \cite{friedrich_cauchy_1983} using conformal methods, and in the early 1990's, the global nonlinear stability of Minkowski spacetime was established in the important work of D. Christodoulou and S. Klainerman  \cite{christodoulou_global_1993-1}\footnote{See also \cite{christodoulou_global_1993, christodoulou_stability_1999} for a summary of the proof. A revisit of the proof can be found in \cite{klainerman_evolution_2003}}.  The main tool they used for the energy estimates is the vector field method developed by Klainerman which generalizes the multiplier method in the works of C.S. Morawetz. They first obtain precise decay estimates \cite{christodoulou_asymptotic_1990} for the Bianchi equations (spin-2 zero rest-mass fields) on Minkowski which model linearized gravity on Minkowski spacetime. Then they prove that the same decay rates are still valid for the full Einstein equations. The perturbed spacetimes they construct has global features resembling those of Minkowski spacetime: a foliation by maximal spacelike slices given by the level hypersurfaces of a time function; an optical function whose level hypersurfaces describe the structure of future null infinity; a family of almost Killing and conformal Killing vector fields related to the time and optical functions. They use symmetries and almost symmetries to get conserved and almost conserved quantities and to define the basic energy norms. These symmetries and almost symmetries are generated by the almost Killing and conformally Killing vector fields. These vector fields are the generator of time translation $T$, the generators of the Lorentz group ( generators of rotations and boosts) $\Omega_{ij}$ for $i,j=0,\dots, 3$, the generator of scaling (dilation) $S$, and the generator of inverted time translation (conformal Morawetz vector field) $K$. In fact, the Lie derivatives along these vector fields are used to define the basic quantities, which give better control in the estimates. This work was the first step towards the proof of the stability of Minkowski spacetime, which is a crucial question for the understanding of the large time evolution of black holes.
Currently, many groups are concentrating on the stability of Kerr black holes\footnote{Works addressing the question of the stability of the Schwarzschild manifold can be found in \cite{dimock_classical_1987, kay_linear_1987, finster_linear_2006, regge_stability_1957,holzegel_ultimately_2010, dafermos_proof_2005}.}. Asymptotically flat vacuum initial data for the evolution problem in general relativity are expected to give rise to spacetimes that can be decomposed into regions each of which approaches a Kerr black hole. The Kerr black hole spacetime is expected to be the unique, stationary, asymptotically flat, vacuum spacetime containing a nondegenerate Killing horizon\footnote{A Killing horizon is a null hypersurface defined by the vanishing of the norm of a Killing vector field.} \cite{alexakis_uniqueness_2010}. This is relevant in the context of some main  problems of the theory, such as the weak cosmic censorship conjecture. Proving the Kerr black hole stability is a major step towards solving these problems. The multiplier method, the vector field method, and its generalizations, are being employed to obtain the required uniformly bounded energies and to prove Morawetz estimates for solutions of the wave equation on black hole spacetimes, motivated by the fact that, proving boundedness and decay in time for solutions to the scalar wave equation on the asymptotically flat exterior of the Kerr spacetime, is an important model problem for the full black hole stability problem. However, there are some fundamental difficulties in the Kerr case, mostly because of the lack of symmetries, the trapping effect ranging over a radial interval, and there is no positive conserved quantity since Kerr black holes do not admit global timelike Killing vector fields.

\subsection*{An Overview of Decay}

Many of the decay results in the literature are for solutions of wave equations. Let alone that waves themselves are interesting as physical and mathematical objects, the reason behind the extensive study is that these prototype equations are not only important model problems, but also appear as a fundamental part of the structure in many problems where their analysis is essential. For example, they model the propagation of several systems including Schrödinger, gravity, and of course Maxwell's equations that we are studying in the present work, in which as we shall see, a wave analysis on the middle spin component of the field plays a central role in obtaining the decay results. We present an overview on the history of decay estimates summarizing some methods used to obtain them and how they evolved to become more adaptable to different geometries.

\subsubsection*{Basic Notions}

\begin{figure}
	\centering
	\includegraphics[scale=1]{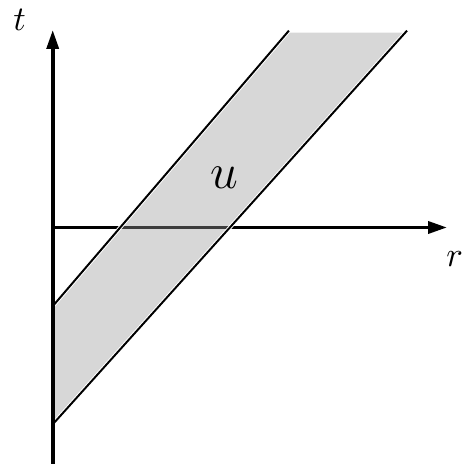}
	\caption{\emph{The support of $u$ is contained in this diagonal strip and therefore $u$ propagates exactly at speed 1. Here $r=|x|$.}}
	\label{fig:Huygensprinciple}
\end{figure}

Consider the simple scalar wave equation on $\R^{1+3}$,

\begin{equation*}
\dl^2_t u -\dl^2_{x_1} u - \dl^2_{x_2} u -\dl^2_{x_3} u =0 \; . 
\end{equation*}
The explicit formula for solutions to the wave equation in one space dimension was due to D'Alembert. In three space dimensions, the wave equation admits radial solutions of the form $$u=\frac{h(|x|-t)}{|x|}$$ where $x\in\R^3$, $|x|=(x_1^2+x_2^2+x_3^2)^\hf$, and $h$ is any twice differentiable function on $\R$. Let $u$ be such a solution and say $h$ is a smooth function which is compactly supported in $]0,+\infty[$. This solution radiates away form the origin at speed 1 as $t$ increases, and for some $R>0$, it identically vanishes\footnote{Actually these particular solutions, i.e. $h$ having such support, vanish for $t\ge |x|$.} for $t>|x|+R$ (figure \ref{fig:Huygensprinciple}). Such a solution models a disturbance starting in a bounded region which then spreads outward and reaches every point in space, but for each point and after a finite amount of time, there is no disturbance left at all. Perhaps not seen as directly as in the previous simple case, but in fact, this is true for all the solutions of the above wave equation on $\R^{1+3}$ which start in confined regions.  This infinite fall off rate follows from Kirchhoff's formula (19th century) which can be proved by the method of spherical means. An equation having this property is said to satisfy the strong Huygens principle:
\begin{thm}[Huygens Principle]
	If the initial data, $(u(0,x),\dl_t u(0,x))$ with $x\in \R^3$, for the above wave equation are supported in the ball $B(0,R)$, then the associated solution $u$ satisfies,
	\begin{equation*}
	u(t,x)=0 \qquad for ~ all ~ |t|>|x|+R
	\end{equation*}
\end{thm}
For equations satisfying the strong Huygens principle, if we start from compactly supported initial data, the field decays infinitely fast in time because at each point in space it vanishes identically after a certain time. Despite this \emph{pointwise decay}, there is a quantity determined by the solution which is conserved for all times due to the time translation symmetry of the system. This can be seen by multiplying the wave equation with the time derivative of the solution, called the  \emph{multiplier}, and rearranging,
$$\dl_t u \left(\dl^2_t u -\sum_{i=1}^3 \dl^2_{x_i} u\right)=\hf \dl_t\left((\dl_t u)^2+\sum_{i=1}^3(\dl_{x_i} u)^2\right)-\sum_{i=1}^3 \dl_{x_i}(\dl_t u \dl_{x_i} u),$$
if we now integrate the right hand side over a spacetime slab $[t_1,t_2]\times \R^3$ and using the fact that $u$ is a solution for the equation and that it has a compact support in space for all $t$, we arrive at the following identity:
$$\int\limits_{\{t_1\}\times \R^3}\left((\dl_t u)^2+\sum_{i=1}^3(\dl_{x_i} u)^2\right) \d x=\int\limits_{\{t_2\}\times \R^3}\left((\dl_t u)^2+\sum_{i=1}^3(\dl_{x_i} u)^2\right) \d x \; .$$
This quantity is called the (total) \emph{energy} $E[u](t)$ and it is conserved: $E[u](t)=E[u](0)$. However, the \emph{local energy} 
$$E[u](D,t)=\int\limits_{\{t\}\times D}\left((\dl_t u)^2+\sum_{i=1}^3(\dl_{x_i} u)^2\right) \d x$$
in any bounded region of space $D=\{x\in\R^3 ; |x|\le R\;  \}$ is clearly not conserved and becomes zero after the wave leaves the region $D$.  

The strong Huygens principle for the wave equation on flat spacetime is only valid in odd space dimensions starting at three. More general wave equations with a potential or on curved spacetimes satisfy a weak Huygens principle which says essentially that the local energy decays. One may then ask at what rate the local energy decays for, say, smooth compactly supported data. By a rate of decay in time for the local energy we mean a function $d(t)$ that tends to zero when $t$ tends to infinity and such that,
$$E[u](D,t)\le d(t) E[u](0)\; .$$
For example, in two space dimensions the rate of pointwise decay of solutions to the above wave equation can be exactly $t^{-1}$, and thus the local energy decays as $t^{-2}$. The question is then when to expect that a solution should tend pointwise with time to zero and at which rate. The main obstruction to decay is the existence of finite energy stationary solutions, i.e. of the form $e^{i\lambda t}v(x)$. Provided we avoid such solutions, vector field and multiplier methods can be applied to obtain decay rates in fairly general situations.

\subsubsection*{The Multiplier Method}

The method of multipliers originated from the so-called Friedrichs' ABC method that dates back to K.O. Friedrichs in the 1950's. The method was used first to obtain decay and energy estimates in non-relativistic situations of geometrical optics, possibly outside of an obstacle whose geometry is known, and then was later applied to relativistic theories. The idea of this method is to multiply the equation with a factor $Mu$, where $M$ is a linear first-order differential operator, defined as $$Mu = Au + B\cdot \nabla u + C\dl_t u \; .$$ Then we try to express the product as a divergence or at least as an identity of the form $$divergence~term ~+~ remaining~terms = 0\; .$$ Then if we integrate over a domain in $\R^{1+n}$, the required estimates are derived by controlling the remainder. The multiplier method generalizes the ABC method: Suppose $L$ is a differential operator of order $k$ and consider the expression $Mu  Lu$, where $M$ is a differential operator of order $k-1$. With the right mix of derivatives, one hopes that $Mu Lu$ can be written as $\mathrm{div}(Qu) + Ru$, where $Qu$ and $Ru$ are quadratic expressions in the derivatives up to order $k-1$. The method of multipliers was used in the 1960's and 1970's to prove uniform decay results for the homogeneous linear wave equation ($\square u = 0$) outside obstacles. C.S. Morawetz was the first to succeed in proving local energy decay for star-shaped obstacles with Dirichlet boundary condition using this method  in 1961 \cite{morawetz_decay_1961}. In this work, the effects of scaling and of the spreading into space on the solution for the wave equation and its local energy, is captured using the scaling multiplier $$Su=t\dl_t u +r\dl_r u + u \; ,$$ and the following local energy decay is established
$$E[u](R,t)\le \frac{C}{t} E[u](0)\; ,$$ where $R$ is a region bounded between the obstacle and an outside sphere, and $C>0$ depends on the obstacle and the support of the initial data. This estimate then gives a pointwise decay of rate $t^{-\hf}$. A year later, Morawetz used the multiplier $$K=\xi^2\dl_\xi u + \eta^2 \dl_\eta u + (\xi+\eta)u= (t^2+r^2)\dl_t u + 2tr\dl_r u + 2tu \qquad ; \; \xi=t-r \; , \eta=t+r\; ,$$ in her work \cite{morawetz_limiting_1962} to improve on the results of \cite{morawetz_decay_1961} and get faster decay rates of $t^{-1}$ for the pointwise decay and $t^{-2}$ for the local energy. She was motivated by the fact that for large times the disturbance is expected to be radiating outwards, and there will be little dependence on the angles. So, $ru$ will approach a solution of $\dl_r^2 w - \dl_t^2 w = 0$ for which an appropriate multiplier is $Nw = p_1 \dl_r w + p_2 \dl_t w$. The multiplier $K$ is in fact related to a ``time'' translation: If we apply the Kelvin transformation on the coordinates $(t,r,\theta,\varphi)$ given by,
$$\hat{t}=\frac{t}{r^2-t^2}\; ,\quad \hat{r}=\frac{r}{r^2-t^2}$$ and leaving the angular variables the same, we can see that  $$\dl_{\hat{t}} = 2tr\dl_r + (r^2+t^2)\dl_t\; .$$ This transformation is conformal and takes the cone $r^2=t^2$ at the origin to a cone at infinity and vice verse. It is not a surprise then that this vector field is appropriate for studying the asymptotic behaviour of the solution.    
Moreover, in 1968 \cite{morawetz_time_1968} Morawetz used a radial multiplier of the form $$\zeta(r)(\dl_r u + r^{-1}u) $$ where $-\zeta(r)$ is a bump function around the origin, to obtain uniform integrated local energy estimates for the non-linear Klein-Gordon equation $\square u + mu + P(u)=0$,
$$\int_0^T E[u](\Omega,t)\d t\le K E[u](0)\; ,$$ where $\Omega$ is a finite region in space and $K$ a positive constant depending only on $\Omega$ (or its volume). She then uses this estimate to prove that the local energy decays but without giving a rate. Also in the same paper, it is proven that the $L^2$-norm of the solution decays. Before this work, also in  1968, a similar but complex radial multiplier was used by C.S. Morawetz and D. Ludwig \cite{morawetz_inequality_1968} on a wave operator.

These multipliers and their corresponding vector fields have all found many important applications, most notably in General Relativity\footnote{Of course, the results of Morawetz in other fields were as important, especially in the field of geometrical optics, and have been built upon and improved: Better decay rates have been achieved, as in odd dimensions $n \ge 3$, Huygens principle has been shown to imply an exponential rate of decay whenever there is some sort of decay by P.D. Lax, C.S. Morawetz, and R.S. Phillips in 1963 \cite{lax_exponential_1963}, and then by Morawetz \cite{morawetz_exponential_1966} in 1966. Moreover, the class of obstacles under consideration has been enlarged using the method of multipliers after generalizing the multipliers to suit the geometry of the obstacle. Other wider generalizations  later followed: W.A. Strauss \cite{strauss_dispersal_1975} proved uniform local energy decay for the homogeneous linear wave equation using the Straussian vector fields. Then these Straussian vector fields were generalized by C.S. Morawetz, J.V. Ralston, and W.A. Strauss \cite{morawetz_decay_1977},  by constructing a pseudo-differential operator $P(x, D)$ (coming from a function $p(x,\xi)$ called the ``escape function''.), and finally setting $P u$ as a multiplier.\label{footnoteMorawetz}}. We mention an interesting work of Morawetz and W.A. Strauss \cite{morawetz_decay_1972} on decay and scattering for a nonlinear relativistic wave equation using these methods. Morawetz also established decay properties for Maxwell fields in \cite{morawetz_notes_1975}. For more on her work, we also refer to \cite{morawetz_energy_1966}.

\subsubsection*{Vector Field Method}

The vector field method is a flexible tool generalizing the multiplier method by making use of well adapted vector fields, related to symmetries or approximate symmetries of the equations, to derive decay estimates and thus to control the long time behaviour of solutions. The basics of this method has two aspects: The vector fields are used to define generalized energy norms, and, if they commute with the equations then one can derive identities for the energy norms considered. In the mid 1980's, S. Klainerman introduced the notion of generalized energy norms defined from the conformal group, which is generated by the vector fields $T,S,K,$ and the $\Omega_{ij}$'s and whose elements have useful commuting properties among themselves and with the D'Alembertian. He used them to obtain energy estimates and prove decay for solutions of the wave equation on $\R^{n+1}$ \cite{klainerman_null_1986,klainerman_remarks_1987,klainerman_uniform_1985}. These works of Klainerman were in essence a combination of the  local energy  decay  estimates of C.S. Morawetz \cite{morawetz_limiting_1962} and the conformal method of Y. Choquet-Bruhat and D. Christodoulou \cite{choquet-bruhat_existence_1981}. If $\mathbb{A}$ is a set of vector fields and $s\in \mathbb{N}$ we define the following norm of a function $u$ on $\R^{n+1}$ by 
$$\Vert u(t)\Vert_{\mathbb{A},s,p}=\left(\sum_{k=0}^s \sum_{X_{i_j}\in \mathbb{A}}\int_{\R^{n}}\vert(X_{i_1} \dots X_{i_k} u)(t,x)\vert^p \d x\right)^{\frac{1}{p}}\; .$$
Klainerman uses such norms for different subsets of the conformal group in place of $\mathbb{A}$ to get what he calls global Sobolev inequalities (which are now known as Klainerman-Sobolev inequalities) of the form 
$$\vert u(t,x)\vert\le h(t,|x|)\Vert u(t)\Vert_{\mathbb{A},s,p}$$ for functions $u$ with $$\Vert u\Vert_{\mathbb{A},s,p}^{\#}=\sup_{t\ge 0}\Vert u(t)\Vert_{\mathbb{A},s,p}<\infty\; .$$
In the same papers he also gets decay estimates of the form 
$$\vert u(t,x)\vert\le d(t)\Vert u\Vert_{\mathbb{A},s,p}^{\#}\; .$$
Many results concerning the long-time and global existence were subsequently obtain using the methods of Klainerman. Klainerman himself used the results we hinted at above to prove long-time existence for a family of nonlinear wave equations \cite{klainerman_uniform_1985}. And using the same methods, he obtained existence and decay results (of rates $t^{\frac{-5}{2}}$) for nonlinear Klein-Gordon equations on Minkowski spacetime \cite{klainerman_global_1987}. Other works in the domain include L. H\"{o}rmander \cite{hormander_lifespan_1987} in 1987 on nonlinear hyperbolic equations; A. Bachelot \cite{bachelot_probleme_1988} in 1988 on Dirac-Klein-Gordon systems; J. Ginibre, A. Soffer, and G. Velo \cite{ginibre_global_1992} in 1992 for the critical non-linear wave equation; and of course, the important work of D. Christodoulou and Klainerman on the stability of Minkowski spacetime \cite{christodoulou_global_1993-1}.

Another important work of D. Christodoulou and S. Klainerman in 1990 is their paper \cite{christodoulou_asymptotic_1990} which studies the asymptotics of linear field equations in Minkowski spacetime. This paper was in fact the preparatory foundation for the proof of the nonlinear stability of Minkowski spacetime and in it, the vector field method took its standard current form. This method was used in many works on decay estimates that came later, and initiated, along with the proof of the stability of Minkowski spacetime, the project of proving Kerr stability. In \cite{christodoulou_asymptotic_1990} they derive uniform decay estimates for solutions to linear field equations in Minkowski spacetime which give
precise information on the asymptotic behaviour of the solutions. It is based on geometric considerations of energy and generalized energy estimates. Their method relies on Klainerman's systematic use of the invariance properties of the field equations with respect to the conformal group of the Minkowski spacetime\footnote{Thus differing from the previous methods of analysing the fundamental solution.}, and was then extended to nonlinear cases, in particular to Einstein's vacuum equations \cite{christodoulou_global_1993-1}. 

The usefulness of the vector field method is best seen, although not exclusively\footnote{See for example \cite{andersson_hidden_2015}.}, in view of Noether's theorem in the case of general field equations derived from a quadratic action in the context of a Lagrangian theory. Let $\phi$ be a general field on a general spacetime $(M,g)$ and assume there is a scalar Lagrangian $L$ which depends on the field and its derivatives and possibly position in spacetime. $L$ is used to define an action $S$ as the integral of $L$ on $M$. The field equations governing the behaviour of the field are derived by the ``principle of least action'', that is to say that $\phi$ satisfies the field equations if it is a minimizer (or a critical point) of the action. These field equations are then a simple relation between the variation of the Lagrangian with respect to the field and its variation with respect to the field's derivatives, and are called the Euler-Lagrange equations:
$$\frac{\delta L}{\delta \phi}-\nabla^a \frac{\delta L}{\delta \nabla^a \phi}=0\; .$$ 
One can then define from the field and the Lagrangian a symmetric 2-tensor $\mathbf{T}$ called the {energy-momentum tensor} (or stress-energy tensor) depending on the field and its derivatives (usually quadratically), which by the Euler-Lagrange equations turns out to be divergence-free\footnote{Although the natural way of obtaining an energy-momentum tensor is by means of a Lagrangian, one can as well directly consider 2-tensors with the desired properties and which might not be derived from a Lagrangian.}(see \cite{hawking_large_1973,wald_general_2010} for example)\footnote{For electromagnetic fields represented by 2-forms on the manifold, we actually vary the (local) potential and not the field (the 2-form) itself.}. The energy associated with a vector field $X$ and evaluated on a hypersurface $\Sigma$ is,
$$E_X(\Sigma)=\int_{\Sigma} \mathbf{T}_{ab}X^b \d \sigma^a\; , $$
where $\alpha_a \d \sigma^a$ is the 3-form $\star \alpha$ given by the Hodge star operator for any 1-form $\alpha$. $E_X(\Sigma)$ is sometimes referred to as the geometric energy. If $\Sigma$ is spacelike and $M$ is time-orientable we choose the normal on $\Sigma$ to be future-oriented\footnote{Depending on the sign conventions.}, since if the energy-momentum tensor satisfies the {dominant energy condition}:
$$\mathbf{T}(V,W)\ge 0 \; , \quad \mathrm{whenever}~V~\mathrm{and}~W~\mathrm{are~future~oriented,~causal~vectors},$$
the above expression of the energy will be positive definite if $X$ is timelike and  future-oriented. When the spacetime is globally hyperbolic\footnote{Admitting a global Cauchy hypersurface. See \cite{hawking_large_1973}.} or foliated by hypersurfaces $\Sigma_t$ of constant time, then by Stokes' theorem, or more precisely by the divergence theorem, we have the following law: If $\Omega_{t_1,t_2}$ is the region enclosed between $\Sigma_{t_1}$ and $\Sigma_{t_2}$, then by the properties of $\mathbf{T}$ we have,
$$E_X(\Sigma_{t_2})-E_X(\Sigma_{t_1})=\int_{\Omega_{t_1,t_2}} \nabla^a(\mathbf{T}_{ab}X^b) \d \mathrm{Vol}_g=\int_{\Omega_{t_1,t_2}} {}^{(X)}\pi_{ab}\mathbf{T}^{ab} \d \mathrm{Vol}_g \; ,$$
where ${}^{(X)}\pi_{ab}=\lie_X g_{ab}=2\nabla_{(a} X_{b)}$ is called the deformation tensor of $X$. This law is called the deformation law.  A vector field $X$ is a conformal Killing vector field if the deformation tensor of $X$ is proportional to the metric by a scalar factor $\lambda$, and $X$ is a Killing vector field when $\lambda=0$. We see then that when $X$ is Killing the deformation law entails that the energy is conserved. The same happens when $X$ is conformal Killing \emph{and} the energy-momentum tensor is trace-free. In general, energy estimates are obtained by controlling the deformation term ${}^{(X)}\pi_{ab}\mathbf{T}^{ab}$, and in that case one says that one has an (almost) conservation law. A symmetry operator for an equation is defined to be a differential operator that takes solutions to solutions; in simple cases, the symmetry operator commutes with the equations. When $Y$ is Killing, the Lie differentiation $\lie_Y$ with respect to $Y$ is a symmetry operator for the wave equation and for Maxwell's equations among others. This means that when $\mathbb{Y}$ is a set of Killing vector fields, one has identities for the energies defined using these vector fields and also for all Lie derivatives of the solutions with respect to these vector fields, at all orders. This adds on the control of the energies and allows better estimates and rates of decay.

In Minkowski spacetime, the conformal group is generated by conformal Killing vector fields, but only time and time inverted translation (or time acceleration as called in \cite{christodoulou_asymptotic_1990}) generators are timelike and thus can be associated with a positive definite energy. D. Christodoulou and S. Klainerman use arguments similar to the one above with the symmetry generators $T,S$, and $K$ in  \cite{christodoulou_asymptotic_1990} to obtain uniform bounds on the generalized energies and then, by means of Klainerman's global Sobolev inequalities, obtain the decay estimates for Maxwell and spin-2 equations. The latter are formally identical to the  Bianchi identities for the Riemann curvature tensor and thus relevant to the understanding of the Einstein field equations. In fact, the methods they developed in \cite{christodoulou_asymptotic_1990} in the study of the spin-2 equations in Minkowski spacetime prepared for the subsequent study of the nonlinear stability of the Minkowski metric, as mentioned at the beginning of this introduction.

\subsubsection*{Some Recent Works}

The literature centred around decay estimates in general relativity is vast, so we  refer to some recent works where additional references can be found.  In particular, Blue's paper \cite{blue_decay_2008} about the decay of Maxwell fields in Schwarzschild in 2008 is central to our work, in fact, we show that the methods used in \cite{blue_decay_2008} can be applied to the case of RNdS black holes. Furthermore, this work shows that already existing methods of vector fields and Morawetz estimates can be applied to generic spherically symmetric black holes including the case of positive cosmological constant, with no real modifications (see section \ref{genericsphericalf}). In their paper of 1999 on nonlinear Schr\"{o}dinger equation \cite{laba_global_2000}, I. Laba and A. Soffer introduced a Morawetz vector field on the Schwarzschild spacetime. They also introduce a modified radial Morawetz multiplier, known as Soffer-Morawetz multiplier, based on the work of C.S. Morawetz, J.V. Ralston, and W.A. Strauss \cite{morawetz_decay_1977} (also see footnote \ref{footnoteMorawetz}). Through the 2000's, these tools were used on Schwarzschild's spacetime with further adaptations in the works of P. Blue, A. Soffer, and J. Sterbenz \cite{blue_semilinear_2003,blue_decay_2008,blue_wave_2005, blue_errata_2006, blue_uniform_2006, blue_decay_2004} and in this present work, to help control the trapping terms. The radial Morawetz vector field that made these estimates possible is centred about the orbiting null geodesics. In 2000, in a paper on Maxwell fields on Schwarzschild's spacetime \cite{inglese_asymptotic_2000}, W. Inglese and F. Nicolò give specific asymptotic estimates for different components of the field. A variant of the problem considered by P. Blue and J. Sterbenz in 2006 \cite{blue_uniform_2006}, about the uniform decay of local energy for wave equations, was independently studied by M. Dafermos and I. Rodnianski \cite{dafermos_red-shift_2009} in 2009 with a stronger estimate obtained near the event horizon (see also \cite{dafermos_lectures_2008} by the same authors). M. Dafermos and I. Rodnianski in 2008 proved decay results for the wave equation on Schwarzschild-de Sitter spacetimes \cite{dafermos_wave_2008}. The same authors proved uniform boundedness for the wave equation on slow Kerr backgrounds in 2011 \cite{dafermos_proof_2011}. In the same year, and using the same methods, D. Tataru and M. Tohaneanu  obtained local decay for energy also on Kerr \cite{tataru_local_2011}, and later in 2013 D. Tataru extended the results to asymptotically flat stationary spacetimes \cite{tataru_local_2013}. A paper by J.-F. Bony and D. H\"{a}fner in 2008 \cite{bony_decay_2008} addresses the decay and non-decay of the local energy for the wave equation on the de Sitter-Schwarzschild metric.  Several decay estimates with rates were obtained in the early 2010's: J. Luk in 2012 \cite{luk_vector_2012}; M. Tohaneanu \cite{tohaneanu_strichartz_2012}; M. Dafermos and I. Rodnianski in 2010 \cite{dafermos_decay_2010, dafermos_new_2010} and in 2014  with Y. Shlapentokh-Roth\-man \cite{dafermos_decay_2014}. There is also a paper in 2013 by L. Andersson P. Blue, and J.-P. Nicolas on wave equations with trapping and complex potential that appear in the Maxwell and linearized Einstein systems on the exterior of a rotating black hole \cite{andersson_decay_2013}.  Two recent papers in 2015 were published by L. Andersson and P. Blue: \cite{andersson_uniform_2015} proving uniform energy bounds for Maxwell fields on Schwarzschild, and \cite{andersson_hidden_2015} in which they generalize the vector field method to take the hidden symmetries of Kerr spacetime into account (also see \cite{andersson_second_2014} for second order symmetry operators) and obtain an integrated Morawetz estimate and uniform bounds for a model energy for the wave equation. M. Dafermos, I. Rodnianski, and Y. Shlapentokh-Rothman's work on scattering for the wave equation on Kerr \cite{dafermos_scattering_2014} contains decay results and the uniform energy equivalence needed for conformal scattering (see J.-P. Nicolas \cite{nicolas_conformal_2015}). There is a more recent paper by L. Andersson, T. B\"{a}ckdahl, and P. Blue  \cite{andersson_decay_2016} in 2016 proving a new integrated local energy decay estimate for Maxwell fields outside a Schwarzschild black hole using a new superenergy tensor ${{\bf{H}}}_{{ab}}$ defined in terms of the Maxwell field and its first derivatives. There have been works on Price's law (see \cite{price_nonspherical_1972, price_nonspherical_1972-1}), such as \cite{metcalfe_prices_2012} in 2012 by J. Metcalfe, D. Tataru, and M. Tohaneanu. Finally, using different techniques (an integral representation of the propagators, see \cite{finster_integral_2005}) F. Finster, N. Kamran, J. Smoller, and S.-T. Yau obtain decay estimates for: Dirac on the Kerr-Newman spacetime \cite{finster_decay_2002} in 2002, and \cite{finster_long-time_2003} in 2003; for the wave equation on Kerr \cite{finster_decay_2006, finster_erratum_2008} in 2006 (corrected in 2008), and by  F. Finster and J. Smoller
\cite{finster_time-independent_2008} in 2008 also for the wave equation on Kerr.

\subsubsection*{Maxwell Fields}

The Maxwell field is a 2-form $F$ on the spacetime satisfying Maxwell's equations:
$$\d F=0 \qquad;\qquad \d \star F=0$$ 
where $\d$ is the exterior differentiation and $\star$ is the Hodge star operator, or in abstract index notation,
\begin{eqnarray}
\nabla^a F_{ab} &=& 0 \; , \\
\nabla_{[ a} F_{bc ] }&=& 0 \; ,
\end{eqnarray}
Interest in the asymptotic behaviour of solutions to Maxwell's equations goes back at least to the 1970's \cite{morawetz_notes_1975}, yet most of the literature is on scalar wave equations. It turns out that some features of the Maxwell field can be captured in the behaviour of its components which are governed by wave-like equations, and results on the latter can be applied to study Maxwell systems. The behaviour of Maxwell fields is well-known in flat spacetime, at any point in space the effect of a signal dies off. But the total energy carried by the signal is preserved, carried off in fact to infinity, as seen for example in the works of C.S. Morawetz in 1974 \cite{morawetz_decay_1974}, and D. Christodoulou and S. Klainerman  \cite{christodoulou_asymptotic_1990} in 1990 with rates of $t^{-5/2}$ obtained using the full conformal group. In Schwarzschild,  a rate of $t^{-3}$ was obtained in regions bounded away from the horizon and null infinity, by R.H. Price in 1972 \cite{price_nonspherical_1972-1}, and later by R.H. Price and L.M. Burko in 2004 \cite{price_late_2004}. Only time and spherical symmetries are available in this case, so the vector field method produces a slower rate of $t^{-1}$, as in P. Blue work \cite{blue_decay_2008}, however, the conformal energy associated to the conformal Morawetz vector field can be used to control all the components of the field, and no spherical harmonic decomposition is required. We prove here that this is also the case for generic spherically symmetric static black holes by working out the details on RNdS black holes; the results can be extended to more general situations including cosmological black holes. In 2015 J Sterbenz and D Tataru \cite{sterbenz_local_2015} obtained local energy decay for Maxwell fields on a general spherically symmetric spacetime but which is required to be asymptotically flat, thus they do not cover cases with positive cosmological constant. Also in 2015, J. Metcalfe, D. Tataru, and M. Tohaneanu studied the pointwise decay properties of solutions to the Maxwell system on a class of non-stationary asymptotically flat backgrounds \cite{metcalfe_pointwise_2014}. Decay of waves and non-scalar fields (including Maxwell) on cosmological backgrounds with a de Sitter character were recently treated in the works of A. Vasy and P. Hintz \cite{hintz_analysis_2015, hintz_global_2015, hintz_semilinear_2015} in 2015, using methods from microlocal analysis and it seems that their work needs positive cosmological constant (maybe with the exception of flat spacetime), whereas the vector field method which we use applies equally well with or without a cosmological constant. 

Before discussing the method we use, it is worth mentioning that there is a resemblance between Maxwell's equations and the spin-2 equations. A spin-2 field can be seen as a covariant 4-tensor with the following symmetries:
\begin{eqnarray*}
	W_{abcd}&=&-W_{bacd} \; , \\
	W_{abcd}&=&-W_{abdc}\; ,\\
	W_{[abc]d}&=& 0\; ,\\
	{W_{abc}}^a&=& 0 \; ,\\
\end{eqnarray*}
satisfying the equations,
\begin{eqnarray*}
	\nabla^a W_{abcd}&=&0 \; ,\\
	\nabla_{[e}W_{ab]cd}&=&0 .\\
\end{eqnarray*}
The symmetries of a spin-2 field extend the antisymmetry of a Maxwell field,
and the two systems of equations have similarities. If the Einstein vacuum equations are satisfied, then the Ricci curvature vanishes, and the Weyl curvature satisfies the spin-2 field equations. In Minkowski spacetime, the spin-2 field equations models the linearization of Einstein's equations about the Minkowski solution. If one
introduces a perturbed metric on Minkowski spacetime and treats the Weyl tensor as
a tensor field on the original space-time, then, using the flatness of the background and the vanishing of the Christoffel symbols in cartesian coordinates, the difference between the covariant
derivative of the Weyl tensor with respect to the perturbed metric and the original
metric will be second order in the perturbation. Thus, ignoring second and higher order terms, the perturbed Weyl tensor satisfies the spin-2 field equations on the original metric. In this sense, the spin-2 field equations are the linearization of the Einstein vacuum equation about Minkowski spacetime. This is the motivation for studying the spin-2 field in \cite{christodoulou_global_1993-1}. However, this is not true for the linearization about other solutions. When linearizing around a curved space-time, the Christoffel symbols do not vanish, and the linearized Einstein equations do not reduce to the spin-2 field equations. Nevertheless, we expect that an analysis using the vector field method and Morawetz estimates will apply to the linearized gravity system. The linearized gravity equations are more complicated than the spin-2 field equations because there are terms involving the perturbed Christoffel symbols contracted against the unperturbed and non-vanishing Weyl tensor.

The arguments used in our work follow the same philosophy as in the works \cite{blue_phase_2009, blue_uniform_2006, blue_decay_2008,dafermos_red-shift_2009} using the vector field method.  The major obstacle is the trapping effect:

\paragraph{Trapping Effect.} The conformal vector field $$K=(t^2+r_*^2)\dl_t + 2tr_* \dl_{r_*}$$ where $r_*$ is the Regge-Wheeler coordinate, is timelike away from the $t=\pm r_*$ hypersurfaces where it is null. It is used to introduce a positive definite quantity, a conformal energy. This quantity is not conserved because of trapping. The presence of null geodesics at the photon sphere manifests itself through the trapping terms which are positive around the photon sphere. They appear as a contribution governing the growth of the conformal energy. It can be seen as the main ``error'' which is generated by the divergence of the conformal energy density. This effect can be overcome by introducing a radial vector field which points away from the photon sphere. This is a modified Morawetz radial multiplier of the form $A\dl_{r_*}$, where $A$ is a continuously differentiable function of $r_*$ that changes sign at the photon sphere, marked at $r_*=0$.

The work can be divided into three main steps. In the first step, the conformal energy, defined by the conformal Morawetz vector field, of a Maxwell field is not conserved but can be controlled by the conformal charge\footnote{This is the conformal energy of the solution to the wave equation, but to avoid confusion with the conformal energy we call it a conformal charge.} of the middle (or spin-weight zero) component of the field which satisfies a wave-like equation decoupled from the other components. This reduces the problem from spin-1 to spin-0, this is the so-called ``spin-reduction''. 

\paragraph{Wave Analysis.} The conformal charge of the solutions to the wave-like equation is not conserved either. The second step is to control the error term using a radial Soffer-Morawetz multiplier which allows us to obtain a uniform bound on the conformal charge of the wave. Because this wave-like equation is actually simpler than the covariant wave equation, the usual analysis on the local energy of the wave equation is replaced  by an analysis of an energy localized inside the light cone, and no decomposition on the spherical harmonics is required. Through some Hardy estimates, the trapping term is controlled by the energy generated by the radial multiplier and the integral of the energy localized inside the light cone. Since the trapping term controls the growth of the conformal charge, and since the energy (generated by time translation) is conserved, this gives a linear bound on the conformal charge. Using the Cauchy-Schwarz estimate and an integration by parts, the linear bound is improved to a uniform one. This also gives a uniform bound on the trapping term.

The third step is to use the conformal energy to control norms of the Maxwell field. The generalized energy and conformal energy of the Maxwell field generated by the rotation group are conserved and control the energy and the conformal charge of the middle spin component which in turn control the trapping term by the uniform bound. Thus, the conformal energy of the Maxwell field is controlled by the generalized energy and conformal energy of the initial data through a uniform bound. Since the integral of the trapping term has been controlled in the entire $r_*$-range, we have a uniform bound on the energy flux through any achronal hypersurface. This can be improved to a uniform decay rate of $t^{-2}$. The integrand in the conformal energy behaves like $t^2$ times the Maxwell field components squared. Since the conformal energy is bounded, the field components decay in $L^2_{\mathrm{loc}}$ like $t^{-1}$. Control on radial derivatives is the main thing that we need to improve this into pointwise decay. Sobolev estimates can be used to convert $L^2_{\mathrm{loc}}$ decay for derivatives into $L^{\infty}_{\mathrm{loc}}$ decay. For this, we need decay on the spatial derivatives of the Maxwell field. From spherical symmetry, the Lie derivative of the Maxwell field along angular Killing vectors also satisfies the Maxwell equations and has the same type of local $L^2$ decay as the field. Since differentiating in the radial direction does not generate a symmetry, the Lie derivative in that direction will not solve the Maxwell equations. To control the radial derivatives, we use the structure of the Maxwell equations. Using the time translation symmetry, we can control Lie time derivatives in $L^2_{\mathrm{loc}}$. In a fixed and compact range of $r_*$-values, the covariant derivatives of the coordinate basis vectors are linear combinations of coordinate basis vectors with bounded smooth coefficients. We are working in $L^2$ where we already control all the components. Thus, we control
the difference between components of the covariant derivative in a direction and the
covariant derivative of the components of the Maxwell tensor.

\subsection*{Summary of Sections}

The main aim of this paper is to prove decay results for the Maxwell field on the Reissner-Nordstr{\o}m-de Sitter Black Hole. We prove two types of decay: The first is a decay of the energy of the Maxwell field on achronal hypersurfaces in the static region as the hypersurfaces approach timelike infinity, with quadratic decay rate. This is Theorem \ref{decayonhypersurfacesTHEOREM}. The second decay result is Theorem \ref{PointwisedecayTHEOREM}. It is a pointwise decay in time with a rate of $t^{-1}$, also in the static region of the spacetime. Both results are consequences of the bounds on the conformal energy obtained from the wave analysis on the middle spin component of the field where we follow the work in \cite{blue_decay_2008}. 

The rest of the paper and the above general outline of the work is detailed in the sections of the paper as follows:

\paragraph{Section \ref{sec:geometricframework}:} This section contains an overview on the Reissner-Nordstr{\o}m-de Sitter black holes discussing the photon sphere which is the most relevant part of the geometry in this study of decay. We also fix some notations regarding the Maxwell system, in addition to some preliminary results we need.

\paragraph{Section \ref{sectionenergyestimates}:} This section is devoted to the analysis of the wave equation satisfied by the middle spin component of the field. We show that the energy for the wave equation is conserved and derive estimates for the conformal charge. Following \cite{blue_decay_2008}, we use these estimates and a Morawetz estimate using a radial multiplier to obtain a uniform bound on the conformal charge where the Hardy estimates are needed. We note that this is where the exclusion of stationary solutions becomes necessary so that we can control the $L^2$-norm of the wave solution by the norm of its angular derivatives. The uniform bound we get controls the integral of the trapping term multiplied by the angular derivative of the wave solution.

\paragraph{Section \ref{Decaysection}:} The fourth section of the paper is the decay results. We introduce some norms on the Maxwell 2-form and discuss the energies of the field. The stress-energy tensor is used to define the energies on a hypersurface, we then write them for the Cauchy hypersurface $\{t=0\}$ in terms of the spin components. We get an almost conservation law, describing quantitatively the influence of the trapping effect on the conformal energy defined by the Morawetz vector field $K$, where the significance of the photon sphere is manifested. We then relate the wave energy of the middle component and that of the full field, bounding the former by the energy and the conformal energy of derivatives of the latter. Using these results a uniform bound on the conformal energy is obtained. At this point we state and prove the decay results in section \ref{sec:decayresults}. Finally, in section \ref{genericsphericalf} we specify under what conditions this work and these decay results can be extended to other spacetimes, which include a wide range of spherically symmetric spacetimes.

\section{Geometric Framework}\label{sec:geometricframework}

We start by setting up with some properties of the Reissner-Nordstr{\o}m-de Sitter metric and the Maxwell field.

\subsection*{Reissner-Nordstr{\o}m-de Sitter Spacetime}

One of the spherically symmetric solutions of Einstein-Maxwell Field equations in the presence of a positive cosmological constant is the Reissner-Nordstr{\o}m-de Sitter solution (RNDS). It models a non-rotating spherically symmetric charged black hole with mass and a charge, in a de Sitter background. The de Sitter background means that there is a cosmological horizon beyond which lies a dynamic region that stretches to infinity, while the Reissner-Nordstr{\o}m nature entails that near the singularity, depending on the relation between the mass and the charge, one has a succession of static and dynamic regions separated by horizons. We shall recall some of its properties that are relevant for our purposes in this paper. For more on this spacetime we refer to \cite{mokdad_reissner-nordstrom-sitter_2017}.

The Reissner-Nordstr{\o}m-de Sitter metric is given in spherical coordinates by 
\begin{equation}\label{RNdSmetric}
g_\mathcal{M}=f(r)\d t^2-\frac{1}{f(r)}\texttt{d}r^2-r^2\d \omega^2,
\end{equation}
where
\begin{equation}\label{f(r)}
f(r)=1-\frac{2M}{r}+\frac{Q^2}{r^2}-\Lambda r^2 \; ,
\end{equation}
and $\d \omega^2$ is the Euclidean metric on the $2$-Sphere, $\mathcal{S}^2$, which in spherical coordinates is,
\begin{equation*}
\d \omega^2=\d \theta^2 + \sin(\theta)^2\d \varphi^2 \; ,
\end{equation*}
and $g_\mathcal{M}$ is defined on \index{$\mathcal{M}$}$\mathcal{M}=\R_t \times ]0,+\infty[_r \times \mathcal{S}^2_{\theta,\varphi}$ . Here $M$ is the mass of the black hole, $Q$ is its charge, and $\Lambda$ is the cosmological constant. We assume that $Q$ is real and non zero, and $M$ and $\Lambda$ are positive.

The metric in these coordinates appear to have singularities at $r=0$ and at the zeros of $f$. Only the singularity at $r=0$ is a real geometric singularity at which the curvature blows up. The apparent singularities at the zeros of $f$ are artificial and due to this particular choice of coordinates. The regions of spacetime where $f$ vanishes are essential features of the geometry of the black hole, they are the event horizons or {horizons}\index{Horizon} for short, and $f$ is called the {horizon function}\index{Horizon function}. If $f$ has three positive zeros and one negative, then the zeros in the positive range corresponds respectively, in an increasing order, to the {Cauchy horizon} or {inner horizon}, the {horizon of the black hole}  or the {outer horizon}, and the {cosmological horizon}. In this case, $f$ changes sign at each horizon and one has static and dynamic regions separated by these horizons.

In this work, we are interested in the decay in time of test (decoupled) Maxwell fields in the static region between the horizon of the black hole and the cosmological horizon, which we refer to as the {exterior static region}\index{Exterior static region}. This part of spacetime contains a {photon sphere}\index{Photon sphere}, i.e. null geodesics orbiting the black hole at fixed $r$. This is a priori an obstacle for the decay but we shall see that the field decays in spite of the existence of a photon sphere. 

Under the the following conditions,
\begin{equation}\label{GC}
Q\neq0 \quad \textrm{and} \quad 0<\Lambda < \frac{1}{12Q^2} \quad \textrm{and} \quad M_1 < M < M_2 \; ,
\end{equation}
where 
\begin{gather}\label{naming}
R=\frac{1}{\sqrt{6\Lambda}} \quad;\quad \Delta=1-12Q^2\Lambda \quad ; \quad m_1=R\sqrt{1-\sqrt{\Delta}} \quad ;\quad m_2=R\sqrt{1+\sqrt{\Delta}}\\ M_1=m_1-2\Lambda m_1^3 \quad ; \quad M_2=m_2-2\Lambda m_2^3 \label{naming2}\; , 
\end{gather}
we have
\begin{prop1}[Three Positive Zeros and One Photon Sphere]\label{1photonsphere}
	The function $f$ has exactly three positive distinct zeros if and only if (\ref{GC}) holds. In this case, there is exactly one photon sphere in the static exterior region of the black hole defined by the portion between the largest two zeros.
\end{prop1}
\begin{proof}
This is proved in \cite{mokdad_reissner-nordstrom-sitter_2017}.
\end{proof}

With the assumption of (\ref{GC}), let the zeros of $f$ be $r_0<0<r_1<r_2<r_3$. For $r>0$, we define the Regge-Wheeler coordinate function $r_*$ by requiring 
$$\frac{\d r_*}{\d r}=\frac{1}{f(r)}>0.$$
The Regge-Wheeler radial coordinate have the following expression:
\begin{equation*}
r_*(r)=\sum_{i=0}^3 a_i \ln |r-r_i| + a \quad ; \; a_i=-\frac{r_i^2}{\Lambda}\prod_{j\ne i}\frac{1}{(r_i-r_j)} \quad ; \; a=-\sum_{i=0}^3 a_i \ln |P_2-r_i|
\end{equation*}
where $$\left\{ r=P_2=\frac{3M+\sqrt{9M^2-8Q^2}}{2} \right\}$$ is the photon sphere hypersurface.

We now introduce the chart $(t,r_*,\theta,\varphi)$ over the exterior static region  $\mathcal{N}=\R_t \times]r_2,r_3[\times\mathcal{S}^2_{\omega}$. We see that $r_*$ is a strictly increasing continuous function of $r$ (thus a bijection) over the interval $]r_2,r_3[$, and ranges from $-\infty$ to $+\infty$. We also have $\partial_{r_*}=f\partial_r$ and $\d r = f \d r_*$. The RNdS metric in these coordinates is:
\begin{equation}\label{Reggewheelermetric}
g_\mathcal{N}=f(r)(\d t^2-\d r_*^2)-r^2\d \omega^2 \- .
\end{equation}

It will be useful for us in calculations to have the Christoffel symbols in the coordinates $(t,r_*,\theta,\varphi)=(\tilde{x}^0,\tilde{x}^1,\tilde{x}^2,\tilde{x}^3)$. The non zero symbols are:
\begin{gather}
\tilde{\Gamma}^0_{01}=\tilde{\Gamma}^1_{00}=\tilde{\Gamma}^1_{11}=\frac{f'}{2}~~~ ;~~~ \tilde{\Gamma}^1_{22}=-r ~~~ ;~~~\tilde{\Gamma}^1_{33}=-r\sin(\theta)^2  \nonumber       \\
\tilde{\Gamma}^3_{13}=\tilde{\Gamma}^2_{12}=\frac{f}{r}~~~ ;~~~\tilde{\Gamma}^2_{33}=-\cos(\theta)\sin(\theta)~~~ ;~~~\tilde{\Gamma}^3_{23}=\cot(\theta).\label{Crstflsymbinr*}
\end{gather}

	\begin{figure}
	\centering
	\includegraphics[scale=0.7]{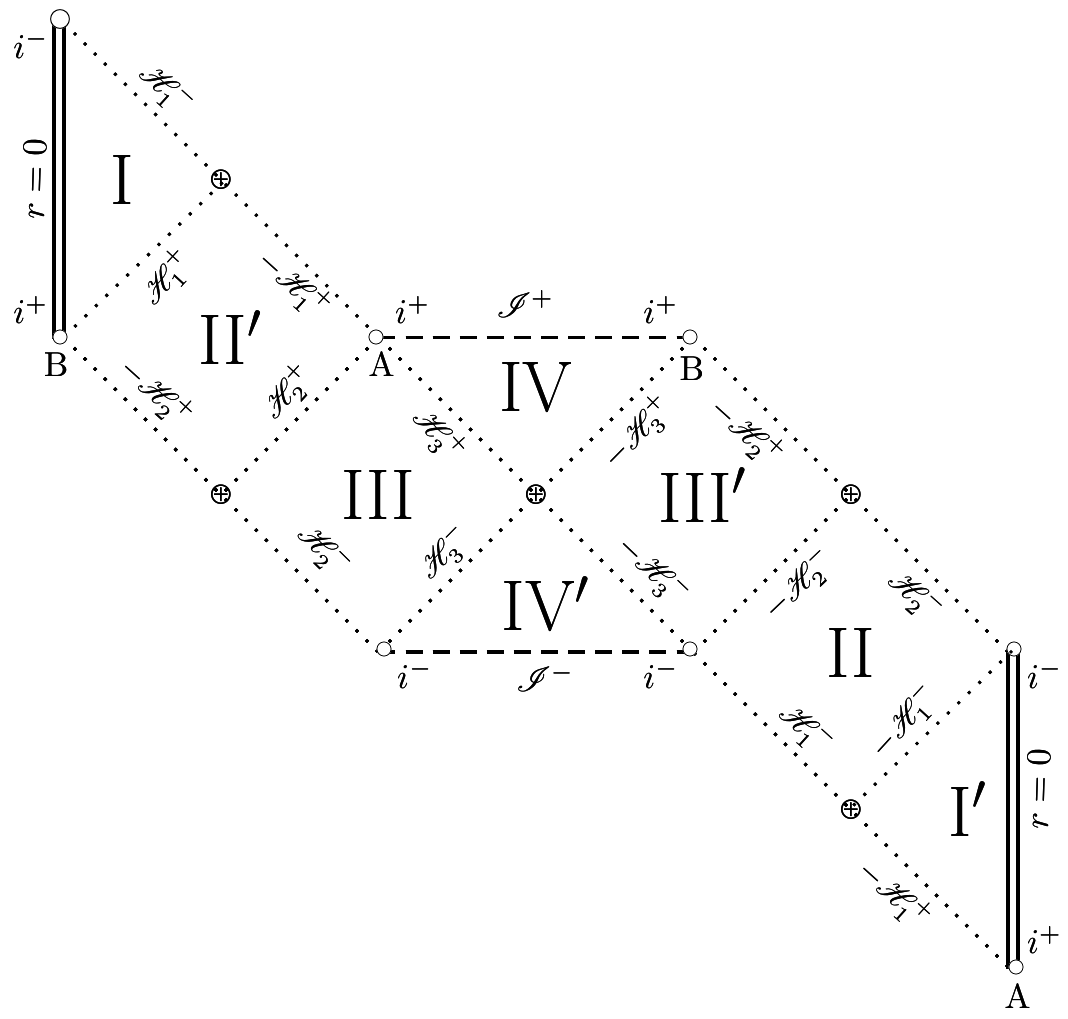}
	\caption{\emph{A building block of the Penrose-Carter conformal diagram of the maximal analytic extension of the RNDS spacetime in the case of three horizons. The maximal extension can be reconstructed from this building block by aligning infinitely many copies of the building block. The alignment is done without rotation and by overlapping $i^+$ points with the same label, A or B. The null hypersufaces $\pm\mathscr{H}_i^\pm$ are the horizons at $r_i$, and $\scri^\pm$ correspond to $r=\infty$.}}
	\label{fig:buildingblock}
\end{figure} 

\subsection*{Maxwell's Equations}

Let $F$ be a 2-form on the RNdS manifold $\mathcal{N}$. As we saw, the source free Maxwell's equations can be written as
\begin{eqnarray}
\delta F&=&0 \; , \label{Maxeq1}\\
\d F&=& 0 \; ,\label{Maxeq2}
\end{eqnarray}
where $\delta=\star \d \star$, and $\star$ is the Hodge star operator, or in abstract index notation,
\begin{eqnarray}
\nabla^a F_{ab} &=& 0 \; , \label{Maxeqabst1}\\
\nabla_{[ a} F_{bc ] }&=& 0 \; ,\label{Maxeqabst2}
\end{eqnarray}
and in coordinate form these translate to the following two sets of equations,
\begin{eqnarray}
g^{ac}\left(\partial_c F_{ab}-F_{db}\Gamma^d_{ca}-F_{ad}\Gamma^d_{cb}\right) &=& 0 \qquad \qquad\forall b,\label{Maxeqcoor1} \\
\partial_c F_{ab}+\partial_b F_{ca}+\partial_a F_{bc}&=& 0 \qquad \qquad \forall a,b,c. \label{Maxeqcoor2}
\end{eqnarray}
If taken in the coordinates $(t,r_*,\theta,\varphi)=(\tilde{x}^0,\tilde{x}^1,\tilde{x}^2,\tilde{x}^3)$, (\ref{Maxeqcoor1}) becomes respectively for $b=0,..,3$:
\begin{eqnarray}
&&\dl_1F_{10}+V\dl_2F_{20}+V\sin(\theta)^{-2}\dl_3F_{30}+2 rVF_{10}+VF_{20}\cot(\theta)+f'F_{01}=0, \label{Maxeqcoor10}\\
&&\dl_0F_{10}+V\dl_2F_{21}+V\sin(\theta)^{-2}\dl_3F_{31}+VF_{21}\cot(\theta)=0, \label{Maxeqcoor11} \\
&&\dl_1F_{12} + \dl_0F_{20}+V\sin(\theta)^{-2}\dl_3F_{32}=0, \label{Maxeqcoor12}\\
&&\dl_1F_{13} + \dl_0F_{30}+V\dl_2F_{23}+VF_{32}\cot(\theta)=0.\label{Maxeqcoor13}
\end{eqnarray}
where $V=fr^{-2}$\index{$V$}.

As much as equations (\ref{Maxeq1})-(\ref{Maxeqabst2}) are elegant and simple they are not the most convenient form for us to use in all arguments and calculations, and evidently neither are their expressions in coordinates. We shall use expressions that depends more on the geometry of the spacetime. This is the tetrad formalism.



Instead of working with the components of the Maxwell field in a coordinate basis, it is more convenient to use the components of the field in a general basis of the tangent space which might not be the canonical basis given by the coordinates. At each point, one defines a set of four vectors, called the tetrad, that forms a basis for the tangent space at that point. One  can then reformulate the field equations using this tetrad. In general relativity, it is natural to project on a null tetrad, which consists of two real null vectors and two conjugate null complex vectors usually defined as $X\pm iY$ for $X$ and $Y$ two spacelike real vectors. Here, we use a null tetrad on $\mathcal{N}$ given by two null real vectors and a two conjugate null complex vector tangent to the 2-Sphere $\mathcal{S}^2$:
\begin{eqnarray}\label{nulltetrad}
L &=& \dl_t +\dl_{r_*} \nonumber\\
N &=& \dl_t -\dl_{r_*} \nonumber\\
M &=& \dl_\theta +\frac{i}{\sin(\theta)}\dl_{\varphi} \\
\bar{M} &=&  \dl_\theta -\frac{i}{\sin(\theta)}\dl_{\varphi}\nonumber
\end{eqnarray}

We shall call this tetrad the ``{stationary tetrad}''. Using this tetrad, we can represent the Maxwell field by three complex scalar functions called the {spin components} of the Maxwell field associated to the given tetrad, and defined by:
\begin{eqnarray}\label{spincomponentsPhis}
\Phi_1&=&F\left(L,M\right) \nonumber\\
\Phi_0 &=& \frac{1}{2}\left(V^{-1}F(L,N)+F\left(\bar{M},M\right)\right) \\
\Phi_{-1}&=&F\left(N,\bar{M}\right) \nonumber
\end{eqnarray}

The conventional definition of spin components of an anti-symmetric tensor is slightly different. One normally defines it without the extra factor of $V^{-1}$ in the middle component $\Phi_0$. However the way we project the field on the tetrad is more convenient for studying decay in the region $\mathcal{N}$. Our conventions are the same as P. Blue's \cite{blue_decay_2008}\footnote{The conformal weight and the spin weight are respectively related to the way the component change when we rescale the complex vector of the tetrad by a complex constant and the conjugate vector by the conjugate constant, and when rescaling the first null vector of the tetrad by a real constant and the second by the inverse constant. More precisely , the components transform as powers of the real rescaling constant, the power being the index of the component.\label{footnotespinweight}}. Also the usual way to label the components is different, conventionally, they are indexed by $0$, $1$, and $2$. For more on these notations see \cite{penrose_spinors_1987,penrose_spinors_1988}. 

Note that the tetrad we use, unlike those in the {Newman-Penrose formalism}\index{Newman-Penrose formalism}, are not normalized: A {normalized tetrad}\index{Normalized tetrad} is such that the inner product of the two null real vectors of the tetrad with each other equals $1$, and the product of the null complex  vector with its conjugate is $-1$, while all other products are zero. The formalism we use is a form of Geroch–Held–Penrose formalism (GHP), which does not require normalization. The form of Maxwell's equations in this formalism is usually referred to as Maxwell's {compacted equations}\index{Compacted equations} (see \cite{penrose_spinors_1987}).

A straight forward coordinate calculation shows that in this framework, Maxwell's equations translate as follows.

\begin{lem1}[Maxwell Compacted Equations]\label{NewmanPenrose} $F$ satisfies Maxwell's equations (\ref{Maxeq1}) and (\ref{Maxeq2}) if and only if its spin components $(\Phi_1,\Phi_0,\Phi_{-1})$ in the stationary tetrad satisfy the compacted equations 
	\begin{eqnarray}
	N\Phi_1 &=& V M\Phi_0, \label{NewmanPenrose1} \\
	L\Phi_0 &=& \bar{M}_1\Phi_1, \label{NewmanPenrose2}\\
	N\Phi_0 &=& -M_1\Phi_{-1},  \label{NewmanPenrose3}\\
	L\Phi_{-1} &=& -V\bar{M}\Phi_0. \label{NewmanPenrose4}
	\end{eqnarray}
	where $M_1=M+\cot(\theta)$ and $\bar{M}_1$ is its conjugate.
\end{lem1}

We define the energy flux of a Maxwell field accross the hypersurface $\Sigma_t=\{t\}\times\R\times\mathcal{S}^2$ of constant $t$ by
\begin{equation*}
E_T[F](t)=\frac{1}{4}\int_{\Sigma_t} |\Phi_1|^2 +2V|\Phi_0|^2 + |\Phi_{-1}|^2 ~\d r_* \d^2 \omega \; .  
\end{equation*}
This norm is the natural energy associated with Maxwell's equations and it can be defined geometrically (see (\ref{enegrysurfacegeneralform}) and (\ref{energyofmaxwellont=cst}) of section \ref{Energysection}).

Evidently not all solutions of Maxwell's equations decay in time, take for example the case where $\bm{\Phi}$ is a non zero constant vector, then it satisfies (\ref{NewmanPenrose}) and clearly does not decay as it does not change with time. Even solutions having finite energy, may not decay in time: Consider the constant vector $\bm{\Phi}=(\Phi_1=0~,~{\Phi}_0=C\neq0~,~{\Phi}_{-1}=0)$, it has finite energy, yet it does not decay. Since Maxwell's equations are linear, the last example shows that solutions (even with finite energy) having charge do not decay, where by the charge of a Maxwell field we mean the constant $\Psi^0_{00}$, i.e. solutions having non zero $l=0$ part in the spin-wieghted sperical harmonic decomposition. So, we need to exclude such solutions in order to prove decay. 

More generally, {time-periodic solutions}\index{Time-periodic solution}, also called {stationary solutions}\index{Stationary solutions}, do not decay. These are solutions of the form $\bm{\phi}(r_*,\omega)e^{i\lambda t}$ ($\lambda$ is real), so they are solutions to $\dl_t \bm{\Phi}=i\lambda\bm{\Phi}$. Such solutions of the Maxwell's equations should be excluded in order to prove decay.
%
Adding the requirement of finite energy assumption:
\begin{equation}\label{Finiteenergy}
\Phi_1 ,\frac{\sqrt{f}}{r}\Phi_0, \Phi_{-1} \in L^2(\Sigma).
\end{equation} 
it is then known in the literature\footnote{The detials can be found in \cite{mokdad_maxwell_2016} for instance.}, that the only admissible time-periodic solutions are exactly the pure charge solutions, that is: 

\begin{prop1}[Stationary Solutions]\label{stationarysolutionstheorem}
	If $\mathbf{\Phi}=(\Phi_1,\Phi_0,\Phi_{-1})$ is a finite energy stationary solution of Maxwell's equations, then
	\begin{equation}
	\mathbf{\Phi}=\left(
	\begin{array}{c}
	0 \\
	C \\
	0 \\
	\end{array}
	\right)                \qquad \textit{where C is a complex constant}.
	\end{equation}
\end{prop1}

The solutions we will consider from now on are finite energy solutions with no stationary part, and by Proposition \ref{stationarysolutionstheorem}, these are the finite energy solutions in the orthogonal complement of the  $l=0$ subspace, that is solutions of the form:

\begin{eqnarray}
\Phi_{\pm 1}(t,r_*,\theta,\varphi)&=& \sum\limits_{l=1}^{+\infty}\sum\limits_{n=-l}^l\Psi_{\pm 1n}^l(t,r_*) W^l_{\pm 1n}(\theta,\varphi) \qquad\qquad \Psi_{\pm1n}^l  \in L^2(\R_{r_*})\; , \label{nonstationarysolutions1}\\
\Phi_0(t,r_*,\theta,\varphi)&=& \sum\limits_{l=1}^{+\infty}\sum\limits_{n=-l}^l\Psi_{0n}^l(t,r_*) W^l_{0n}(\theta,\varphi) \qquad\qquad \frac{\sqrt{f}}{r}\Psi_{0n}^l \in L^2(\R_{r_*})\; . \label{nonstationarysolutions}
\end{eqnarray}

where $$\{W^l_{mn}(\theta,\varphi) ;l,m,n \in \mathbb{Z}; l\geq0, -l\leq m,n\leq l \}$$ form an orthonormal basis of {spin-weighted spherical harmonics} of $L^2(\mathcal{S}^2)$. For simplicity, and by density, it is enough to consider only smooth compactly supported solution of the above form. 

If $F$ is a Maxwell field on $\mathcal{N}$ with spin components $\bm\Phi$ whose spin-weighted harmonic coefficients are $\bm{\Psi}^l_n$, then $F$ has a global potential if and only if the imaginary part of $\Psi^0_{00}$ vanish. Thus, as a consequence of the form (\ref{nonstationarysolutions1}) and (\ref{nonstationarysolutions}), the solutions we consider here have global potentials.

\subsection*{Divergence Theorem}

One important tool that we use frequently is the divergence theorem. We present a version of this theorem which we think is better suited for Lorentzian manifolds than the usual one found in textbooks on Riemannian geometry. 
The divergence of $X$ is defined to be the unique function $divX$ such that,
\begin{equation}\label{divergenceinternsicexpression}
\lie_X\omega=(divX)\omega \; .
\end{equation}
If the orientation on $\mathcal{U}$ is given by a pseudo-Riemannian metric $g$, i.e. $\omega=\d V_g$, then the above definition of $divX$ coincides with the more familiar one, which is locally defined as:
\begin{equation}\label{divergencelocalexp1}
\frac{1}{\sqrt{|g|}}\dl_i \left(\sqrt{|g|}X^i\right)\; ,
\end{equation}
where $|g|$ is the absolute value of the determinant of the metric $g$.

\begin{lem1}[Divergence Theorem]\label{divergencetheorem}
	Let $\mathcal{U}$ be an oriented smooth n-manifold with boundary (possibly empty), with $\omega$ a positively oriented volume form , i.e. determining the orientation on $\mathcal{U}$, and the boundary $\dl \mathcal{U}$ is outward oriented (Stokes' orientation), and let $X$ be a smooth vector field on $\mathcal{U}$. If $\mathcal{U}$ is compact or $X$ is compactly supported then,
	
	\begin{equation}\label{divergencetheoremgeneralformula}
	\int_{\dl \mathcal{U}} i_X \omega = \int_{\mathcal{U}}  \lie_X \omega= \int_{\mathcal{U}}  divX \omega \; ,
	\end{equation}
	
	Moreover, if the orientation on $\mathcal{U}$ is given by a pseudo-Riemannian metric $g$, i.e. $\omega=\d V_g$, then (\ref{divergencetheoremgeneralformula}) can be reformulated as:
	
	\begin{equation}\label{divergencetheoremnormal-transverseformula}
	\int_{\dl \mathcal{U}} N (X)  i_L \d V_g = \int_{\mathcal{U}}   \lie_X \d V_g = \int_{\mathcal{U}}   divX \d V_g\; ,
	\end{equation}
	where $N$ is a conormal field to $\dl \mathcal{U}$, i.e. $N^\sharp$ is a normal vector field, and $L$ is a vector field transverse to $\dl \mathcal{U}$, such that $N(L)=1$.
\end{lem1}

Note that if the normal vector field can be normalized, which is always the case if the metric is Riemannian, and is only true if the hypersurface is timelike in the Lorentzian case, one can then choose the transverse vector to be the normal itself and thus recovering the well known form of this theorem:
\begin{equation}\label{divergencetheoremwellknownformula}
\int_{\dl \mathcal{U}} N(X) \d V_{\tilde{g}} = \int_{\mathcal{U}} div X \d V_g \; , \hspace{1cm} \textrm{or,}  \hspace{1cm} \int_{\dl \mathcal{U}} N_a X^a \d V_{\tilde{g}} = \int_{\mathcal{U}} \nabla_a X^a \d V_g
\end{equation}
$\tilde{g}$ being the induced metric on $\dl \mathcal{U}$, and $\d V_{\tilde{g}}=i_{N^\sharp}\d V_g$.

Killing vector fields have the nice property of vanishing divergence. A vector field $X$ is said to be Killing if the metric is conserved along the flow of $X$, i.e. $\lie_X g=0$. 
Thus, $\lie_X g_{ab}=\nabla_a X_b + \nabla_b X_a= 2\nabla_{(a} X_{b)}$, and so for Killing fields,
\begin{equation}\label{killingequation}
\nabla_a X_b - \nabla_b X_a=0 \; ,
\end{equation}
consequently,
\begin{equation}\label{killingvanishdivergence}
0=g_{ab}\left(\nabla_a X_b - \nabla_b X_a\right)=2\nabla_a X^a \; ,
\end{equation}
hence $divX=0$. Equation (\ref{killingequation}) is called the Killing equation, and the $(0,2)$-tensor involved is sometimes called the deformation tensor or Killing tensor, denoted
\begin{equation}\label{deformationtensor}
{}^X \pi_{ab}= 2\nabla_{(a} X_{b)} \; .
\end{equation}

Sometimes it is useful to see more directly the dependence of the integral on the hypersurface as we did in (\ref{divergencetheoremnormal-transverseformula}) for a boundary hypersurface using a normal and a transverse vector field. The existence of these vector fields is a consequence of the fact that $S$ is an oriented hypersurface of a pseudo-Riemannian oriented manifold as the following lemma guarantees. 

\begin{lem1}
	If $S$ is a smooth orientable hypersurface of a smooth orientable $n$-manifold $M$, then $S$ admits a nowhere vanishing smooth 1-form $\alpha$ defined on a neighbourhood of $S$ in $M$ with the property that $\alpha_p \equiv 0$ on $T_p S$, for all $p\in S$. Such a 1-form is unique up to a multiplication by a smooth function that does not vanish.   
\end{lem1}

\section{Energy Estimates}\label{sectionenergyestimates}

In this section we derive several estimates that will help us prove decay. As we said in the introduction, we need to obtain a uniform bound on the middle component, which satisfies a wave-like equation, in order to control the energies of the Maxwell field. Accordingly, we shall first analyse the wave-like equation by proving different estimates on the energy and the conformal charge of the solutions. We then use Morawetz estimates to control the conformal charge and obtain the uniform bound.

\subsection{The Wave-Like Equation}\label{sectionwaveanalysis}

From Maxwell's compacted equations (\ref{NewmanPenrose1})-(\ref{NewmanPenrose4}) we can see that $\Phi_0$ satisfies the following wave equation:
\begin{equation}\label{WavePhi0}
\dl_t^2\Phi_0=\dl_{r_*}^2\Phi_0+ V \Delta_{\mathcal{S}^2}\Phi_0.
\end{equation}
One way this can be obtianed is by applying $L$ to (\ref{NewmanPenrose3}) then use (\ref{NewmanPenrose4}).

We use $\slashed{\Delta}$ as another symbol for the operator $\Delta_{\mathcal{S}^2}$. Similarly, $\slnb$ will designate the Levi-Civita connection $\nabla_{g_{\mathcal{S}^2}}$ on the sphere.

Since the box notation
\begin{equation}\label{waveequation}
\square u=\square_g u=\nabla^\alpha \nabla_\alpha u=g^{\mathbf{ab}}(\dl_\mathbf{a}\dl_\mathbf{b}u
- \Gamma^\mathbf{c}_\mathbf{ab}\dl_\mathbf{c}u)
\end{equation}
is reserved for the geometric wave equation on the RNdS manifold, we will use the \index{$\widetilde{\square}$} $\widetilde{\square}$ symbol for our wave
equation so that
\begin{equation}\label{boxtimesoperator}
\widetilde{\square} u= \dl_t^2 u -\dl_{r_*}^2 u- V \sldlt u
\end{equation}
and the equation $\widetilde{\square} u=0 $ is nothing but (\ref{WavePhi0}).

We also use the dot notation to designate the scalar product on $1$-forms on $\mathcal{S}^2$ as well as for the divergence of a $1$-form on $\mathcal{S}^2$. So we denote, for $\alpha$ and $\beta$, $1$-forms on $\mathcal{S}^2$,
\begin{equation}\label{dotnotation1}
\alpha \cdot \beta:=\alpha^a \beta_b=g^{-1}_{\mathcal{S}^2}(\alpha ,\beta)\qquad ; \qquad |\alpha|^2=g^{-1}_{\mathcal{S}^2}(\alpha ,\alpha) \qquad ; \qquad \  \slnb  \cdot \alpha :=\slnb^a \alpha_a
\end{equation}
We can readily see that for a smooth function $u$,
\begin{equation}\label{dotnotation2}
\slnb  \cdot (u\alpha) =\slnb u  \cdot \alpha + u\slnb  \cdot \alpha:= g^{-1}_{\mathcal{S}^2}(\slnb u, \alpha) + u\slnb  \cdot \alpha \; .
\end{equation}
and
\begin{equation}\label{dotnotation3}
\slnb \cdot \slnb u = \slnb^a  \slnb_a u= \sldlt u.
\end{equation}

Finally, throughout this section, $C$ will designate a constant that may change from a line to another and which depends \textit{neither} on $t,r_*,\omega$, \textit{nor} on the solution $u$.

\subsubsection{Energy and Conformal Charge}

For solutions of $\widetilde{\square} u =0$, there are two important quantities, one of which is conserved and the other is controlled. They are related to the time-translation vector field $T=\dl_t$ and Morawetz vector field $K=(t^2+r_*^2)\dl_t + 2tr_*\dl_{r_*} u$. These are the associated energy and the conformal charge, and they are given by the integral of their respective densities:

\begin{eqnarray}
e&=&(\dl_t u)^2 + (\dl_{r_*}u)^2 + V|\slnb u|^2 \; , \label{energydensity}\\
e_\mathcal{C}&=&\frac{1}{2} (t^2+r_*^2)e + 2tr_*\dl_t u\dl_{r_*}u + e \; .\label{conformaldensity}\\
E[u](t)&=&\frac{1}{2}\int_{\Sigma_t} e \d r_* \d^2 w \; , \label{energywavelike}\\
E_{\mathcal{C}}[u](t)&=&\frac{1}{2}\int_{\Sigma_t} e_{\mathcal{C}} \d r_* \d^2 w \; , \label{conformalchargewave}
\end{eqnarray}
where $\Sigma_t = {\{t\}\times \R \times \mathcal{S}^2}$.\index{$E_{\mathcal{C}}[u](t)$}\index{$E[u](t)$}

These densities are positive quantities noting that the conformal charge density $e_\mathcal{C}$ can be written as the sum of squares,

\begin{equation}\label{conformaldensitypositive}
e_\mathcal{C}=\frac{1}{4}[(t+r_*)(\dl_t+\dl_{r_*})u)]^2 + \frac{1}{4}[(t-r_*)(\dl_t-\dl_{r_*})u)]^2 + \frac{1}{2}(t^2+r_*^2)V|\slnb u|^2 + e \; .
\end{equation}

The conservation laws can be obtained using a multiplier method, which is what we do in Lemma  \ref{controlonconformalchargeTHEOREM} and its proof, nonetheless, we obtain the conservation law first using a geometrical approach in which we use the Lagrangian method on the wave equation $\square u =0$, thus relating $T$ and $E[u](t)$ is a geometric way.

The natural energy associated to the wave equation $\square u =0$ is generated by a stress-energy tensor contracted with the Killing vector field $T=\dl_t$ which describes a static observer at infinity. The stress-energy tensor for the wave equation $\square \phi=0$ in the abstract index notation is given by

\begin{equation*}
\mathbf{S}_{ab}=\nabla_a \phi \nabla_b \phi -\frac{1}{2} g_{ab}\nabla^c \phi \nabla_c \phi \; .
\end{equation*}

For solutions of the wave equation, the stress-energy is divergence-free. More generally,

\begin{equation}\label{Divergenceofstressenergy}
\nabla^a \mathbf{S}_{ab}=\nabla_b \phi \square \phi \; .
\end{equation}

\begin{lem1}\label{conservationofenergy}
	The energy of solutions to $\widetilde{\square} u=0$ is conserved, i.e.
	\begin{equation*}
	E[u](t)=E[u](0) \quad \forall t \in\R \; .
	\end{equation*}
\end{lem1}

\begin{proof}
	First we show that
	\begin{equation}\label{wave-wavelikeeq}
	\widetilde{\square} u = rf\square \left(\frac{u}{r}\right) - \frac{ff'}{r} u \; .
	\end{equation}
	Calculating directly from (\ref{waveequation}), we have
	\begin{eqnarray*}
		\square u &=& \frac{1}{f}\left(\dl^2_t u - \dlr^2 u \right) -\frac{2}{r}\dlr u  -\frac{1}{r^2}\left(\dl_\theta^2 u + \frac{1}{\sin(\theta)^2}\dl_\varphi^2 u +\cot(\theta)\dl_\theta u  \right)  \\
		&=& \frac{1}{f}\left(\dl^2_t u - \dlr^2 u \right) -\frac{2}{r}\dlr u  -\frac{1}{r^2}\sldlt u \; .
	\end{eqnarray*}
	The only terms that still need to be computed are those with partial derivative with respect to $r_*$,
	\begin{eqnarray}
	\dlr\left(\frac{u}{r}\right) &=& \frac{1}{r}\dlr u - \frac{f}{r^2} u \; , \label{dlr u over r} \\
	\dlr^2\left(\frac{u}{r}\right) &=& \frac{1}{r}\dlr^2 u - 2\frac{f}{r^2} \dlr u + u\frac{f}{r^3}(2f -r f') \; .  \label{dlr2 u over r}
	\end{eqnarray}
	Putting these terms in the above expression of $\square u$ and using (\ref{boxtimesoperator}) we get (\ref{wave-wavelikeeq}).
	
	The natural energy associated to the wave equation is:
	\begin{equation}\label{Energywaveeq}
	\tilde{E}[\phi](t)=\int_{\Sigma_t} \mathbf{S}_{ab} \hat{T}^a T^b \d \sigma
	\end{equation}
	where $\d \sigma$ is the induced measure defined on $\Sigma_t$.
	If $\phi$ is a solution of the wave equation then the energy defined above is conserved. In general, for any $\phi$
	\begin{equation}\label{Energydiffrence}
	\tilde{E}[\phi](t)-\tilde{E}[\phi](0)=  \!\!\!\! \int\limits_{[0,t] \times \R \times \mathcal{S}^2} (\nabla^a \mathbf{S}_{ab}) T^b \d^4 x \; .
	\end{equation}
	where $T=\dl_t$, $\hat{T}=f^{-1}T$, and $\d^4 x$ is the $4$-volume measure on the RNdS manifold. To see why (\ref{Energydiffrence}) is true, we start with the difference and apply the divergence theorem,
     noting that $\hat{T}^a$ is the unit normal vector field on $\Sigma_t$,
	
	\begin{eqnarray*}
		\tilde{E}[\phi](t)-\tilde{E}[\phi](0) &=&\int_{\Sigma_t} \mathbf{S}_{ab} T^b\hat{T}^a \d \sigma -
		\int_{\Sigma_0} \mathbf{S}_{ab} T^b\hat{T}^a \d \sigma \\
		&=& \!\!\!\! \int\limits_{[0,t] \times \R \times \mathcal{S}^2}\nabla^a (\mathbf{S}_{ab}T^b) \d^4 x \\
		&=& \!\!\!\! \int\limits_{[0,t] \times \R \times \mathcal{S}^2} \left[(\nabla^a \mathbf{S}_{ab})T^b +(\nabla^a T^b)\mathbf{S}_{ab}\right] \d^4 x\; ,
	\end{eqnarray*}
	but $\mathbf{S}_{ab}$ is symmetric, so
	$$(\nabla^a T^b)\mathbf{S}_{ab}=\frac{1}{2}(\nabla^a T^b + \nabla^b T^a )\mathbf{S}_{ab} \; ,$$
	and $T$ is Killing i.e.
	$$0=\mathcal{L}_T g^{ab}=\nabla^a T^b + \nabla^b T^a \; .$$
	
	On the other hand,
	
	\begin{equation*}
	\d \sigma = \hat{T} \lrcorner \d^4 x = \left(f^{-\hf}\dl_t\right) \lrcorner
	\left(f r^2 \d t \wedge \d r_* \wedge \d^2 \omega  \right)= f^{\hf} r^2 \d r_* \wedge \d^2 \omega \; ,
	\end{equation*}
	where $\d^2 \omega=\sin(\theta)\d \theta \wedge \d \varphi$ is the Euclidean area element on $\mathcal{S}^2$. And so,
	\begin{equation*}
	\tilde{E}[\phi](t)=\int_{\Sigma_t} \mathbf{S}_{00} r^2 \d r_* \d^2 \omega
	\end{equation*}
	with
	\begin{equation*}
	\mathbf{S}_{00}=\hf \left((\dl_t \phi)^2 + (\dl_{r_*} \phi)^2 + \frac{f}{r^2}|\slnb \phi|^2 \right) \; ,
	\end{equation*}
	
	\begin{equation*}
	\slnb \phi = \dl_\theta \phi \d \theta + \dl_\varphi \phi  \d \varphi \, , \qquad \mathrm{and} \qquad
	|\slnb \phi|^2=\slnb \phi \cdot \slnb \phi = \left(\dl_\theta \phi \right)^2 + \sin(\theta)^{-2}\left(\dl_\varphi \phi \right)^2 \; .
	\end{equation*}
	We can write
	\begin{equation*}
	\tilde{E}[\phi](t)=\hf \int_{\Sigma_t} \left(\left(\dl_t (r\phi  )\right)^2 + \left(r\dl_{r_*} \phi\right)^2 +
	\frac{f}{r^2}\left|\slnb (r\phi  )\right|^2 \right) \d r_* \d^2 \omega \; ,
	\end{equation*}
	and 
	we can apply a double integration by parts on the middle term of the energy, and making use of the previous calculations in (\ref{dlr u over r}) and (\ref{dlr2 u over r}) we get
	\begin{eqnarray*}
		\int_{\Sigma_t} (r\dl_{r_*} \phi)^2  \d r_* \d^2 \omega &=&-\int_{\Sigma_t}\phi \dl_{r_*}(r^2\dl_{r_*} \phi)
		\d r_* \d^2 \omega\\
		&=&-\int_{\Sigma_t}r \phi \left(2f\dlr\left(\frac{1}{r}r\phi \right) + r \dlr^2 \left(\frac{1}{r}r\phi \right)\right)\d r_*\d^2 \omega\\
		&=&-\int_{\Sigma_t}r \phi \left(-2\frac{f^2}{r^2}r\phi -r\left(\frac{ff'}{r^2}r\phi-2\frac{f^2}{r^3}r\phi\right) +\dlr^2 (r\phi)\right)\d r_*\d^2 \omega\\
		&=&-\int_{\Sigma_t}r\phi  \dl^2_{r_*}( r\phi)\d r_* \d^2 \omega + \int_{\Sigma_t}(r\phi  )^2 \frac{ff'}{r}
		\d r_* \d^2 \omega\\
		&=&\int_{\Sigma_t} \left(\dl_{r_*} (r\phi)\right)^2 \d r_* \d^2 \omega+ \int_{\Sigma_t}(r\phi)^2 \frac{ff'}{r}
		\d r_* \d^2 \omega \; .
	\end{eqnarray*}
	Therefore, if we set $u=r\phi$ then,
	\begin{equation}\label{Energyphiu}
	\tilde{E}[\phi](t)= E[u](t) + \hf \int_{\Sigma_t} \frac{ff'}{r}u^2\d r_* \d^2 \omega \; .
	\end{equation}
	By (\ref{wave-wavelikeeq}), $u$ is a solution of (\ref{WavePhi0}) if and only if  $\phi$ satisfies
	\begin{equation*}
	\square \phi = \frac{f'}{r} \phi \; ,
	\end{equation*}
	It follows from (\ref{Divergenceofstressenergy}) and (\ref{Energydiffrence}) that
	
	\begin{eqnarray*}
		\tilde{E}[\phi](t)-\tilde{E}[\phi](0) &=&  \!\!\!\! \int\limits_{[0,t] \times \R \times \mathcal{S}^2}
		(\nabla^a \mathbf{S}_{ab}) T^b \d^4 x \\
		&=& \!\!\!\! \int\limits_{[0,t] \times \R \times \mathcal{S}^2}\square \phi (\nabla_b \phi) T^b \d^4 x  \\
		&=& \!\!\!\! \int\limits_{[0,t] \times \R \times \mathcal{S}^2}\left(f'r^{-1} \phi\right) (\dl_t \phi) fr^2
		\d t \d r_* \d^2 \omega \\
		&=&\int\limits_{\R \times \mathcal{S}^2}r f f'\int\limits_{[0,t]}\phi \dl_t \phi \, \d t\quad
		\d r_* \d^2 \omega \\
		&=&\hf\int\limits_{\R \times \mathcal{S}^2}r f f' (\phi(t)^2-\phi(0)^2) \d r_* \d^2 \omega  \\
		&=&\hf \int_{\Sigma_t}r f f' (\phi^2)\d r_* \d^2 \omega -\hf \int_{\Sigma_0}r f f' (\phi^2)\d r_* \d^2 \omega \\
		&=&\hf \int_{\Sigma_t} \frac{f f'}{r} u^2\d r_* \d^2 \omega -\hf \int_{\Sigma_0}\frac{f f'}{r} u^2\d r_* \d^2 \omega \; .
	\end{eqnarray*}
	Evaluating the left hand side using (\ref{Energyphiu}) we obtain the conservation law.
\end{proof}

In order to control the conformal charge $E_{\mathcal{C}}[u]$ we use the \textit{Morawetz multiplier}\index{Morawetz multiplier} $K(u)$: developing $K(u)\widetilde{\square} u =0$ and then integrating by parts. In fact, the conformal charge is more associated with $T+K$ rather than just $K$. As we shall see, near the end of the proof of the next lemma, we use the fact $E[u]$ is conserved to obtain control on the error term of the conformal charge, which is essentially the same as using the multiplier $(T+K)(u)$ in place of $K(u)$.

\begin{lem1}\label{controlonconformalchargeTHEOREM}
	If $\widetilde{\square} u =0$, then there is a non-negative compactly supported smooth function of $r_*$ , $\chi_{trap}$, such that for all $t_1 , t_2>0$,
	\begin{equation}\label{controlonconformalcharge}
	E_\mathcal{C}[u](t_2)-E_\mathcal{C}[u](t_1)\le \int\limits_{[t_1,t_2] \times \R \times \mathcal{S}^2} t \chi_{trap}  |\slnb u|^2 \d t \d r_* \d^2 \omega \; .
	\end{equation}
\end{lem1}
\begin{proof}
	We develop $K(u)\widetilde{\square} u =0$. First, we have
	\begin{equation}\label{conservationofenergysecondpoof}
	\dl_t u \widetilde{\square} u = \hf \dl_t e - \dl_{r_*} (\dl_t u \dl_{r_*} u) - \slnb \cdot \left(V\dl_t u  \slnb u \right),
	\end{equation}
	since,
	\begin{eqnarray*}
		\dl_{t}(|\slnb u |^2)&=&2\left( \slnb u \cdot \slnb(\dl_{t} u)\right) \; ; \\
		\hf \dl_t e &=& \dl_t u \dl^2_t  u + \dl_{r_*} u \dl^2_{t r_* } u +V(\slnb u \cdot \slnb(\dl_t u)) \; ; \\
		\dl_{r_*} (\dl_t u \dl_{r_*} u) &=&  \dl_{r_*} u \dl^2_{t r_* } u + \dl_t u \dl^2_{r_*}  u \; ; \\
		\slnb \cdot \left(V\dl_t u  \slnb u \right) &=& V \left( \slnb u \cdot \slnb(\dl_t u)\right)+ V \dl_t u \slnb \cdot \slnb u \; .
	\end{eqnarray*}
	Similarly, 
	\begin{equation*}
	\dl_{r_*} u \widetilde{\square} u = -\hf \dl_{r_*} e + \dl_{t} (\dl_{r_*} u \dl_t u ) - \slnb \cdot \left(V\dl_{r_*} u  \slnb u \right) + \hf \dl_{r_*}\left(V\right)|\slnb u |^2 + V\dl_{r_*}(|\slnb u |^2).
	\end{equation*}
	
	Next we integrate  $K(u)\widetilde{\square} u =0$ over the domain $\Omega=[t_1,t_2] \times \R \times \mathcal{S}^2$ and use integration by parts. For simplicity, let,
	\begin{equation*}
	\d x = \d t \d r_* \d^2 \omega \; ,  \qquad \mathrm{and} \qquad \d \varsigma =  \d r_* \d^2 \omega \; .
	\end{equation*}
	We divide the integral into two parts,
	\begin{equation*}
	\int_\Omega K(u)\widetilde{\square} u \d x =  \int_\Omega (t^2 + r_*^2)\dl_t u \widetilde{\square} u \d x + \int_\Omega 2 t r_*\dl_{r_*} u \widetilde{\square} u \d x \; .
	\end{equation*}
	For the first term of the first part, by integration by parts in the $t$ variable, we have:
	\begin{eqnarray*}
		\hf \int_\Omega (t^2 + r_*^2)\dl_t e \d x 
		&=& \hf\left(\int_{\Sigma_{t_2}} (t^2+r_*^2)e \d \varsigma -\int_{\Sigma_{t_1}} (t^2+r_*^2)e \d \varsigma \right) -\int_\Omega t e \d x \; .
	\end{eqnarray*}
	We do a similar integration by parts for the second term but this time in the $r_*$ variable noting that $u$ is taken to be compactly supported in $r_*$.
	\begin{equation*}
	\int_\Omega (t^2 + r_*^2)\dl_{r_*} (\dl_t u  \dl_{r_*} u )\d x= -2 \int_\Omega r_* \dl_t u  \dl_{r_*} u \d x \; .
	\end{equation*}
	The last term in the first part of the integral is zero by the divergence theorem on $\mathcal{S}^2$,
	\begin{equation*}
	\int_\Omega (t^2 + r_*^2)\slnb \cdot \left(V\dl_t u  \slnb u \right) \d x=0\; .
	\end{equation*}
	
	Next, and using the same technique, we have:
	\begin{equation*}
	\hf\int_\Omega 2 t r_* \dl_{r_*}e \d x = -\int_\Omega t e \d x \; ;
	\end{equation*}
	\begin{equation*}
	2\int_\Omega t r_*\dl_{t} (\dl_{r_*} u \dl_t u ) \d x = 2\left(\int_{\Sigma_{t_2}}t r_* \dl_{r_*}u \dl_t u \d \varsigma -\int_{\Sigma_{t_1}} t r_* \dl_{r_*}u \dl_t u  \d \varsigma \right) - 2\int_\Omega r_* \dl_{r_*}u \dl_t u \d x \; ;
	\end{equation*}
	\begin{equation*}
	2\int_\Omega t r_*\slnb \cdot \left(V\dl_t u  \slnb u \right) \d x = 0\; .
	\end{equation*}
	Finally, an integration by parts on the last term yields,
	\begin{eqnarray*}
		\int_\Omega \left(   t r_*\dl_{r_*}\left(V\right)|\slnb u |^2 + 2 t r_*V\dl_{r_*}(|\slnb u |^2)\right) \d x \hspace{-0.3cm}&=&- \int_\Omega \left( t r_*\dl_{r_*}\left(V\right)|\slnb u |^2 + 2t V|\slnb u |^2\right)\d x \\
		\hspace{-0.3cm}&=&- \int_\Omega t \left( r_*\dl_{r_*}V + 2 V \right)|\slnb u |^2 \d x \; .
	\end{eqnarray*}
	
	Putting everything together, we get:
	
	\begin{eqnarray*}
		\int_\Omega K(u)\widetilde{\square} u &=& 0 \\
		&\Longleftrightarrow&  \\
		\hf \int_{\Sigma_{t_2}} \left(\hf (t^2+r_*^2)e + 2tr_*\dl_t u\dl_{r_*}u\right) \d \varsigma \!\!\!\! &-& \!\!\!\! \hf \int_{\Sigma_{t_1}} \left( \hf (t^2+r_*^2)e + 2tr_*\dl_t u\dl_{r_*}u \right)\d \varsigma   \\
		&& \qquad =\int_\Omega t \left( r_*\dl_{r_*}V + 2 V \right)|\slnb u |^2 \d x \; .
	\end{eqnarray*}
	By conservation of the energy $E[u](t)$,
	\begin{equation*}
	\int_{\Sigma_{t_2}} e \d \varsigma -  \int_{\Sigma_{t_1}} e \d \varsigma =0,
	\end{equation*}
	thus,
	\begin{equation*}
	E_\mathcal{C}[u](t_2)-E_\mathcal{C}[u](t_1) = \int_\Omega t \left(  r_*\dl_{r_*}V + 2 V \right)|\slnb u |^2 \d x= \int_\Omega 2tV\mathscr{T} |\slnb u |^2 \d x \; .
	\end{equation*}
	where
	\begin{equation}\label{trappingterm}
	\mathscr{T}= \frac{r_*}{2V}\dl_{r_*}V = 1+ r_*\left(\frac{f'}{2} - \frac{f}{r}\right) + 1 \; ,
	\end{equation}
	is called the {trapping term}. It is easy to see that as $r_*$ goes to minus infinity, $r$ approaches $r_2$ the middle zero of $f$, and $f'(r_2)>0$. Also, we see that as $r_*$ goes to plus infinity, $r$ approaches $r_3$ the largest zero of $f$, and $f'(r_3)<0$. This means that the limits of $\mathscr{T}$ are negative at both infinities, and so, $\mathscr{T}$ is positive only on a compact interval of $r_*$. Now since $V>0$, there exists some non-negative compactly supported function $\chi_{trap}$ of $r_*$ which dominates $2V\mathscr{T}$. This proves (\ref{controlonconformalcharge}).
\end{proof}

As we mentioned before, Lemma \ref{conservationofenergy} can also be proved using a multiplier method. From (\ref{conservationofenergysecondpoof}) we  see directly that
\begin{equation*}
\int\limits_{[0,t] \times \R \times \mathcal{S}^2} \dl_t u \widetilde{\square} u \d t \d r_* \d^2 \omega=E[u](t)-E[u](0)
\end{equation*}
which is the conservation law since the integrals of the last two terms of the right hand side of (\ref{conservationofenergysecondpoof}) are zero by the divergence theorem.

\subsection{Morawetz Estimate}

In order to obtain a uniform bound on the conformal charge, we use a weighted radial  Soffer-Morawetz multiplier that points away form the photon sphere at $r_*=0$ so its weight changes sign there. The error terms coming from the multiplier method can be controlled by energy localized inside the light cone,
\begin{equation}\label{localenergy}
E_{\ell}[u](t)=\int\limits_{\{t\}\times \{|r_*|\le \frac{3}{4}t\}\times \mathcal{S}^2} \!\!\!\!\!\!\!\!\!\!e \d r_* \d^2 \omega \; .
\end{equation}
Using this multiplier and suitable Hardy estimates, we control the error term of the conformal charge by the energy localized inside the light cone. This in turn is controlled by the conformal charge multiplied by $t^2$. This is because, inside the light cone, the conformal charge density controlled the energy density times a factor of $t^2$. This factor of $t^2$ compensate for the factor in the linear bound on the conformal charge, allowing us to obtain the uniform bound we need.

For the rest of this section, $u$ will be a smooth compactly supported function of the form (\ref{nonstationarysolutions}). Also, as before, let
\begin{equation*}
\d x = \d t \d r_* \d^2 \omega \; ,  \qquad \mathrm{and} \qquad \d \varsigma =  \d r_* \d^2 \omega \; .
\end{equation*}

We say that two functions $v$ and $w$ are equivalent over a set $A$ and write $v \sim w$, if there exists a positive constant $C>0$ such that for all $x\in A$,
\begin{equation*}
\frac{1}{C}v(x) \le w(x) \le C v(x) \; .
\end{equation*}

\subsubsection{Useful Estimates and Identities}

We will need a couple of estimates for the solutions of (\ref{WavePhi0}). Since we exclude the stationary solutions with finite energy of the Maxwell field equations, we can then benefit from the following estimates.
\begin{lem1}\label{angularboundlemma}
	If $u$ is of the form (\ref{nonstationarysolutions}), then
	\begin{equation}\label{angularbound}
	\int_{\mathcal{S}^2} u^2 \d^2 \omega \le \hf \int_{\mathcal{S}^2} |\slnb u|^2 \d^2 \omega \; .
	\end{equation}
\end{lem1}
\begin{proof}
This can be found in the literature. The reader may refer to \cite{mokdad_maxwell_2016}.

%
	
\end{proof}

We now establish some Hardy-like estimates.

\begin{lem1}\label{hardylemma}
	Let $t\ge 1$, $0<\sigma $, 
	and let $\xi$ be a non-negative function of $r_*$ which is positive in an open non-empty subinterval of $|r_*|\le \hf $, and $u$ be a smooth compactly supported function. Then,
	\begin{eqnarray}
	\int\limits_{\{t\}\times \{|r_*|\le \hf t\}\times \mathcal{S}^2} \frac{u^2}{(1+( r_*)^2)^{\sigma + 1}} \d\varsigma &\le& C \int\limits_{\{t\}\times \{|r_*|\le \hf t\}\times \mathcal{S}^2} \left(\frac{(\dl_{r_*} u)^2}{(1+( r_*)^2)^{\sigma}} + \xi u^2\right) \d\varsigma  \; . \label{hardy1} 
	\end{eqnarray}
	Moreover, if $u$ is given by (\ref{nonstationarysolutions}) then,
	\begin{eqnarray}
	\int\limits_{\{t\}\times \{|r_*|\le \frac{3}{4}t\}\times \mathcal{S}^2} \frac{u^2}{1+r^2_*} \d\varsigma &\le& C E_\ell[u](t)\; , \label{hardy2} 
	\end{eqnarray}
\end{lem1}
\begin{proof}
	For simplicity, we denote $u(t,s=r_*,\omega)$ by $u(s)$ since $t$ is given and fixed, and $\omega$ is irrelevant for this calculation. Let $\alpha \ge 0$ for now. For $s_1 > 0$ we have,
	\begin{eqnarray*}
		\frac{u(s_1)^2}{(1+ s_1)^{\alpha +1}} - u(0)^2 &=& \int_0^{s_1} \dl_s \left(\frac{u^2}{(1+ s)^{\alpha +1}}\right) \d s \\
		&=&\int_0^{s_1}\left(\frac{2u \dl_s u}{(1+ s)^{\alpha +1}}- \frac{(\alpha +1)u^2}{(1+ s)^{\alpha +2}}\right)\d s \; .
	\end{eqnarray*}
	But for $s\ge 0$,
	\begin{equation*}
	\frac{1}{(1+ s)^\alpha}\left(\frac{u}{1+ s}\sqrt{\frac{(\alpha +1)}{2}} - \dl_s u\sqrt{\frac{2}{(\alpha +1)}}\right)^2\ge 0 \; ,
	\end{equation*}
	hence,
	\begin{equation*}
	\left(\frac{\alpha + 1}{2}\right) \left(\frac{u^2}{(1+s)^{\alpha +2}}\right) + \left(\frac{2}{\alpha + 1}\right)\left(\frac{(\dl_s u)^2}{(1+s)^{\alpha}}\right)\ge \frac{2u \dl_s u}{(1+s)^{\alpha +1}}\; .
	\end{equation*}
	Thus,
	\begin{equation*}
	\frac{u(s_1)^2}{(1+s_1)^{\alpha +1}} - u(0)^2 \le -\frac{\alpha +1}{2}\int_0^{s_1} \frac{u^2}{(1+s)^{\alpha+2}} \d s + \frac{2}{\alpha +1}\int_0^{s_1} \frac{(\dl_s u)^2}{(1+s)^\alpha} \d s \; ,
	\end{equation*}
	and since
	\begin{equation*}
	\frac{u(s_1)^2}{(1+s_1)^{\alpha +1}} \ge 0 \; ,
	\end{equation*}
	we have,
	\begin{equation*}
	\int_0^{s_1} \frac{u^2}{(1+s)^{\alpha+2}} \d s \le \frac{4}{(\alpha +1)^2}\int_0^{s_1} \frac{(\dl_s u)^2}{(1+s)^\alpha} \d s + \frac{2}{\alpha +1}u(0)^2 \; .
	\end{equation*}
	Since the Max norm and the Euclidean norm are equivalent over $\R^2$, there are some $a,b>0$ such that for all $s\in \R$,
	\begin{equation*}
	a(1+|s|)\le \sqrt{1+s^2} \le b(1+|s|)\; ,
	\end{equation*}
	and so,
	\begin{equation}\label{normequivalenceMaxandEuclid}
	\frac{1}{C(\beta)}(1+|s|)^\beta \le (1+s^2)^{\beta/2} \le C(\beta)(1+|s|)^\beta \qquad\qquad \beta \ge 0 \; ,
	\end{equation}
	implying that,
	\begin{equation}\label{estimateinhardy1}
	\int_0^{s_1} \frac{u^2}{(1+s^2)^{1+\alpha/2}} \d s \le C\int_0^{s_1} \frac{(\dl_s u)^2}{(1+s^2)^{\alpha/2}} \d s + Cu(0)^2 \; .
	\end{equation}
	The same estimate holds true for the function $v$ defined by $v(s)=u(-s)$. Thus, for all $s_1,s_2>0$ and all $u$ smooth,
	\begin{equation*}
	\int_{-s_2}^{s_1} \frac{u^2}{(1+s^2)^{1+\alpha/2}} \d s \le C\int_{-s_2}^{s_1} \frac{(\dl_s u)^2}{(1+s^2)^{\alpha/2}} \d s + Cu(0)^2 \; .
	\end{equation*}
	In particular, this holds over $[-s_1-s_0,s_1 - s_0]$ with $|s_0|<s_1$, and for the function $v(s)=u(s+s_0)$, i.e.
	\begin{equation}\label{estimateinhardy2}
	\int_{-s_1}^{s_1} \frac{u^2}{(1+(s-s_0)^2)^{1+\alpha/2}} \d s \le C\int_{-s_1}^{s_1} \frac{(\dl_s u)^2}{(1+(s-s_0)^2)^{\alpha/2}} \d s + Cu(s_0)^2 \; .
	\end{equation}
	Since the function 
	$$\left(\frac{1+(s-s_0)}{1+s^2}\right)^\beta$$
	is positive and tends to 1 at both infinities, we have the following equivalence: For all $s,s_0\in \R$ and $|s_0|<a$ with $a>0$ and $\beta \ge 0$, there exists a constant $\tilde{C}(a,\beta)>0$, depending on $a$ and $\beta$ only, such that,
	\begin{equation*}
	\frac{1}{\tilde{C}(a,\beta)}(1+s^2)^\beta \le (1+(s-s_0)^2)^\beta \le \tilde{C}(a,\beta)(1+s^2)^\beta \; .
	\end{equation*}
	Using this equivalence, (\ref{estimateinhardy2}) becomes: For $|s_0|<a \le s_1 \in \R $ and $\alpha \ge 0$, there exists a constant $C(a,\alpha)>0$ such that,
	\begin{equation}\label{estimateinhardy3}
	\int_{-s_1}^{s_1} \frac{u^2}{(1+s^2)^{1+\alpha/2}} \d s \le C(a,\alpha) \int_{-s_1}^{s_1} \frac{(\dl_s u)^2}{(1+s^2)^{\alpha/2}} \d s + C(a,\alpha)u(s_0)^2 \; .
	\end{equation}
	
	Set $s_1=\hf t$, $a=\hf$, and $\sigma=\hf\alpha>0$ , let $\xi$ be any non-negative function of $r_*$ which is positive on  $[c,d]\subset \:]-a,a[$ with $c < d$. We integrate over $\mathcal{S}^2$ and $s_0\in [c,d]$ where $\xi$ is bounded below away from zero. Since $[c,d]\subset \:]-a,a[$, we can extend the integration domain in the $s_0$ variable to $]-a,a[$. Thus, there is some $C>0$, which depends on $\xi$ and $\sigma$ only, such that (\ref{hardy1}) holds.
	
	Similarly, in (\ref{estimateinhardy3}), let $\alpha=0$, $a=\frac{3}{4}$, and $s_1=\frac{3}{4} t$ ($t \ge 1$). Integrating over $s_0\in [-\frac{3}{4},\frac{3}{4}]$ and then over the $2$-sphere $\mathcal{S}^2$, we get
	\begin{equation*}
	\int\limits_{\{t\}\times \{|r_*|\le \frac{3}{4}t\}\times \mathcal{S}^2}\hspace{-1cm} \frac{u^2}{1+r^2_*} \d\varsigma \quad \le\qquad C\hspace{-1cm} \int\limits_{\{t\}\times \{|r_*|\le \frac{3}{4}t\}\times \mathcal{S}^2} \hspace{-1cm} (\dl_{r_*}u)^2\d\varsigma \qquad +\qquad C\hspace{-1cm} \int\limits_{\{t\}\times \{|r_*|\le \frac{3}{4}\}\times \mathcal{S}^2}\hspace{-1cm} u^2 \d\varsigma \; .
	\end{equation*}
	Since the continuous function $fr^{-2}$ is positive for all $r_* \in ]-\infty,+\infty[$, then there is some $\delta>0$ such that $\frac{1}{\delta} fr^{-2}>1$ for $|r_*|\le \frac{3}{4}$. If $u$ is of the form (\ref{nonstationarysolutions}), then by Lemma \ref{angularboundlemma} we have,
	\begin{equation*}
	\int\limits_{\{t\}\times \{|r_*|\le \frac{3}{4}t\}\times \mathcal{S}^2}\hspace{-1cm} \frac{u^2}{1+r^2_*} \d\varsigma \quad \le\qquad C\hspace{-1cm} \int\limits_{\{t\}\times \{|r_*|\le \frac{3}{4}t\}\times \mathcal{S}^2} \hspace{-1cm} (\dl_{r_*}u)^2\d\varsigma \qquad +\qquad C\hspace{-1cm} \int\limits_{\{t\}\times \{|r_*|\le \frac{3}{4}\}\times \mathcal{S}^2}\hspace{-1cm} \frac{f}{r^2}|\slnb u|^2 \d\varsigma \; ,
	\end{equation*}
	form which (\ref{hardy2}) follows.
	
\end{proof}

Finally, we derive some identities on $\mathcal{S}^2$. Let $\Theta_1$ , $\Theta_2$ , and $\Theta_3$ be the generators of rotation around the  $x$, $y$, and $z$ axes in $\R^3$ :
\begin{eqnarray*}
	\Theta_1 &=& z\dl_y - y\dl_ z \; , \\
	\Theta_2 &=& z\dl_x - x\dl_ z \; , \\
	\Theta_3 &=& x\dl_y - y\dl_x \; .
\end{eqnarray*}

\begin{lem1}
	Let $u$ be a smooth function on $\mathcal{S}^2$ then,
	\begin{eqnarray}
	|\slnb u|^2 &=& \sum_{i=1}^3 (\Theta_i u)^2 \; , \label{rotationformofslnb} \\
	\sldlt u &=&\sum_{i=1}^3 \Theta_i^2 u \; ,\label{rotationformofsldlt} \\
	\int_{\mathcal{S}^2} |\sldlt u|^2 \d^2 \omega &=& \int_{\mathcal{S}^2}\sum_{i=1}^3 |\slnb \Theta_i u|^2 \d^2 \omega \; . \label{rotationformofintegralsldlt^2}
	\end{eqnarray}
\end{lem1}
\begin{proof}
	We prove the first two identities in spherical coordinates $(\theta,\varphi)$. On $\mathcal{S}^2$ we have,
	\begin{eqnarray*}
		x &=& \sin(\theta)\cos(\varphi) \; , \\
		y &=& \sin(\theta)\sin(\varphi) \; , \\
		z &=& \cos(\theta)\; ,
	\end{eqnarray*}
	so,
	\begin{eqnarray*}
		\Theta_1 &=& \sin(\varphi)\dl_\theta + \cot(\theta)\cos(\varphi)\dl_\varphi \; , \\
		\Theta_2 &=& \cos(\varphi)\dl_\theta - \cot(\theta)\sin(\varphi)\dl_\varphi \; , \\
		\Theta_3 &=& \dl_\varphi \; .
	\end{eqnarray*}
	Thus, noting that
	\begin{eqnarray*}
		|\slnb u|^2 &=& (\dl_\theta u)^2 +\frac{1}{\sin(\theta)^2}(\dl_\varphi u)^2 , \\
		\sldlt u &=& \dl_\theta^2 u +\frac{1}{\sin(\theta)^2}\dl^2_\varphi u +\cot(\theta)\dl_\theta u \; ,
	\end{eqnarray*}
	a straightforward calculation gives (\ref{rotationformofslnb}) and (\ref{rotationformofsldlt}). To show (\ref{rotationformofintegralsldlt^2}) we use the following properties of the commutators of the generators of rotation which are a direct calculation,
	\begin{eqnarray*}
		[\Theta_1  , \Theta_2] &=& \Theta_3 \; ,\\
		{[\Theta_2 , \Theta_3]} &=& \Theta_1 \; ,\\
		{[\Theta_3 , \Theta_1]} &=& \Theta_2 \; ,
	\end{eqnarray*}
	or using the Levi-Civita symbol,
	\begin{equation}\label{rotationcommutators}
	\Theta_i   \Theta_j - \Theta_j \Theta_i=[\Theta_i  , \Theta_j] ={\varepsilon_{ij}}^k \Theta_k=\sum_{k=1}^3\varepsilon_{ijk} \Theta_k \; ,
	\end{equation}
	where the Levi-Civita symbol is
	\begin{equation}\label{Levi-Civitasymbol}
	{\varepsilon_{ij}}^k=\varepsilon_{ijk}= \left\{
	\begin{array}{ll}
	+1 & \hbox{if } (i, j, k) \hbox{ is an even permutation of } (1,2,3); \\
	-1 & \hbox{if } (i,j,k) \hbox{ is an odd permutation of } (1,2,3); \\
	0 & \hbox{otherwise}
	\end{array}
	\right.
	\end{equation}
	We also need the skew-symmetry of the rotation generators:
	\begin{equation}\label{rotationantiselfadjointness}
	\int_{\mathcal{S}^2}v\Theta_i u \d^2\omega= - \int_{\mathcal{S}^2}u\Theta_i v \d^2\omega \; .
	\end{equation}
	To prove this, it is enough to show that it holds for $\dl_\varphi$, since any other $\Theta_i$ can be changed into $\dl_\varphi$ via a permutation of the variables $(x,y,z)$. Indeed, let $u$ be a smooth function on ${\mathcal{S}^2}$,
	\begin{equation*}
	\int_{\mathcal{S}^2}\Theta_i u \d^2\omega= \int\limits_{\theta=0}^\pi \int\limits_{\varphi=0}^{2\pi} (\dl_\varphi u) \sin (\theta) \d\theta \d \varphi=\int\limits_{\theta=0}^\pi (u(2\pi)-u(0)) \sin (\theta) \d\theta=0\; ,
	\end{equation*}
	since $u(2\pi)=u(0)$ as $\varphi=0$ and $\varphi=2\pi$ correspond to the same points on the sphere. Applying this to the product $vu$ in place of $u$ we get (\ref{rotationantiselfadjointness}). Now we prove (\ref{rotationformofintegralsldlt^2}). We have,
	\begin{eqnarray*}
		\int_{\mathcal{S}^2} |\sldlt u|^2 \d^2 \omega  &=& \int_{\mathcal{S}^2} \left(\sum_{i=1}^3 \Theta_i^2 u\right)\left(\sum_{j=1}^3 \Theta_j^2 u\right)\d^2 \omega\qquad\mathrm{by}\; (\ref{rotationformofsldlt})\\
		&=&-\sum_{i,j=1}^3\int_{\mathcal{S}^2} \Theta_i u \Theta_i \Theta_j^2 u \d^2 \omega \qquad\mathrm{by}\; (\ref{rotationantiselfadjointness})\\
		&=&-\sum_{i,j=1}^3\int_{\mathcal{S}^2} \Theta_i u ({\varepsilon_{ij}}^k \Theta_k+\Theta_j \Theta_i) \Theta_j u \d^2 \omega  \qquad\mathrm{by}\; (\ref{rotationcommutators})\\
		&=&-\sum_{i,j=1}^3\int_{\mathcal{S}^2}( {\varepsilon_{ij}}^k\Theta_i u \Theta_k\Theta_j u+\Theta_i u \Theta_j \Theta_i \Theta_j u) \d^2 \omega \\
		&=&-\sum_{i,j=1}^3\int_{\mathcal{S}^2} {\varepsilon_{ij}}^k\Theta_i u \Theta_k\Theta_j u \d^2 \omega +\sum_{i,j=1}^3\int_{\mathcal{S}^2} \Theta_j \Theta_i u \Theta_i \Theta_j u \d^2 \omega \\
		&=&-\sum_{i,j=1}^3\int_{\mathcal{S}^2} {\varepsilon_{ij}}^k\Theta_i u \Theta_k\Theta_j u \d^2 \omega +\sum_{i,j=1}^3\int_{\mathcal{S}^2}(\Theta_i \Theta_j u)^2 \d^2 \omega \\
		&&\hspace{2cm} -\sum_{i,j=1}^3\int_{\mathcal{S}^2} {\varepsilon_{ij}}^k\Theta_k u \Theta_i\Theta_j u \d^2 \omega
	\end{eqnarray*}
	By the antisymmetric nature of the Levi-Civita symbol, the first and the last integrals cancel out. So, we are left with
	\begin{equation*}
	\int_{\mathcal{S}^2} |\sldlt u|^2 \d^2 \omega=\sum_{i,j=1}^3\int_{\mathcal{S}^2}(\Theta_i \Theta_j u)^2 \d^2 \omega \; .
	\end{equation*}
	By (\ref{rotationformofslnb}) the right hand side is nothing but
	\begin{equation*}
	\int_{\mathcal{S}^2}\sum_{i=1}^3 |\slnb \Theta_i u|^2 \d^2 \omega \; .
	\end{equation*}
\end{proof}

It is an easy calculation to show that the $\Theta_i$'s are Killing vector fields of the full metric (and of the Euclidean metric on $\mathcal{S}^2$) and thus satisfy the Killing equation,
\begin{equation*}
\nabla_\alpha (\Theta_i)_\beta + \nabla_\beta (\Theta_i)_\alpha=\dl_\alpha  (\Theta_i)_\beta + \dl_\beta (\Theta_i)_\alpha -2\tilde{\Gamma}_{\alpha \beta}^\gamma (\Theta_i)_\gamma = 0 \; .
\end{equation*}

Consequently, we can make use of the following important property of Killing vector fields.
\begin{lem1}
	Let $X$ be any smooth vector field on a Lorentzian manifold endowed with the Levi-Civita connection. In abstract index notation, let the deformation tensor of $X$ be
	\begin{equation*}
	{}^{(X)}\pi_{ab}=\nabla_a X_b + \nabla_b X_a \; .
	\end{equation*}
	Let $u$ be a smooth function then,
	\begin{equation*}
	\nabla^a\nabla_a (X^b\nabla_b u)=(\nabla^a{}^{(X)}\pi_{ab}) \nabla^b u +  \hf(\nabla_b{}^{(X)}{\pi_a}^a) \nabla^b u +  {}^{(X)}\pi_{ab} \nabla^a \nabla^b u + X^b\nabla_b\nabla^a \nabla_a u \; .
	\end{equation*}
	Thus, if $X$ is Killing, i.e. satisfies the Killing equation:
	\begin{equation*}
	{}^{(X)}\pi_{ab}=0 \; ,
	\end{equation*}
	then $X$  commutes with the wave operator,
	\begin{equation}\label{Killingcommuteswithbox}
	\nabla^a\nabla_a (X^b\nabla_b u)=\square(X(u))=X(\square u)=X^b\nabla_b\nabla^a \nabla_a u\; .
	\end{equation}
\end{lem1}
\begin{proof}
	The proof is a direct computation and then using the symmetries of the curvature tensor of the metric. We have,
	\begin{eqnarray*}
		\nabla^a\nabla_a (X^b\nabla_b u) &=&\nabla^a((\nabla_aX^b)\nabla_b u+ X^b\nabla_a\nabla_b u)  \\
		&=&(\nabla^a\nabla_a X^b)\nabla_b u +2(\nabla^a X^b)(\nabla_a \nabla_b u)+ X^b\nabla^a\nabla_a \nabla_b u\; .
	\end{eqnarray*}
	Since $\nabla$ is torsion-free then,
	\begin{equation*}
	\nabla_a \nabla_b u=\nabla_b \nabla_a u \; ,
	\end{equation*}
	and so,
	\begin{equation*}
	2(\nabla^a X^b)(\nabla_a \nabla_b u)={}^{(X)}\pi_{ab} \nabla^a \nabla^b u \; .
	\end{equation*}
	The connection being torsion-free also implies also that for  a vector field $Y$,
	\begin{equation*}
	{R_{abc}}^d Y^c=(\nabla_a \nabla_b-\nabla_b \nabla_a)Y^d \; ,
	\end{equation*}
	where ${R_abc}^d$ is the curvature tensor. Thus,
	\begin{equation*}
	X^b\nabla^a\nabla_a \nabla_b u=X^b\nabla_a\nabla_b\nabla^a u = {R_{abc}}^a X^b \nabla^c u +X^b\nabla_b\nabla^a \nabla_a u \; .    \end{equation*}
	Adding and subtracting terms we also have,
	\begin{eqnarray*}
		(\nabla^a\nabla_a X^b)\nabla_b u &=& (\nabla^a{}^{(X)}\pi_{ab}) \nabla^b u   - (\nabla_a\nabla_b X^a)\nabla^b u\\
		&=& (\nabla^a{}^{(X)}\pi_{ab}) \nabla^b u - {R_{abc}}^a X^c \nabla^b u +  \hf(\nabla_b{}^{(X)}{\pi_a}^a) \nabla^b u \; .
	\end{eqnarray*}
	Everything is there except that we have additional curvature terms. Due to the compatibility of the connection with the metric,
	\begin{equation*}
	R_{abcd}=-R_{bacd}=R_{badc}=R_{dcba} \; .
	\end{equation*}
	and the terms involving the curvature tensor cancel out.
\end{proof}

Therefore, as $\Theta_i$'s are Killing, they commute with the wave operator. Since they also commute with $r$, we have
\begin{equation}\label{Thetacommutewithbox}
\widetilde{\square}(\Theta_i(u))=\Theta_i(\widetilde{\square} u) \; ,
\end{equation}
where $\widetilde{\square}$ was defined in (\ref{boxtimesoperator}).
This implies that $\Theta_i u$ is a solution of $\tilde{\square}v=0$ when $u$ is. Observing that the $\Theta_i$'s commute with $\sldlt$ since they are also Killing on the Sphere, under the assumptions of Lemma \ref{angularboundlemma} we have,
\begin{equation}\label{angularboundThetau}
\int_{\mathcal{S}^2} (\Theta_i u)^2 \d^2 \omega \le \hf \int_{\mathcal{S}^2} |\slnb \Theta_i u|^2 \d^2 \omega \; ,
\end{equation}
which can be directly seen from the proof of the Lemma by commuting with $\sldlt$ after the first integration by parts. Summing over $i$ then using (\ref{rotationformofslnb}) and (\ref{rotationformofintegralsldlt^2}) adds the following inequality to Lemma \ref{angularboundlemma},
\begin{equation}\label{angularboundslnb}
\int_{\mathcal{S}^2} u^2 \d^2 \omega \le \hf \int_{\mathcal{S}^2} |\slnb u|^2 \d^2 \omega \le \frac{1}{4}\int_{\mathcal{S}^2} |\sldlt u|^2 \d^2  \omega \; .
\end{equation}

\subsubsection{Uniform Bound on the Conformal Charge}

We use a radial multiplier:
\begin{equation*}
\gamma_{(u)} = g\dl_{r_*}u + \hf(\dl_{r_*} g)u\; ,
\end{equation*}
where $u$ is a smooth solution of (\ref{WavePhi0}) and $g$ in this paragraph denotes a function of the $t$ and $r_*$ variables only. We start with the following lemma.
\begin{lem1}\label{Egamma}
	Let $\gamma_{(u)}$ be as above. Set
	\begin{equation*}
	E_\gamma[u] (t)= \int_{\Sigma_t} \gamma_{(u)}\dl_t u \d \varsigma \; ,
	\end{equation*}
	then,
	\begin{eqnarray}\label{Egammaestimate}
	2E_\gamma[u] (t_2) - 2E_\gamma[u] (t_1)&=&-\hspace{-0.5cm}\int\limits_{[t_1,t_2]\times \R \times \mathcal{S}^2}(2(\dl_{r_*}g)(\dlr u)^2 -\hf u^2(\dlr^3 g)-g |\slnb u|^2 \dlr V \nonumber\\
	&& \qquad\quad\qquad  -2(\dl_t u)(\dl_t g)(\dlr u) -  u(\dl_t u) \dl_{tr_*}^2 g )\d x \; .
	\end{eqnarray}
\end{lem1}
\begin{proof}
	We use dot for time derivative and prime for derivative with respect to the radial variable $r_*$. Since
	\begin{equation*}
	2E_\gamma[u] (t_2) - 2E_\gamma[u] (t_1)=\int\limits_{[t_1,t_2]} 2\dot{E}_\gamma[u] (t) \d t= \int\limits_{[t_1,t_2]\times \R \times \mathcal{S}^2} \dl_t(2\dot{u}\gamma_{(u)})\d x \; .
	\end{equation*}
	Using (\ref{WavePhi0}) we calculate,
	\begin{equation*}
	\dl_t(2\dot{u}\gamma_{(u)})= 2\gamma_{(u)} u'' +2 \frac{f}{r^2}\gamma_{(u)} \sldlt u + 2 \dot{u} \dot{\gamma}_{(u)}\; .
	\end{equation*}
	But
	\begin{eqnarray*}
		\gamma_{(u)} \sldlt u &=&  \slnb \cdot (\gamma_{(u)} \slnb u) - \slnb\gamma_{(u)} \cdot \slnb u  \; ,   \\
		2\gamma_{(u)} u'' &=& \dlr(gu''+u u'g') -2u''g'-u u'g'' \; , \\
		2 \dot{u}\dot{\gamma}_{(u)}&=& \dlr(g\dot{u}^2) + u \dot{u} \dot{g}' + 2 \dot{u}\dot{g}u' \; ,
	\end{eqnarray*}
	and
	\begin{eqnarray*}
		\slnb\gamma_{(u)} \cdot \slnb u &=& \hf \dlr(g|\slnb u|^2) \; , \\
		u u' g'' &=& \hf \dlr(u^2 g'') -\hf u^2 g''' \; .
	\end{eqnarray*}
	Thus
	\begin{eqnarray*}
		\dl_t (2\gamma_{(u)} \dot{u})&=& \dlr\left(u''g+uu'g'-\frac{f}{r^2}g|\slnb u|^2 - \hf u^2 g''+g\dot{u}^2\right)+2\slnb \cdot \left(\frac{f}{r^2} \gamma_{(u)} \slnb u \right) \\
		&&\qquad+\hf u^2 g''' + u \dot{u} \dot{g}' + 2 \dot{u}\dot{g}u' -\dlr\left(\frac{f}{r^2}\right)g|\slnb u|^2 - 2 g' u'^2  \; .
	\end{eqnarray*}
	Integrating over $[t_1,t_2]\times \R \times \mathcal{S}^2$ we get (\ref{Egammaestimate}) since the integrals of the first two terms are zero.
\end{proof}

We are now ready to establish the uniform bound estimates of this section.

\begin{prop1}[Uniform Bound]\label{unifomboundTHEOREM}
	Let $u$ be a smooth solution of (\ref{WavePhi0}) of the form (\ref{nonstationarysolutions}), and $\chi$ any compactly supported smooth function. There exists a constant $C>0$ such that for all $t\ge0$ we have,
	\begin{eqnarray}
	E_\mathcal{C}[u](t) &\le& C(E_\mathcal{C}[u](0)+E[\sldlt^2 u](0))\; , \label{unifomboundenergyconformal}\\
	\int\limits_{[0,+\infty[\times \R \times \mathcal{S}^2} t\chi |\slnb u|^2 \d t \d r_* \d^2 \omega  &\le&  C(E_\mathcal{C}[u](0)+E[\sldlt^2 u](0)) \; . \label{unifomboundtrappingterm}  
	\end{eqnarray}
\end{prop1}
\begin{proof}
	Recall that
	$$V=\frac{f}{r^2} \;.$$
	Let $t\ge 1$. Following \cite{blue_decay_2008}, we set
	\begin{eqnarray*}
		h(r_*)&=&\int^{r_*}_0 \frac{1}{(1+(\epsilon y)^2)^{\sigma}}\d y \; ,\\
		g&=&t\mu h \; ,
	\end{eqnarray*}
	where, $\sigma \in [1,2]$, $\epsilon>0$, and $\mu=\tilde{\mu}(\frac{r_*}{t})$ with $\tilde{\mu}$ a smooth  function with compactly supported in $]-\frac{3}{4},\frac{3}{4}[$ which is identically $1$ on $[-\hf ,\hf]$. Note that $h$ is bounded. 
	
	The idea now is bound  $E_\gamma[u]$ and the error terms in Lemma \ref{Egamma} by the local energy  $E_\ell[u]$. The bounds on the error terms are uniform, while the bound on $E_\gamma[u]$ is linear. Using the Hardy estimates and these bounds we get a linear bound on the error term of the conformal charge and thus on the conformal charge itself. This linear bound can be improved to a uniform one using the fact that inside the light cone, the  conformal charge density bounds the energy density times $t^2$.
	As before, dot and prime indicate differentiation with respect to $t$ and $r_*$ respectively. 
	
	We calculate the terms in (\ref{Egammaestimate}) to get,
	\begin{eqnarray}
	2E_\gamma[u] (t_2) - 2E_\gamma[u] (t_1)=&-&\hspace{-0.5cm}\int\limits_{[t_1,t_2]\times \R \times \mathcal{S}^2} (2t\mu h' u'^2 - t \mu V' h |\slnb u|^2) \d x  \label{Egammaestimate1}\\
	&+&\hf\hspace{-0.5cm}\int\limits_{[t_1,t_2]\times \R \times \mathcal{S}^2} t\mu h''' u^2 \d x \label{Egammaestimate2}\\
	&-&\hspace{-0.5cm}\int\limits_{[t_1,t_2]\times \R \times \mathcal{S}^2} 2th\mu' u'^2 \d x \label{Egammaestimate3}\\
	&+&\hf\hspace{-0.5cm}\int\limits_{[t_1,t_2]\times \R \times \mathcal{S}^2} tu^2(3\mu' h'' +3 \mu'' h' + \mu''' h)\d x \label{Egammaestimate4}\\
	&+&\hspace{-0.5cm}\int\limits_{[t_1,t_2]\times \R \times \mathcal{S}^2} (2\dot{u}\dot{g}u' + u\dot{u}\dot{g}')\d x \; . \label{Egammaestimate5}
	\end{eqnarray}
	All these integrals are in fact over the domain $[t_1,t_2]\times \{|r_*|\le \frac{3}{4}t\}\times \mathcal{S}^2$ because of $\mu$ and its derivatives.
	
	We start by bounding (\ref{Egammaestimate3})-(\ref{Egammaestimate5}) by the integral of local energy $E_\ell$. To do so, we use (\ref{normequivalenceMaxandEuclid}): 
	
	Since $\tilde{\mu}$ is constant on $[-\hf,\hf]$, its derivatives (denoted $\frac{\d^n \tilde{\mu}}{{\d \tau}^n}$) are supported away from zero, namely in $\left]-\frac{3}{4},-\hf\right] \cup \left[\hf,\frac{3}{4}\right[$. For $\sigma \in [1,2]$ we have,
	\begin{eqnarray*}
		|h|&\le&C \; , \\
		|h'|&=&\frac{1}{(1+(\epsilon r_*)^2)^{\sigma}}\sim \frac{1}{(1+r_*^2)^{\sigma}} \le \frac{1}{1+ r_*^2} \le C \; , \\
		|h'|&=&\frac{1}{(1+(\epsilon r_*)^2)^{\sigma}}\sim \frac{1}{(1+|r_*|)^{2\sigma}} \le \frac{1}{1+ |r_*|} \le C \; ,
	\end{eqnarray*}
	and for $t\ge 1$,
	\begin{equation*}
	\dlr^n \mu =\frac{1}{t^n}\frac{\d^n\tilde{\mu}}{{\d \tau}^n}\left(\frac{r_*}{t}\right)\le \frac{C}{t^n}\le C\;  ,
	\end{equation*}
	and if in addition $\left|\frac{r_*}{t}\right|\in \left[\hf,\frac{3}{4}\right[$  then,
	\begin{eqnarray*}
		\frac{t}{|r_*|}&\sim&\frac{t}{1+|r_*|} \; , \\
		1<\frac{4}{3}&\le& \frac{t}{|r_*|}\le 2 \; , \\
		|h''| &=& \frac{2\epsilon^2 \sigma |r_*|}{(1+(\epsilon r_*)^2)^{\sigma+1}} \sim \frac{1}{(1+r_*^2)^{\sigma + \hf}}\le \frac{1}{1+r_*^2}\le C \; .
	\end{eqnarray*}
	
	Therefore, by (\ref{hardy2}),
	\begin{align*}
	\left| \; \int\limits_{\{t\} \times \R \times \mathcal{S}^2} 2th\mu' u'^2 \d \varsigma\right| &\le&& \qquad C \hspace{-1cm}\int\limits_{\{t\} \times \{\hf t\le |r_*|\le \frac{3}{4}t\} \times \mathcal{S}^2}\hspace{-1cm} u'^2 \d \varsigma \le C E_\ell[u](t) \; , \\
	\left|\hf\int\limits_{\{t\} \times \R \times \mathcal{S}^2} 3tu^2\mu' h'' \d \varsigma\right|  &\le&& \qquad C \hspace{-1cm}\int\limits_{\{t\} \times \{\hf t\le|r_*|\le \frac{3}{4}t\} \times \mathcal{S}^2}3tu^2\left(\frac{1}{t}\right)\left(\frac{1}{1+r_*^2}\right) \d \varsigma \le C E_\ell[u](t) \;  \\
	\left|\hf\int\limits_{\{t\} \times \R \times \mathcal{S}^2} 3tu^2\mu'' h' \d \varsigma \right| &\le&& \qquad C \hspace{-1cm}\int\limits_{\{t\} \times \{\hf t\le|r_*|\le \frac{3}{4}t\} \times \mathcal{S}^2}3tu^2\left(\frac{1}{t^2}\right)\left(\frac{1}{1+r_*^2}\right) \d \varsigma \le C E_\ell[u](t) \;  \\
	\left| \hf\int\limits_{\{t\} \times \R \times \mathcal{S}^2} tu^2\mu''' h \d \varsigma \right|&\le&& \qquad C \hspace{-1cm}\int\limits_{\{t\} \times \{\hf t\le|r_*|\le \frac{3}{4}t\} \times \mathcal{S}^2}tu^2\left(\frac{1}{t^3}\right)\frac{t^2}{|r_*|^2} \d \varsigma \qquad \qquad\left(\frac{t^2}{|r_*|^2}>1\right)  \\
	&\le&& \qquad C \hspace{-1cm}\int\limits_{\{t\} \times \{\hf t\le|r_*|\le \frac{3}{4}t\} \times \mathcal{S}^2}tu^2\left(\frac{1}{t^3}\right)\frac{t^2}{1+r_*^2} \d \varsigma \le C E_\ell[u](t) \;  \; .
	\end{align*}
	
	This treats (\ref{Egammaestimate3}) and (\ref{Egammaestimate4}). For (\ref{Egammaestimate5}) we do the same.
	\begin{eqnarray*}
		\dot{g}&=&h\left(\mu-\frac{1}{t}\frac{\d\tilde{\mu}}{{\d \tau}}\left(\frac{r_*}{t}\right)\right)\le C \; , \\
		\dot{g}'&=&h\left(\frac{1}{t}\frac{\d\tilde{\mu}}{{\d \tau}}\left(\frac{r_*}{t}\right)-\frac{1}{t^2}\frac{\d^2\tilde{\mu}}{{\d \tau^2}}\left(\frac{r_*}{t}\right)\right) + h'\left(\mu-\frac{1}{t}\frac{\d\tilde{\mu}}{{\d \tau}}\left(\frac{r_*}{t}\right)\right)\\
		&\le& C\left(\frac{1}{t}\right)\frac{t}{1+|r_*|}+ C\frac{1}{1+|r_*|} \le C\frac{1}{1+|r_*|} \; , \qquad \mathrm{for}\quad \left|\frac{r_*}{t}\right|\in \left[\hf,\frac{3}{4}\right[ \; .
	\end{eqnarray*}
	Now using this and (\ref{normequivalenceMaxandEuclid}) and the inequality $(a-b)^2\ge0$, we have:
	\begin{eqnarray*}
		\left|\int\limits_{\{t\}\times \R \times \mathcal{S}^2} (2\dot{u}\dot{g}u' + u\dot{u}\dot{g}')\d \varsigma\right| &\le& \qquad C\hspace{-1cm}\int\limits_{\{t\} \times \{|r_*|\le \frac{3}{4}t\} \times \mathcal{S}^2}\hspace{-1cm}|\dot{u}u'|\d \varsigma+\qquad C \hspace{-1cm}\int\limits_{\{t\} \times \{\hf t\le|r_*|\le \frac{3}{4}t\} \times \mathcal{S}^2}\hspace{-1cm}|\dot{u}| \frac{|u|}{1+|r_*|}\d \varsigma \\
		&\le& \qquad C\hspace{-1cm}\int\limits_{\{t\} \times \{|r_*|\le \frac{3}{4}t\} \times \mathcal{S}^2}\hspace{-1cm}(|\dot{u}|^2 + |u'|^2)\d \varsigma   \\ 
		&&\hspace{3cm}\qquad C \hspace{-1cm}+\int\limits_{\{t\} \times \{|r_*|\le \frac{3}{4}t\} \times \mathcal{S}^2}\left(|\dot{u}|^2 + \frac{u^2}{(1+|r_*|)^2}\right)\d \varsigma \\
		&\le& C E_\ell[u](t) \; .
	\end{eqnarray*}
	The same calculation gives,
	\begin{equation}\label{EgammalesstElocal}
	|E_\gamma[u](t)|\le C t E_\ell[u](t) \; .
	\end{equation}
	
	Recapitulating, we have shown that for $t\ge 1$,
	
	\begin{eqnarray*}
		2E_\gamma[u] (t_2) - 2E_\gamma[u] (t_1)\le\quad -\hspace{-0.5cm}\int\limits_{[t_1,t_2]\times \R \times \mathcal{S}^2}\hspace{-0.5cm} (2t\mu h' u'^2 &-& t \mu V' h |\slnb u|^2) \d x \quad + \quad\hf\hspace{-0.5cm}\int\limits_{[t_1,t_2]\times \R \times \mathcal{S}^2}\hspace{-0.5cm} t\mu h''' u^2 \d x \quad\\
		&+& \quad C \int_{t_1}^{t_2} E_\ell[u](t) \d t \; .
	\end{eqnarray*}
	We still need to treat the first two terms on the right hand side. Since $h'>0$, $h$ is increasing and has the sign of $r_*$. By the argument used in the proof of Proposition \ref{1photonsphere},
	$$V'=\frac{f}{r^3}(rf'-2f)$$
	has the opposite sign of $r_*$. So, $hV'\le0$ with the equality holding only at $r_*=0$ where both functions vanish. Thus, the first term is negative.
	We need to control the second. We have,
	\begin{equation*}
	|h'''| = \frac{2\epsilon^2 \sigma |(2\sigma +1)(\epsilon r_*)^2-1|}{(1+(\epsilon r_*)^2)^{\sigma+2}}\le \frac{2\epsilon^2 \sigma ((2\sigma +1)(\epsilon r_*)^2+2\sigma +1)}{(1+(\epsilon r_*)^2)^{\sigma+2}}= \frac{2\epsilon^2 \sigma (2\sigma +1)}{(1+(\epsilon r_*)^2)^{\sigma+1}} \; ,
	\end{equation*}
	so,
	\begin{equation}\label{boundingh'''}
	\hf\hspace{-0.5cm}\int\limits_{\{t\}\times \R \times \mathcal{S}^2}\hspace{-0.5cm} t\mu h''' u^2 \d \varsigma \quad\le \quad\epsilon^2 \sigma (2\sigma +1)\hspace{-0.5cm}\int\limits_{\{t\}\times \{|r_*|\le \frac{3}{4}t\} \times \mathcal{S}^2}\hspace{-0.5cm} t\mu  \frac{u^2}{(1+(\epsilon r_*)^2)^{\sigma+1}} \d \varsigma \; ,
	\end{equation}
	We divide the integral into two parts. One which we bound by the local energy and one on which we apply (\ref{hardy1}) and eventually get absorbed by the right hand side of (\ref{Egammaestimate1}). Since $t \ge 1$ and $\sigma \ge 1 $,
	\begin{eqnarray*}
		\epsilon^2 \sigma (2\sigma +1)\hspace{-0.5cm}\int\limits_{\{t\}\times \{\hf t\le|r_*|\le \frac{3}{4}t\} \times \mathcal{S}^2}\hspace{-0.5cm} \mu\frac{t}{(1+(\epsilon r_*)^2)}  \frac{u^2}{(1+(\epsilon r_*)^2)^{\sigma}} \d \varsigma &\le& \quad C \hspace{-0.5cm}\int\limits_{\{t\}\times \{\hf t\le|r_*|\le \frac{3}{4}t\} \times \mathcal{S}^2}\hspace{-0.5cm}\frac{t^2}{r_*^2} \frac{u^2}{(1+(r_*)^2)} \d \varsigma\\
		&\le& C E_\ell[u](t) \; ,
	\end{eqnarray*}
	using (\ref{hardy2}).
	
	For the other part we need to keep track of the constants, in particular those involved in (\ref{hardy1}). $-V'h$ is non-negative and vanishes only at zero, and so does
	\begin{equation*}
	-V'\int_0^{r_*}\frac{1}{(1+y^2)^\sigma}\d y \le -V' h\; .
	\end{equation*}
	So, using Lemma \ref{angularboundlemma} and (\ref{hardy1}),
	\begin{eqnarray}
	\hspace{-2cm}\int\limits_{\{t\}\times \{|r_*| \le \hf t\}\times \mathcal{S}^2} \frac{u^2}{(1+(\epsilon r_*)^2)^{\sigma + 1}} \d\varsigma &\le& C \int\limits_{\{t\}\times \{|r_*|\le \hf t\}\times \mathcal{S}^2} \left(\frac{(\dl_{r_*} u)^2}{(1+(\epsilon r_*)^2)^{\sigma}} \right. \\ &&\hspace{2cm}\left. -V'\left(\int_0^{r_*}\frac{1}{(1+y^2)^\sigma}\d y\right)  u^2\right) \d\varsigma \nonumber \\
	&\le& C \int\limits_{\{t\}\times \{|r_*|\le \hf t\} \times \mathcal{S}^2}\hspace{-0.5cm} (2 h' u'^2  - V' h |\slnb u|^2) \d \varsigma \; , \label{hardyinboundtheorem}
	\end{eqnarray}
	Now, $\mu\equiv 1$ over $\{|r_*|\le \hf t\}$ so,
	\begin{equation*}
	\hf\hspace{-0.5cm}\int\limits_{\{t\}\times \{|r_*|\le \hf t\} \times \mathcal{S}^2}\hspace{-0.5cm} t\mu h''' u^2 \d \varsigma \quad \le \quad\epsilon^2 \sigma (2\sigma +1)C\hspace{-0.5cm} \int\limits_{\{t\}\times \{|r_*|\le \hf t\} \times \mathcal{S}^2}\hspace{-0.5cm} (2t\mu h' u'^2 - t \mu V' h |\slnb u|^2) \d \varsigma \; .
	\end{equation*}
	We choose $\epsilon^2<\frac{1}{2\sigma (2\sigma +1)C}$, so
	\begin{equation*}
	2E_\gamma[u] (t_2) - 2E_\gamma[u] (t_1)\le\quad -\hf\hspace{-0.5cm}\int\limits_{[t_1,t_2]\times \R \times \mathcal{S}^2}\hspace{-0.5cm} (2t\mu h' u'^2 - t \mu V' h |\slnb u|^2) \d x \quad+ \quad C \int_{t_1}^{t_2} E_\ell[u](t) \d t \; ,
	\end{equation*}
	i.e.
	\begin{equation*}
	\hf\hspace{-0.5cm}\int\limits_{[t_1,t_2]\times \R \times \mathcal{S}^2}\hspace{-0.5cm} (2t\mu h' u'^2 - t \mu V' h |\slnb u|^2) \d x \quad \le\quad -2E_\gamma[u] (t_2) + 2E_\gamma[u] (t_1)\quad+ \quad C \int_{t_1}^{t_2} E_\ell[u](t) \d t \; .
	\end{equation*}
	Using (\ref{Thetacommutewithbox}) we see that this estimate holds equally for $\Theta_i u$. Summing over $i$ the estimate for $\Theta_i u$, and using (\ref{rotationformofslnb}) and (\ref{rotationformofintegralsldlt^2}) we get,
	\begin{equation}
	\hf\hspace{-0.5cm}\int\limits_{[t_1,t_2]\times \R \times \mathcal{S}^2}\hspace{-0.5cm} t\mu(2 h' |\slnb u'|^2 - V' h (\sldlt u)^2) \d x  \le  -2E_\gamma[\slnb u] (t_2) + 2E_\gamma[\slnb u] (t_1)+  C \int_{t_1}^{t_2} E_\ell[\slnb u](t) \d t \; . \label{controlonhardy}
	\end{equation}
	where the energy terms denote,
	\begin{equation*}
	E_\gamma[\slnb u] (t):= \sum_{i=1}^3 E_\gamma[\Theta_i u] (t) \qquad ; \qquad E_\ell[\slnb u](t):= \sum_{i=1}^3 E_\ell[\Theta_i u](t)\; ,
	\end{equation*}
	motivated by (\ref{rotationformofslnb}) and (\ref{rotationformofintegralsldlt^2}).
	Now, applying the estimate (\ref{hardy1}) for $\Theta_i u$ then using (\ref{angularboundThetau}), we get an estimate for $\Theta_i u$ analogous to (\ref{hardyinboundtheorem}). Again, summing over $i$ and using (\ref{rotationformofslnb}) and (\ref{rotationformofintegralsldlt^2}), still keeping $t\ge 1$, we have,
	\begin{equation*}
	\int\limits_{\{t\}\times \{|r_*| \le \hf t\}\times \mathcal{S}^2} t\frac{|\slnb u|^2}{(1+(\epsilon r_*)^2)^{\sigma + 1}} \d\varsigma \le C \int\limits_{\{t\}\times \R \times \mathcal{S}^2}\hspace{-0.5cm} t\mu(2 h' |\slnb u'|^2 - V' h (\sldlt u)^2) \d \varsigma \; ,
	\end{equation*}
	Let $a\ge1$ such that $supp(\chi_{trap})\subseteq [-a/2,a/2]$ (see (\ref{controlonconformalcharge})). Then for some constant $C>0$ and all $r_* \in [-a/2,a/2]$
	\begin{equation*}
	\chi_{trap} \le C \frac{1}{(1+(\epsilon r_*)^2)^{\sigma + 1}} \; ,
	\end{equation*}
	and hence for $t\ge a$ ,
	\begin{align}\label{unifomboundtrappingtermproof1}
	\int\limits_{\{t\}\times \R \times \mathcal{S}^2} t\chi_{trap}|\slnb u|^2 \d \varsigma &= \int\limits_{\{t\}\times \{|r_*| \le \hf t\}\times \mathcal{S}^2} t\chi_{trap}|\slnb u|^2 \d \varsigma \\
	&\le C \int\limits_{\{t\}\times \{|r_*| \le \hf t\}\times \mathcal{S}^2} t\frac{|\slnb u|^2}{(1+(\epsilon r_*)^2)^{\sigma + 1}} \d\varsigma \; .
	\end{align}
	Integrating from $a$ to $t'\ge a$, then using (\ref{controlonhardy}), we get
	\begin{equation}\label{boundontrapUNIFORMBOUNDTHEOREM}
	\int\limits_{[a,t']\times \R \times \mathcal{S}^2} t\chi_{trap}|\slnb u|^2 \d x \le C\left|\left.-2E_\gamma[\slnb u] (t)\right|_a^{t'}\right|+C \int_{a}^{t'} E_\ell[\slnb u](t) \d t \; .
	\end{equation}
	Now using (\ref{controlonconformalcharge}), we then have
	\begin{equation*}
	E_\mathcal{C}[u](t')\le E_\mathcal{C}[u](a)+  C\left|\left.-2E_\gamma[\slnb u] (t)\right|_a^{t'}\right|+C \int_{a}^{t'} E_\ell[\slnb u](t) \d t \; .
	\end{equation*}
	By (\ref{EgammalesstElocal}),
	\begin{equation*}
	E_\mathcal{C}[u](t')\le E_\mathcal{C}[u](a)+  C\sup\limits_{t\in[0,t']}(t E_\ell [\slnb u](t))+C \int_{0}^{t'} E_\ell[\slnb u](t) \d t \; .
	\end{equation*}
	Since $V>0$ and $\chi_{trap}$ has a compact support then $\chi_{trap}\le C V$. For $0\le t_2 \le 1$, we have,
	\begin{equation}\label{unifomboundtrappingtermproof2}
	\int\limits_{[t_1,t_2]\times \R \times \mathcal{S}^2} t\chi_{trap}|\slnb u|^2 \d x \le C \int\limits_{[t_1,t_2]\times \R \times \mathcal{S}^2} V |\slnb u|^2 \d x \le C \int_{t_1}^{t_2} E[u](t) \d t \le C \int_{t_1}^{t_2} E_\mathcal{C}[u](t) \d t \; .
	\end{equation}
	And for $t_1 \ge 1$,
	\begin{equation}\label{unifomboundtrappingtermproof21}
	\int\limits_{[t_1,t_2]\times \R \times \mathcal{S}^2} t\chi_{trap}|\slnb u|^2 \d x \le C \int\limits_{[t_1,t_2]\times \R \times \mathcal{S}^2} t^2 V |\slnb u|^2 \d x \le C \int_{t_1}^{t_2} E_\mathcal{C}[u](t) \d t \; .
	\end{equation}
	Therefore, for all $t_1,t_2\ge0$, by (\ref{controlonconformalcharge}) we have
	\begin{equation}\label{exponentialboundontrappingintergal}
	E_\mathcal{C}[u](t_2)-E_\mathcal{C}[u](t_1) \le \int\limits_{[t_1,t_2]\times \R \times \mathcal{S}^2} t\chi_{trap}|\slnb u|^2 \d x \le C \int_{t_1}^{t_2} E_\mathcal{C}[u](t) \d t \; .
	\end{equation}
	By Gronwall's inequality,
	\begin{equation}\label{GronwallEc}
	E_\mathcal{C}[u](t_2) \le E_\mathcal{C}[u](t_1) e^{c(t_2 - t_1)} \; ,
	\end{equation}
	then taking $t_1=0$ and $t_2=a$ we get an exponential bound on the conformal energy,
	\begin{equation*}
	E_\mathcal{C}[u](a) \le C E_\mathcal{C}[u](0)\; ,
	\end{equation*}
	with $C$ depending on $a$ but not on $u$. Thus, for $t' \ge a$,
	\begin{equation}\label{firstestimate}
	E_\mathcal{C}[u](t')\le C \left(E_\mathcal{C}[u](0) + \sup\limits_{t\in[0,t']}(t E_\ell [\slnb u](t))+ \int_{0}^{t'} E_\ell[\slnb u](t) \d t \right) \; .
	\end{equation}
	
	Since $E_\ell[\slnb u](t)\le C E[\slnb u](t)=C E[\slnb u](0)$ by Lemma \ref{conservationofenergy}, this gives a linear bound on the conformal charge,
	\begin{equation}\label{linearbound}
	E_\mathcal{C}[u](t)\le C(1+t)(E_\mathcal{C}[u](0)+ E[\slnb u](0))\; , \qquad\qquad t \ge a \; .
	\end{equation}
	
	Next, we derive an estimate for the local energy. Let $\mathcal{S}_t=\{t\}\times \left[-\frac{3}{4}t,\frac{3}{4}t\right]\times \mathcal{S}^2$, and recall (\ref{dotnotation3}). Integrating by parts then using Cauchy-Schwartz inequality , then a double integration by parts followed by another Cauchy-Schwartz inequality, we have:
	\begin{eqnarray*}
		E_\ell[\slnb u](t)&=& \int_{\mathcal{S}_t}\left( |\slnb \dot{u}|^2 +|\slnb u'|^2 + V (\slnb \cdot \slnb u)^2 \right) \d \varsigma \\
		&=& \int_{\mathcal{S}_t}\left((-\sldlt \dot{u}) \dot{u} +(-\sldlt u')u' + V (-\slnb  \sldlt u)\cdot \slnb u \right) \d \varsigma \\
		&\le&\left(\int_{\mathcal{S}_t}(\sldlt \dot{u})^2\d \varsigma\right)^\hf \left(\int_{\mathcal{S}_t}\dot{u}^2\d \varsigma\right)^\hf +\left(\int_{\mathcal{S}_t}(\sldlt u')^2\d \varsigma\right)^\hf \left(\int_{\mathcal{S}_t}u'^2\d \varsigma\right)^\hf \\
		& &+ \left(\int_{\mathcal{S}_t}V |\slnb  \sldlt u|^2\d \varsigma\right)^\hf \left(\int_{\mathcal{S}_t}V|\slnb u|^2\d \varsigma\right)^\hf \\
		&\le& \left(\left(\int_{\mathcal{S}_t}(\sldlt^2 \dot{u})^2\d \varsigma\right)^\hf \left(\int_{\mathcal{S}_t}\dot{u}^2\d \varsigma\right)^\hf\right)^\hf \left(\int_{\mathcal{S}_t}\dot{u}^2\d \varsigma\right)^\hf \\
	\end{eqnarray*}
	
	\begin{eqnarray*}
		&&+ \left(\left(\int_{\mathcal{S}_t}(\sldlt^2 u')^2\d \varsigma\right)^\hf \left(\int_{\mathcal{S}_t}u'^2\d \varsigma\right)^\hf\right)^\hf \left(\int_{\mathcal{S}_t}u'^2\d \varsigma\right)^\hf \\
		&&+\left(\left(\int_{\mathcal{S}_t}V|\slnb \sldlt^2 u|^2\d \varsigma\right)^\hf \left(\int_{\mathcal{S}_t}V|\slnb u|^2\d \varsigma\right)^\hf\right)^\hf \left(\int_{\mathcal{S}_t}V|\slnb u|^2\d \varsigma\right)^\hf \\
		&=& \left(\int_{\mathcal{S}_t}(\sldlt^2 \dot{u})^2\d \varsigma\right)^\frac{1}{4} \left(\int_{\mathcal{S}_t}\dot{u}^2\d \varsigma\right)^{\frac{3}{4}}+ \left(\int_{\mathcal{S}_t}(\sldlt^2 u')^2\d \varsigma\right)^\frac{1}{4} \left(\int_{\mathcal{S}_t}u'^2\d \varsigma\right)^{\frac{3}{4}}\\
		&&+ \left(\int_{\mathcal{S}_t}V|\slnb \sldlt^2 u|^2\d \varsigma\right)^\frac{1}{4} \left(\int_{\mathcal{S}_t}V|\slnb u|^2\d \varsigma\right)^{\frac{3}{4}} \; .
	\end{eqnarray*}
	Using H\"{o}lder's inequality,
	\begin{equation*}
	\sum_{k} \left|x_k^\frac{1}{4}y_k^\frac{3}{4}\right| \le \left|\sum_{k}x_k\right|^\frac{1}{4}\left|\sum_{k}y_k\right|^\frac{3}{4} \; ,
	\end{equation*}
	on the last line of the above estimate, we arrive at,
	\begin{equation*}
	E_\ell[\slnb u](t)\le (E_\ell[\sldlt^2 u](t))^\frac{1}{4}(E_\ell[u](t))^\frac{3}{4} \; .
	\end{equation*}
	But $E_\ell[\sldlt^2 u](t)\le C E[\sldlt^2 u](t)$, and for $|r_*|< \frac{3}{4}t$,
	\begin{equation*}
	\frac{t-r_*}{t}>\frac{1}{4} \qquad;\qquad \frac{t+r_*}{t}>\frac{1}{4} \; ,
	\end{equation*}
	so using (\ref{conformaldensitypositive}) we have,
	\begin{eqnarray*}
		\frac{e_{\mathcal{C}}}{t^2} &\ge&\frac{1}{64}\left((\dot{u} + u')^2+(\dot{u} - u')^2 \right) +\hf (1+\frac{r_*^2}{t^2})|\slnb u|^2 V + \frac{e}{t^2}  \\
		&=& \frac{1}{32}\left(\dot{u}^2+u'^2 \right)+\frac{1}{32}V|\slnb u|^2+(\frac{15}{32}+\frac{r_*^2}{2t^2})|\slnb u|^2 V + \frac{e}{t^2}  \\
		&\ge& \frac{e}{32} \; .
	\end{eqnarray*}
	Thus,
	\begin{equation*}
	E_\ell[u](t)\le C \frac{E_\mathcal{C}[u](t)}{t^2} \; .
	\end{equation*}
	Therefore,
	\begin{equation}\label{estimateonEell}
	E_\ell[\slnb u](t)\le C\left(E[\sldlt^2 u](t)\right)^\frac{1}{4}\left(\frac{E_\mathcal{C}[u](t)}{t^2}\right)^\frac{3}{4} \; .
	\end{equation}
	Substituting for $E_\ell[\slnb u](t)$ in (\ref{firstestimate}) using (\ref{estimateonEell}) we have, for $t'\ge a$
	\begin{eqnarray*}
		E_\mathcal{C}[u](t') &\le& C\left(E_\mathcal{C}[u](0)+ \sup\limits_{t\in[a,t']}\left(t\left(E[\sldlt^2 u](t)\right)^\frac{1}{4}\left(\frac{E_\mathcal{C}[u](t)}{t^2}\right)^\frac{3}{4}\right)\right.\\
		&& \hspace{6cm}\left.+\int_a^{t'}\left(E[\sldlt^2 u](t)\right)^\frac{1}{4} \left(\frac{E_\mathcal{C}[u](t)}{t^2} \right)^\frac{3}{4} \d t\right)\ \; ,
	\end{eqnarray*}
	by Lemma \ref{conservationofenergy} and the linear bound (\ref{linearbound}),
	\begin{eqnarray*}
		E_\mathcal{C}[u](t') &\le&C\left(E_\mathcal{C}[u](0)+ \left(E[\sldlt^2 u](0)\right)^\frac{1}{4} \Big(E_\mathcal{C}[u](0)+ E[\slnb u](0)\Big)^\frac{3}{4} \sup\limits_{t\in[a,t']}\left(t\left(\frac{1+t}{t^2} \right)^\frac{3}{4}\right) \right.\\
		&&\hspace{2.5cm} \left. + \left(E[\sldlt^2 u](0)\right)^\frac{1}{4} \Big(E_\mathcal{C}[u](0)+ E[\slnb u](0)\Big)^\frac{3}{4} \int_a^{t'} \left(\frac{1+t}{t^2} \right)^\frac{3}{4} \d t \right) \; ,
	\end{eqnarray*}
	but for $t\ge a \ge 1$,
	\begin{equation*}
	\left(\frac{1+t}{t^2} \right)^\frac{3}{4}\le C t^{-\frac{3}{4}} \; ,
	\end{equation*}
	so,
	\begin{equation*}
	E_\mathcal{C}[u](t) \le C\left(E_\mathcal{C}[u](0)+ t^{\frac{1}{4}}\left(E[\sldlt^2 u](0)\right)^\frac{1}{4} \Big(E_\mathcal{C}[u](0)+ E[\slnb u](0)\Big)^\frac{3}{4} \right) \; .
	\end{equation*}
	From (\ref{angularboundslnb}) one has,
	\begin{equation*}
	E[u](t) \le  E[\slnb u](t) \le E[\sldlt u](t) \; ,
	\end{equation*}
	but as it is easily seen that $\sldlt u$ is again a solution of (\ref{WavePhi0}) of the form (\ref{nonstationarysolutions}), then \begin{equation*}
	E[\sldlt u](t) \le  E[\slnb \sldlt u](t) \le E[\sldlt^2 u](t) \; .
	\end{equation*}
	Thus, replacing $E[\slnb u](0)$ by $E[\sldlt^2 u](0)$ in the above estimate and upon adding positive terms we have for $t\ge a$,
	\begin{equation}\label{t^1/4bound}
	E_\mathcal{C}[u](t) \le C(1+t^{\frac{1}{4}})\left(E_\mathcal{C}[u](0)+ E[\sldlt^2 u](0)\right) \; .
	\end{equation}
	Doing the same again but using this estimate gives the uniform bound, that is, substituting for $E_\ell[\slnb u](t)$ by (\ref{estimateonEell}) in (\ref{firstestimate}) then using (\ref{t^1/4bound}) this time yields,
	\begin{equation*}
	E_\mathcal{C}[u](t) \le C\left(E_\mathcal{C}[u](0)+ E[\sldlt^2 u](0)\right) \; ,
	\end{equation*}
	for $t\ge a$. One can use the exponential bound to get the estimate over $0\le t \le a$. So, using (\ref{GronwallEc}):
	\begin{equation*}
	E_\mathcal{C}[u](t) \le E_\mathcal{C}[u](0) e^{Ct} \le C E_\mathcal{C}[u](0)\; .
	\end{equation*}
	This proves (\ref{unifomboundenergyconformal}), and in doing so, as a matter of fact, we have essentially proved (\ref{unifomboundtrappingterm}) also. If $\chi$ is any function of $r_*$ with a compact support, and with $b\ge 1$ depending on $\chi$ and playing the role of $a$, then, similar to the case of $\chi_{trap}$, one has the corresponding version of (\ref{unifomboundtrappingtermproof21}). From that, and after applying (\ref{EgammalesstElocal}), we obtain for $t' \ge b \ge 1$ that
	\begin{equation}\label{firstestimatetrappingterm}
	\int\limits_{[b,t']\times \R \times \mathcal{S}^2} t\chi_{trap}|\slnb u|^2 \d x \le C  \left( \sup\limits_{t\in[0,t']}(t E_\ell [\slnb u](t))+ \int_{0}^{t'} E_\ell[\slnb u](t) \d t \right)\; .
	\end{equation}
	Then the fact that $E_\ell[\slnb u](t)\le C E[\slnb u](t)=C E[\slnb u](0)$ gives a linear bound on the integral over the interval $[b,t']$. As before, using (\ref{estimateonEell})we substitute for $E_\ell[\slnb u](t)$ in (\ref{firstestimatetrappingterm}),  after which we use the linear bound to obtain a $t^{1/4}$ bound. Then repeat the same thing again using the $t^{1/4}$ bound we get the uniform bound over $[b,t']$. Finally, an inequality similar to the second part of (\ref{exponentialboundontrappingintergal}) holds true for $\chi$ over $[t_1,t_2]=[0,b]$. From it and from the uniform bound (or any bound) on the conformal energy, the uniform bound follows.
\end{proof}

\section{Decay of the Maxwell Field}\label{Decaysection}

We start this section by introducing some useful notations. Consider the following sets of smooth vector fields:
\begin{eqnarray}\label{vectorfields}
\mathbb{O}=\{\Theta_1,\Theta_2,\Theta_3\} \; , \nonumber\\
\mathbb{T}=\{T=\dl_t\}\cup\mathbb{O}\; ,\\
\mathbb{X}=\{R=\dl_{r_*}\}\cup\mathbb{T}\; . \nonumber
\end{eqnarray}
We also use the hatted version of a letter to indicate the corresponding normalized vector field.

We define the following pointwise norms on smooth tensors fields: Let $F$ be a $(0,m)$-tensor field, and $A$ a set of vector fields. We set
\begin{equation}\label{norm0}
|F|^2_A= \sum_{Y_1,\dots,Y_m \in A} |F(Y_1,\dots,Y_m)|^2 \;
\end{equation}
and it is a norm when $A$ is a spanning set\footnote{A spanning set $A$ of a vector space is a subset of the space with the property that every vector in the space can be written as a linear combination of vectors in $A$ only}. If in addition $F$ is smooth, $n$ is a non-negative integer, and $B$ is set of vector fields, we can define the higher order quantity:
\begin{equation}\label{normN}
|F|^2_{A,n,B}=\sum_{k=0}^n \; \sum_{X_1,\dots,X_k \in B}\; \sum_{Y_1,\dots,Y_m \in A} |(\mathcal{L}_{X_k} \dots \mathcal{L}_{X_1}F)(Y_1,\dots,Y_m)|^2 \; ,
\end{equation}
where $\mathcal{L}_X$ is the Lie derivative with respect to the vector field $X$. When $A$ is a spanning set, $|F|^2_{A,n,B}$ is a norm.  

When working with inequalities and $A$ is a spanning set, we sometimes use the same notation for the uniformly equivalent norm defined by dropping the squares in (\ref{normN}).

Another norm for the Maxwell field is the spin norm defined by $|\Phi_1|+|\Phi_0|+|\Phi_{-1}|$ or the equivalent norm $(|\Phi_1|^2+|\Phi_0|^2+|\Phi_{-1}|^2)^{1/2}$ (See (\ref{spincomponentsPhis})). It is easy to see that this norm is uniformly equivalent to the norm $|F|_{\mathbb{X}}$ by noting that the frames $\{\dl_\theta , \frac{1}{\sin(\theta)}\dl_\varphi\}$ and $\mathbb{O}$ can be expressed as linear combinations of each other with bounded smooth coefficient functions. For example, the vector field
\begin{equation*}
X= \frac{1}{\sin(\theta)}\dl_\varphi= \frac{1}{\sin(\theta)}\Theta_3
\end{equation*}
can be written as
\begin{equation}\label{normalizedTheta3}
X=\frac{\rho (\theta,\varphi)}{\sin(\theta)}\Theta_3 + \frac{(1-\rho(\theta,\varphi))}{\cos(\theta)}(\cos(\varphi)\Theta_1 - \sin(\varphi)\Theta_2) \; ,
\end{equation}
with $0\leq\rho\leq1$ a smooth function on the sphere compactly supported away from the pole $(0,0,1)$ i.e. $\theta=0$.

An essential property of Maxwell's equations is that if $F$ is a solution and $X$ is a Killing vector field then $\mathcal{L}_X F$, the Lie derivative of $F$ with respect to $X$, is again a solution of the equations. To see why, we start by (\ref{Maxeq2}). Using Cartan's identity for Lie derivatives on differential forms we immediately get
\begin{equation*}
\d(\mathcal{L}_X F)=\d(i_X \d F + \d i_X F)=\d i_X \d F=\mathcal{L}_X (\d F)=0 \; ,
\end{equation*}
in other words, the $[\lie_X,\d]=0$ on forms, and this is true for any vector field $X$, not necessarily Killing. For (\ref{Maxeq1}), we can use the expression of the Hodge star in the abstract index formalism:
\begin{equation*}
(\star \alpha)_{a_{k+1}\dots a_n}=\frac{1}{k!} \alpha^{a_1\dots a_k}\omega_{a_1\dots a_n} \; .
\end{equation*}
where $\alpha$ is a differential $k$-form and $\omega$ is the $n$-volume form given by the metric. When $X$ is Killing, $[\lie_X, \star]=0$ since the divergence of a Killing vector field vanishes by (\ref{killingvanishdivergence}). To see that this is true, we note that
\begin{equation*}
\left(\lie_X \star \alpha \right)_{a_{k+1}\dots a_n}=\frac{1}{k!} \lie_X \left(\alpha^{a_1\dots a_k}\omega_{a_1\dots a_n} \right)=\frac{1}{k!} \left(\left(\lie_X\alpha\right)^{a_1\dots a_k}\omega_{a_1\dots a_n} + \alpha^{a_1\dots a_k} \left(\lie_X\omega\right)_{a_1\dots a_n} \right) \; .
\end{equation*}
The last term is zero because of (\ref{divergenceinternsicexpression}), and the term before it is $\star \left(\lie_X \alpha\right)$, which means that the two operators commute. Hence $\lie_X F$ is a solution of (\ref{Maxeq1}) if $X$ is Killing.

We keep the assumptuion that our Maxwell field $F$ is a non-stationary solution with finite energy, i.e. satisfying (\ref{nonstationarysolutions1}) and (\ref{nonstationarysolutions}).

\subsection{Energies of the Maxwell Field }\label{Energysection}

Motivated by the Lagrangian theory, we consider an energy-momentum tensor $\mathbf{T}_{ab}$ that is a $(0,2)$-symmetric tensor i.e. $\mathbf{T}_{ab}=\mathbf{T}_{(ab)}$, and that is divergence-free i.e. $\nabla^a\mathbf{T}_{ab}=0$. Let $X$ be a vector field and  ${}^{(X)}\pi_{ab}=\nabla_a X_b + \nabla_b X_a$ be its deformation tensor. If $\mathcal{U}$ is an open submanifold of $\mathcal{N}$ with a piecewise $\mathcal{C}^1$-boundary $\dl \mathcal{U}$, then by the divergence theorem and the properties of $\mathbf{T}$, we have for $\eta$ a normal vector to $\dl \mathcal{U}$ and $\tau$ a transverse one such that $\eta^a \tau_a=1$:
\begin{eqnarray}
\int_{\dl \mathcal{U}} \mathbf{T}_{ab}X^b \eta^a i_{\tau} \d^4x &=& \int_{\mathcal{U}} \nabla^a\left(\mathbf{T}_{ab}X^b\right)  \d^4x \nonumber \\
&=&  \int_{\mathcal{U}} \left(X^b\nabla^a\mathbf{T}_{ab} + \mathbf{T}_{ab}\nabla^a X^b \right)  \d^4x \nonumber \\
&=&\hf\int_{\mathcal{U}} {}^{(X)}\pi_{ab} \mathbf{T}^{ab} \d^4x \; . \label{divergencetheoremSEtensor}
\end{eqnarray}
This is exactly what we did in the proof of Lemma \ref{conservationofenergy}, and we saw that it is particularly interesting when $X$ is Killing and thus its deformation tensor vanishes. Motivated by this, we define the energy of a general $2$-form $F$ which the energy-momentum tensor depends on (aside from the metric), on an oriented smooth hypersurface $S$ to be:
\begin{equation}\label{enegrysurfacegeneralform}
E_X [F](S)=\int_S (X\hook \mathbf{T})^\sharp \hook \d^4 x \; .
\end{equation}
Of course, it is understood that we are not integrating the $3$-form but its restriction on $S$, which is the pull back of the form by the inclusion map. We can choose $\eta_S$ and $\tau_S$ to be respectively vector fields normal and  transverse to $S$, such that their scalar product is one, then 
we have
\begin{equation}\label{energyonasurface}
E_X [F](S)=\int_S \mathbf{T}_{ab}X^b \eta_S^a i_{\tau_S} \d^4x \; .
\end{equation}
The definition is independent of the choice of $\eta_S$ and $\tau_S$. As we can see, this will lead to a conservation law when $X$ is Killing. There is a more physical and a more ``natural'' motivation for this quantity to be called energy. In fact, if $S$ is a spacelike hypersurface and $X=\dl_t$, then (\ref{energyonasurface}) is indeed the energy measured at an instant of time by an observer whose frame of reference is defined by the integral curves of the vector field $X$. For our purpose, we shall mainly consider spacelike or null slices  defined as level-hypersurfaces of a smooth function.

We know that if the Maxwell field is given as an exterior derivative of some $1$-form, one can then define a Lagrangian, and by varying the $1$-form, the Euler-Lagrange equations will be Maxwell's equations. Using this Lagrangian, it is possible to define an energy-momentum tensor which by the Euler-Lagrange equations is divergence-free. However, since in general not all Maxwell fields admit a global potential, we shall take the same energy-momentum tensor and show by direct calculations that it is divergence-free using Maxwell field equations.
Let $F$ be a Maxwell field and consider the $(0,2)$-symmetric tensor
\begin{equation}\label{energymomentumtensormaxwell}
\mathbf{T}_{ab}=\frac{1}{4} g_{ab}F^{cd}F_{cd}-F_{ac}{F_b}^c  \; .
\end{equation}

In the following lemma we summarize some of the properties of this tensor that are most important to us in this work.

\begin{lem1}\label{propertiesofenergymomentumtensorMaxwell}
	Let $\mathbf{T}$ be defined as in (\ref{energymomentumtensormaxwell}), and recall the definition of the vectors $L$ and $N$ from (\ref{nulltetrad}). We have:
	\begin{eqnarray}
	\mathbf{T}_{ab}=\mathbf{T}_{(ab)}\; , \label{symmetryEMtensor}
	g^{ab}\mathbf{T}_{ab}&=& 0 \; , \label{MaxwellEMtensortracefree}\\
	\nabla^a \mathbf{T}_{ab} &=& 0 \; , \label{MaxwellEMtensordivergencefree}\\
	\mathbf{T}_{ab} L^a L^b &=& r^{-2}|\Phi_1|^2 \; , \label{TLL}\\
	\mathbf{T}_{ab} L^a N^b &=& f r^{-4}|\Phi_0|^2 \; , \label{TLN}\\
	\mathbf{T}_{ab} N^a N^b &=& r^{-2}|\Phi_{-1}|^2 \; .\label{TNN}
	\end{eqnarray}
\end{lem1}
\begin{proof}
	The trace-freeness is immediate:
	\begin{equation*}
	g^{ab} \mathbf{T}_{ab}={\mathbf{T}^a}_a= \frac{4}{4}F^{cd}F_{cd}-F_{ac}F^{ac}=0 \; .
	\end{equation*}
	For the divergence, from (\ref{Maxeqabst1}) we have,
	\begin{eqnarray*}
		\nabla^a \mathbf{T}_{ab} &=& \hf g_{ab}F_{cd}\nabla^a F^{cd} - F_b^c \nabla^a F_{ac} - F_{ac}\nabla^a F_b^c\\
		&=&\hf F^{ac}\nabla_b F_{ac} -F^{ac}\nabla_a F_{bc} \; .
	\end{eqnarray*}
	And from (\ref{Maxeqabst2}),
	\begin{equation*}
	\nabla_b F_{ac}= \nabla_c F_{ab} + \nabla_a F_{bc} \; ,
	\end{equation*}
	which entails
	\begin{equation*}
	\nabla^a \mathbf{T}_{ab}= \hf F^{ac}\nabla_c F_{ab} - \hf F^{ac}\nabla_a F_{bc}=0
	\end{equation*}
	upon swapping the indices $a$ and $c$ in the first term, then using the fact that $F$ is antisymmetric.
	
	For the last three properties (\ref{TLL})-(\ref{TNN}), we have,
	\begin{eqnarray}
	|\Phi_1|^2 &=& \left(F_{02}+F_{12} \right)^2 + \frac{1}{\sin(\theta)^2}\left(F_{03}+F_{13}\right)^2 \; , \label{absolutevaluePhi1} \\
	|\Phi_0|^2 &=& \frac{r^4}{f^2} (F_{01})^2 +  \frac{1}{\sin(\theta)^2} (F_{23})^2 \; , \label{absolutevaluePhi0} \\
	|\Phi_{-1}|^2 &=&  \left(F_{02}-F_{12} \right)^2 + \frac{1}{\sin(\theta)^2}\left(F_{03}-F_{13}\right)^2 \; . \label{absolutevaluePhi-1}
	\end{eqnarray}
	Since both $L$ and $N$ are null, then for $X \in \{L,N \}$
	\begin{equation*}
	\mathbf{T}_{ab} X^a X^b= -g^{cd} F_{ac} F_{bd} X^a X^b \; ,
	\end{equation*}
	which in a local coordinate basis expands to
	\begin{equation*}
	\mathbf{T}(X,X)=-\frac{1}{f}\left( F_{\mathbf{a}0}X^\mathbf{a}\right)^2 + \frac{1}{f}\left( F_{\mathbf{a}1}X^\mathbf{a}\right)^2 + \frac{1}{r^2}\left( F_{\mathbf{a}2}X^\mathbf{a}\right)^2 + \frac{1}{r^2 \sin(\theta)^2}\left( F_{\mathbf{a}3}X^\mathbf{a}\right)^2 \; ,
	\end{equation*}
	since $F_{00}=F_{11}=0$, the first two terms cancel out when we evaluate $L$ and $N$ with their coordinate expressions. This proves (\ref{TLL}) and (\ref{TNN}). For the middle equation, we need to calculate $\mathbf{T}(L,N)$. We have
	\begin{equation*}
	\mathbf{T}(L,N)=\mathbf{T}_{00}- \mathbf{T}_{11} \; ,
	\end{equation*}
	but,
	\begin{equation*}
	\mathbf{T}_{00} = \frac{1}{4} f F_{\mathbf{cd}}F^{\mathbf{cd}} - F_{0\mathbf{c}}{F_0}^{\mathbf{c}}= \frac{1}{4} f F_{\mathbf{cd}}F^{\mathbf{cd}} -f F_{0\mathbf{c}} F^{0\mathbf{c}} \; ,
	\end{equation*}
	similarly,
	\begin{equation*}
	\mathbf{T}_{11}= - \frac{1}{4} f F_{\mathbf{cd}}F^{\mathbf{cd}} + f F_{1\mathbf{c}} F^{1\mathbf{c}} \; ,
	\end{equation*}
	thus,
	\begin{eqnarray*}
		\mathbf{T}_{00}- \mathbf{T}_{11} &=& f \left( \hf F_{\mathbf{cd}}F^{\mathbf{cd}} - F_{0\mathbf{c}} F^{0\mathbf{c}} - F_{1\mathbf{c}} F^{1\mathbf{c}} \right) \\
		&=&f \left( \sum_{0 \le c<d \le 3} F_{\mathbf{cd}}F^{\mathbf{cd}} - \sum_{c} F_{0\mathbf{c}} F^{0\mathbf{c}} - \sum_{c} F_{1\mathbf{c}} F^{1\mathbf{c}} \right) \\
		&=&  f \left( \sum_{1 \le c<d \le 3} F_{\mathbf{cd}}F^{\mathbf{cd}} - \sum_{c} F_{1\mathbf{c}} F^{1\mathbf{c}} \right) \\
		&=& f\left( \sum_{2 \le c<d \le 3} F_{\mathbf{cd}}F^{\mathbf{cd}} - F_{10} F^{10} \right)\\
		&=& f \left(  F_{23}F^{23} - F_{10} F^{10} \right)= \frac{1}{f} (F_{01})^2 +\frac{f}{r^4\sin(\theta)^2}(F_{23})^2 \; ,
	\end{eqnarray*}
	which proves (\ref{TLN}).
\end{proof}

Identities (\ref{TLL}), (\ref{TLN}), and (\ref{TNN}) have two important consequences. First, let $V$ be a future oriented causal vector with no angular components, i.e. of the form
\begin{equation}\label{futuretimelikevectornoangular}
V=\alpha \dl_t + \beta \dlr \quad \mathrm{with} \quad \alpha \ge |\beta| \; .
\end{equation}
Since any such vector can be expressed as a linear combination of the vectors $L$ and $N$ with non-negative coefficients, more precisely,
\begin{equation*}
V=\left(\frac{\alpha + \beta}{2}\right) L + \left(\frac{\alpha - \beta}{2}\right) N \; ,
\end{equation*}
the next corollary is immediate.
\begin{cor1}[Dominant Energy Condition]\label{DEC}
	Let $\mathbf{T}$ be the stress-energy tensor of (\ref{energymomentumtensormaxwell}), and let $V$ and $W$ be two future (past) oriented causal vectors with no angular components then,
	\begin{equation}
	\mathbf{T}(V,W)\ge 0 \; .
	\end{equation}
\end{cor1}
The energy-momentum tensor of the Maxwell field satisfies a stronger positivity condition. Actually, the corollary holds true for all future oriented causal vectors, even with non zero angular components. This is called the dominant energy condition. The proof of the full dominant energy condition is much easier to see using spinor notations (see \cite{penrose_spinors_1987,penrose_spinors_1988}).

The second consequence of these identities is related to the definition of the energy of the Maxwell field. We can expand the left hand side of (\ref{TLL})-(\ref{TNN}) to get
\begin{eqnarray}
r^{-2}|\Phi_1|^2 &=& \mathbf{T}_{00}+ 2 \mathbf{T}_{01} + \mathbf{T}_{11} \; , \label{PhisinT1} \\
f r^{-4}|\Phi_0|^2  &=& \mathbf{T}_{00}- \mathbf{T}_{11} \; ,  \label{PhisinT2}  \\
r^{-2}|\Phi_{-1}|^2 &=& \mathbf{T}_{00}- 2 \mathbf{T}_{01} + \mathbf{T}_{11}  \; ,  \label{PhisinT3}
\end{eqnarray}
and from this, we can compute the components of the stress-energy tensor in terms of the spin components,
\begin{eqnarray}
\mathbf{T}_{00} &=& \frac{1}{4r^2}\left(|\Phi_1|^2 + \frac{2f}{r^2}|\Phi_0|^2 + |\Phi_{-1}|^2\right)\label{T00} \; , \\
\mathbf{T}_{01} &=& \frac{1}{4r^2}\left(|\Phi_1|^2 - |\Phi_{-1}|^2\right)\label{T01}  \; , \\
\mathbf{T}_{11} &=& \frac{1}{4r^2}\left(|\Phi_1|^2 - \frac{2f}{r^2}|\Phi_0|^2 + |\Phi_{-1}|^2\right)\label{T11} \; .
\end{eqnarray}
If $S=\Sigma_t=\{t\}\times\R\times\mathcal{S}^2$, then its unit normal is $\hat{T}=f^{-\hf}\dl_t$, and taking the transverse vector to be $\hat{T}$ also, we see that
\begin{eqnarray}
E_T[F](\Sigma_t)&=&\int_{\Sigma_t} \mathbf{T}_{ab}T^b \hat{T}^a i_{\hat{T}} \d^4x \nonumber \\
&=&\int_{\Sigma_t} \mathbf{T}_{00} f^{-\hf} f^{\hf}r^2 \d_{r_*} \d^2\omega \nonumber\\
&=&\frac{1}{4} \int_{\Sigma_t} |\Phi_1|^2 + \frac{2f}{r^2}|\Phi_0|^2 + |\Phi_{-1}|^2 \d r_* \d^2 \omega=E_T[F](t) \; .\label{energyofmaxwellont=cst}
\end{eqnarray}
This gives, 
from (\ref{divergencetheoremSEtensor}) in the introductory argument of this section, that it is a conserved quantity as the vector field $T$ is Killing, i.e.
\begin{equation}\label{conservationofenergyMAXWELL}
E_T[F](t)=E_T[F](0) \; .
\end{equation}
And since the Lie derivative of a Maxwell field $F$ with respect to a Killing vector field $X$ is again a Maxwell field, its energy is conserved as well. So, for $X_1, \dots , X_k$ Killing,
\begin{equation}\label{conservationofenergyLIEDERIVMAXWELL}
E_T[\lie_{X_1}\dots\lie_{X_k} F](t)=E_T[\lie_{X_1}\dots\lie_{X_k} F](0) \; .
\end{equation}

To get decay for the Maxwell field we will need to control another energy on $\Sigma_t$ defined by the vector field
\begin{equation}\label{confromalvectortr}
K=(t^2+r_*^2)\dl_t + 2tr_* \dlr \; .
\end{equation}
Using the advanced and retarded coordinates $u_+=t+r_*$ and $u_- = t-r_*$ it becomes
\begin{equation}\label{confromalvector}
K=\hf(u_+^2 L + u_-^2 N)\; .
\end{equation}
To compute $E_K[F](t):=E_K[F](\Sigma_t)$, which we shall call the conformal energy, we have to compute $\mathbf{T}_{0\mathbf{a}}K^{\mathbf{a}}$. From the above form of $K$, (\ref{T00}), and (\ref{T01}) we have
\begin{eqnarray*}
	\mathbf{T}_{0\mathbf{a}}K^{\mathbf{a}} &=& \hf\left(u_+^2 \mathbf{T}_{0\mathbf{a}}L^{\mathbf{a}}+ u_-^2 \mathbf{T}_{0\mathbf{a}}N^{\mathbf{a}} \right) \\
	&=& \hf\left(u_+^2(\mathbf{T}_{00}+\mathbf{T}_{01})+u_-^2(\mathbf{T}_{00}-\mathbf{T}_{01})\right) \\
	&=& \frac{1}{8r^2}\left(u_+^2\left(2|\Phi_1|^2 +  \frac{2f}{r^2}|\Phi_0|^2\right) + u_-^2\left(2|\Phi_{-1}|^2 +  \frac{2f}{r^2}|\Phi_0|^2\right) \right) \\
	&=&  \frac{1}{4r^2}\left(u_+^2|\Phi_1|^2 +  (u_+^2+u_-^2)\frac{f}{r^2}|\Phi_0|^2 + u_-^2|\Phi_{-1}|^2 \right) \; .
\end{eqnarray*}
Therefore, the conformal energy is
\begin{eqnarray}
E_K[F](t)&=&\int_{\Sigma_t} \mathbf{T}_{ab}K^b \hat{T}^a i_{\hat{T}} \d^4x \nonumber \\
&=&\frac{1}{4}\int_{\Sigma_t} u_+^2|\Phi_1|^2 +  (u_+^2+u_-^2)\frac{f}{r^2}|\Phi_0|^2 + u_-^2|\Phi_{-1}|^2 \d r_* \d^2 \omega \label{conformalenergy} \; .
\end{eqnarray}

The next result is rather important and central to our work, it is sometimes referred to as the Trapping Effect. From it and its proof, a great part of the importance of the photon sphere, the Regge-Wheeler coordinate, and the wave analysis we did in section \ref{sectionwaveanalysis}, are revealed.

\begin{lem1}[Trapping Lemma]\label{trappinglemma}
	There is a non-negative compactly supported function $\chi_{trap}$ of $r_*$, depending only on the geometry of the spacetime, such that for $F$ a solution to Maxwell's equations, we have
	\begin{equation}\label{trappinginequality}
	E_K[F](t_2)-E_K[F](t_1) \le \int\limits_{[t_1,t_2]\times \R\times \mathcal{S}^2} t \chi_{trap} |\Phi_0|^2  \d t \d r_* \d^2 \omega \; .
	\end{equation}
\end{lem1}
\begin{proof}
	This is a matter of applying the divergence theorem (\ref{divergencetheoremSEtensor}). For this we need to calculate the deformation tensor of $K$. Since $T^\flat=f \d t$ and $R^\flat=-f \d r_*$, then from (\ref{confromalvectortr}),
	\begin{equation*}
	K_b=(t^2 +r_*^2)T_b + 2tr_* R_b \; ,
	\end{equation*}
	so,
	\begin{eqnarray*}
		\nabla_a K_b \! \!\!&=&\!\!\! T_b \nabla_a (t^2+r_*^2) + (t^2 + r_*^2)\nabla_a T_b  + R_b\nabla_a(2tr_*) +2tr_*\nabla_a R_b \\
		\!\!\! &=& \!\!\! T_b\left(2tf^{-1}T_a -2r_* f^{-1}R_a \right) + (t^2 + r_*^2)\nabla_a T_b + R_b\left(2r_* f^{-1}T_a - 2tf^{-1}R_a \right) + 2tr_*\nabla_a R_b\\
		\!\!\! &=&\!\!\! 2f^{-1}t\left( T_a T_b - R_a R_b \right)+2tr_*\nabla_a R_b + 2f^{-1}r_* \left(T_a R_b - T_b R_a\right) + (t^2 + r_*^2)\nabla_a T_b \; .
	\end{eqnarray*}
	Symmetrizing, the third term vanishes since it is skew, and the last term also vanishes as $T$ is Killing. Thus, the only term that we need to calculate is the covariant derivative of $R_a$. By (\ref{Crstflsymbinr*}) we have,
	\begin{eqnarray*}
		\nabla R^\flat &=& \left(\dl_\mathbf{a} R_\mathbf{b} - \tilde{\Gamma}^\mathbf{c}_{\mathbf{ab}}R_\mathbf{c} \right)\d x^\mathbf{a} \otimes \d x^\mathbf{b} \\
		&=& \dlr R_1 \d r_*^2 - \sum_a \tilde{\Gamma}^1_{\mathbf{aa}}R_1(\d x^\mathbf{a})^2\\
		&=&  \frac{ff'}{2}\d t^2 - \frac{ff'}{2}\d r_*^2 +\frac{f}{r}\d \omega^2  \\
		&=&  \frac{f}{r} g_\mathcal{M} +  f\left(\frac{f'}{2} - \frac{f}{r}\right)\left(\d t^2 - \d r_*^2\right) \; .
	\end{eqnarray*}
	Now, since this is a symmetric tensor, $\nabla_a R_b =\nabla_{(a} R_{b)}$. Therefore,
	\begin{eqnarray*}
		{}^K \pi_{ab} =2 \nabla_{(a} K_{b)} \!\!\!&=& \!\!\! 4f^{-1}t\left( T_a T_b - R_a R_b \right) +4tr_*\nabla_a R_b \\
		\!\!\! &=&\!\!\! 4f^{-1}t\left( T_a T_b - R_a R_b \right) + 4tr_*fr^{-1} g_\mathcal{M} +  4tr_*f^{-1}\left(\frac{f'}{2} - \frac{f}{r}\right)\left( T_a T_b - R_a R_b\right) \\
		&=& 4tr_*fr^{-1} g_{ab} + 4f^{-1}t\left(1+r_* \left(\frac{f'}{2} - \frac{f}{r}\right)\right)\left( T_a T_b - R_a R_b \right) \; .
	\end{eqnarray*}
	Recall the trapping term defined in (\ref{trappingterm}),
	\begin{equation*}
	\mathscr{T}= 1+ r_*\left(\frac{f'}{2} - \frac{f}{r}\right) \; .
	\end{equation*}
	Contracting the deformation tensor of $K$ with the energy-momentum tensor of the Maxwell field $\mathbf{T}$, and using (\ref{MaxwellEMtensortracefree}) and (\ref{PhisinT2}) we obtain
	\begin{equation*}
	{}^K \pi_{ab}\mathbf{T}^{ab}=4tr^{-4}\mathscr{T}|\Phi_0|^2 \; .
	\end{equation*}
	In virtue of the definition of the conformal energy given by (\ref{conformalenergy}), and through (\ref{divergencetheoremSEtensor}) we arrive at
	\begin{equation*}
	E_K[F](t_2)-E_K[F](t_1) = \int\limits_{[t_1,t_2]\times \R\times \mathcal{S}^2} 2t\frac{f}{r^2} \mathscr{T}  |\Phi_0|^2  \d t \d r_* \d^2 \omega \; .
	\end{equation*}
	Finally, our choice of $\chi_{trap}$ at the end of the proof of Lemma \ref{controlonconformalchargeTHEOREM} ensures that the statement of the lemma holds true.
\end{proof}

The trapping lemma shows that the conformal energy of the Maxwell field is controlled by the integral of the middle component. This is why we did the wave analysis in section \ref{sectionwaveanalysis}. As the middle component satisfies the wave-like equation (\ref{WavePhi0}), then the bounds we obtained in Proposition \ref{unifomboundTHEOREM} can be used to establish bounds on the conformal energy. For this purpose, we need to find a relation between the energy of the middle term and the full Maxwell field. Two $1$-forms are particularly important for the discussion. We set
\begin{equation}\label{alphas}
\alpha=i_L F \; ,   \hspace{2cm} \textrm{and} \hspace{2cm} \underline{\alpha}=i_N F \; .
\end{equation}
We see that
\begin{equation*}
\Phi_1=\alpha(M)=\alpha_{\theta} + \frac{i}{\sin(\theta)}\alpha_{\varphi} \; ,
\end{equation*}
and
\begin{equation*}
\Phi_{-1}=\underline{\alpha}(\bar{M})=\underline{\alpha}_{\theta} - \frac{i}{\sin(\theta)}\underline{\alpha}_{\varphi} \; .
\end{equation*}
Using this, a simple calculation shows that
\begin{eqnarray*}
	(\bar{M}+\cot(\theta))\Phi_1 &=& \dl_\theta \alpha_\theta +\frac{1}{\sin(\theta)^2}\dl_\varphi \alpha_\varphi + \cot(\theta)\alpha_\theta + \frac{i}{\sin(\theta)}(\dl_\theta \alpha_\varphi - \dl_\varphi \alpha_\theta) \; , \\
	(M+\cot(\theta))\Phi_{-1} &=& \dl_\theta \underline{\alpha}_\theta +\frac{1}{\sin(\theta)^2}\dl_\varphi \underline{\alpha}_\varphi + \cot(\theta)\underline{\alpha}_\theta - \frac{i}{\sin(\theta)}(\dl_\theta \underline{\alpha}_\varphi - \dl_\varphi \underline{\alpha}_\theta) \; .
\end{eqnarray*}
Since we are going to work only with operations defined on the sphere $\mathcal{S}^2$ when dealing with these $1$-forms, then the $t$ and $r_*$ variables play no role in our calculations. So one might replace these forms with their projections on the sphere of radius $r_*$. The $r_*$ will only be a scaling factor. Moreover, from (\ref{alphas}) and then (\ref{absolutevaluePhi1}) and (\ref{absolutevaluePhi-1}),
\begin{eqnarray*}
	(\alpha_\theta)^2  + \frac{1}{\sin(\theta)^2}(\alpha_\varphi)^2 &=& \left(F_{02}+F_{12} \right)^2 + \frac{1}{\sin(\theta)^2}\left(F_{03}+F_{13}\right)^2 = |\Phi_1|^2  \; , \\
	(\underline{\alpha}_\theta)^2  + \frac{1}{\sin(\theta)^2}(\underline{\alpha}_\varphi)^2 &=&  \left(F_{02}-F_{12} \right)^2 + \frac{1}{\sin(\theta)^2}\left(F_{03}-F_{13}\right)^2 = |\Phi_{-1}|^2  \; ,
\end{eqnarray*}
and we can consider them as covector fields on the sphere $\mathcal{S}^2$. Therefore we set
\begin{equation*}
\alpha = \alpha_\theta \d \theta + \alpha_\varphi \d \varphi \; , \hspace{0.5cm} \underline{\alpha} =  \underline{\alpha}_\theta \d \theta + \underline{\alpha}_\varphi \d \varphi \; ,
\end{equation*}
\begin{equation*}
|\alpha|^2:=g_{\mathcal{S}^2}(\alpha,\alpha) \; , \hspace{0.5cm} |\underline{\alpha}|^2:=g_{\mathcal{S}^2}(\underline{\alpha},\underline{\alpha})\; ,
\end{equation*}
and so
\begin{equation}\label{alphaandPhi}
|\alpha|^2=|\Phi_1|^2  \; , \hspace{0.5cm} |\underline{\alpha}|^2=|\Phi_{-1}|^2 \; .
\end{equation}
Taking this interpretation of the $1$-forms, we now see that
\begin{eqnarray}
(\bar{M}+\cot(\theta))\Phi_1 &=& \slnb^a \alpha_a + i \epsilon^{ab} \slnb_a \alpha_b \; ,\label{geomtricinterpretation1} \\
(M+\cot(\theta))\Phi_{-1} &=&  \slnb^a \underline{\alpha}_a - i \epsilon^{ab} \slnb_a \underline{\alpha}_b \label{geomtricinterpretation2} \; ,
\end{eqnarray}
where $\epsilon^{ab}$ is the Levi-Civita tensor on $\mathcal{S}^2$ defined by
\begin{equation*}
\epsilon^{\mathbf{ab}}=\frac{1}{\sqrt{|g_{\mathcal{S}^2}|}} \varepsilon^{\mathbf{ab}}=\frac{1}{\sin(\theta)} \varepsilon^{\mathbf{ab}}
\end{equation*}
and $\varepsilon^{ab}$ is the Levi-Civita symbol given in (\ref{Levi-Civitasymbol}).

With these notations, we can now relate the energies of $\Phi_0$ and $F$.
\begin{lem1}\label{EPhi0andEFlemma}
	Let $F$ be a Maxwell field and let $(\Phi_1,\Phi_0,\Phi_{-1})$ be its spin components and consider the energy and the conformal charge defined in (\ref{energywavelike}) and (\ref{conformalchargewave}). Then
	\begin{eqnarray}
	E[\Phi_0](t) &\le& C \sum_{i=1}^3 E_T[\lie_{\Theta_i}F](t) \; , \label{EPhi0andEF1}\\
	E[\Delta_{\mathcal{S}^2}^2 \Phi_0](t) &\le& C \sum_{X_1, \dots, X_5\in \mathbb{O}}E_T[\lie_{X_1} \dots \lie_{X_5}F](t) \; ,  \label{EPhi0andEF2} \\
	E_{\mathcal{C}}[\Phi_0](t) &\le& C \sum_{i=1}^3\left( E_K[\lie_{\Theta_i}F](t) + E_T[\lie_{\Theta_i}F](t)\right) \; . \label{EPhi0andEF3}
	\end{eqnarray}
\end{lem1}

\begin{proof}
	Let $\Sigma_t=\{t\} \times \R \times \mathcal{S}^2$. By (\ref{alphaandPhi}) we have
	\begin{equation*}
	E_T[F](t)=\frac{1}{4} \int_{\Sigma_t} |\alpha|^2 + 2V|\Phi_0|^2 + |\underline{\alpha}|^2 \d r_* \d^2 \omega \; .
	\end{equation*}
	Thus, to find $E_T[\lie_{\Theta_i}F](t)$ we need to calculate the middle spin component of $\lie_{\Theta_i}F$ and the corresponding $1$-forms. For a $2$-form $\beta$ and two vectors $X$ and $Y$ such that $[X,Y]=0$, we have for all vectors $Z$
	\begin{eqnarray*}
		(\lie_Y i_X \beta) (Z)&=& \lie_Y \left((i_X \beta) (Z) \right) -(i_X \beta) (\lie_Y Z) \\
		&=& \lie_Y\left( \beta (X,Z) \right) - \beta \left(X , \lie_Y Z \right)  \\
		&=& (\lie_Y \beta)(X,Z) + \beta(\lie_Y X , Z) + \beta(X,\lie_Y Z) - \beta \left(X , \lie_Y Z \right) \\
		&=&(i_X \lie_Y \beta)(Z) \; .
	\end{eqnarray*}
	Thus,
	\begin{equation*}
	i_L (\lie_{\Theta_i} F )= \lie_{\Theta_i}(i_L F) = \lie_{\Theta_i} \alpha ,
	\end{equation*}
	and the same for $\underline{\alpha}$. This treats the forms, for the middle component we need the commutators of $M$ and $\bar{M}$ with the $\Theta_i$'s. Direct calculation gives,
	\begin{gather*}
	[\Theta_1,M]=\frac{-i\cos(\varphi)}{\sin(\theta)}M \quad ; \quad [\Theta_1,\bar{M}]=\frac{i\cos(\varphi)}{\sin(\theta)}\bar{M} \\
	[\Theta_2,M]=\frac{i\sin(\varphi)}{\sin(\theta)}M \quad ; \quad [\Theta_2,\bar{M}]=\frac{-i\sin(\varphi)}{\sin(\theta)}\bar{M} \\
	[\Theta_3,M]=0 \quad ; \quad [\Theta_3,\bar{M}]=0 \; .\\
	\end{gather*}
	Hence, the middle component of $\lie_{\Theta_i} F$ is
	\begin{eqnarray*}
		\frac{1}{2}\left(V^{-1}(\lie_{\Theta_i}F)(L,N)+(\lie_{\Theta_i}F)\left(\bar{M},M\right)\right) &=& \hf \Theta_i\left(V^{-1} F(L,N) + F \left(\bar{M},M \right)\right)\\
		& & \hspace{1cm} - \hf \left( F \left(\lie_{\Theta_i}\bar{M},M \right) + F \left(\bar{M},\lie_{\Theta_i}M \right) \right)\\
		&=&\hf \Theta_i\left(V^{-1} F(L,N) + F \left(\bar{M},M \right)\right)\\
		&=& \Theta_i (\Phi_0) \; .
	\end{eqnarray*}
	Therefore, and using (\ref{rotationformofslnb}),
	\begin{equation}\label{energylieOF}
	\sum_{i=1}^3 E_T[\lie_{\Theta_i}F](t)=\frac{1}{4} \int_{\Sigma_t}  \sum_{i=1}^3 \left(|\lie_{\Theta_i}\alpha|^2 +|\lie_{\Theta_i}\underline{\alpha}|^2\right)  + 2V|\slnb\Phi_0|^2 \d r_* \d^2 \omega \; .
	\end{equation}
	
	Next we want to substitute for the $t$ and $r_*$ derivatives of $\Phi_0$ in $E[\Phi_0](t)$ in terms of derivatives of the $1$-forms $\alpha$ and $\underline{\alpha}$. From (\ref{NewmanPenrose2}) and (\ref{NewmanPenrose3}) we have
	\begin{equation*}
	|L\Phi_0|^2 + |N\Phi_0|^2 = |\bar{M}\Phi_1 + \cot(\theta)\Phi_1|^2 + |M\Phi_{-1} + \cot(\theta)\Phi_{-1}|^2 \; .
	\end{equation*}
	On the one hand, by the parallelogram rule,
	\begin{equation*}
	|L\Phi_0|^2 + |N\Phi_0|^2 = 2|\dl_t \Phi_0|^2 + 2|\dlr \Phi_0|^2 \; ,
	\end{equation*}
	and on the other hand, by (\ref{geomtricinterpretation1}) and (\ref{geomtricinterpretation2}),
	\begin{equation*}
	|(\bar{M}1 + \cot(\theta))\Phi_1|^2 + |(M + \cot(\theta))\Phi_{-1}|^2=  (\slnb^a \alpha_a)^2 + (\epsilon^{ab} \slnb_a \alpha_b)^2 +  (\slnb^a \underline{\alpha}_a)^2 + (\epsilon^{ab} \slnb_a \underline{\alpha}_b)^2 \; .
	\end{equation*}
	Thus,
	\begin{equation*}
	E[\Phi_0](t)=\frac{1}{4} \int_{\Sigma_t}   (\slnb^a \alpha_a)^2 + (\epsilon^{ab} \slnb_a \alpha_b)^2 +  (\slnb^a \underline{\alpha}_a)^2 + (\epsilon^{ab} \slnb_a \underline{\alpha}_b)^2 + 2V|\slnb\Phi_0|^2 \d r_* \d^2 \omega \; .
	\end{equation*}
	So, (\ref{EPhi0andEF1}) follows if we show that
	\begin{eqnarray}
	(\slnb^a \alpha_a)^2 + (\epsilon^{ab} \slnb_a \alpha_b)^2 &\le& C  \sum_{i=1}^3 |\lie_{\Theta_i}\alpha|^2 \; , \nonumber \\
	&& \label{liegreaterdivpluscurl} \\
	(\slnb^a \underline{\alpha}_a)^2 + (\epsilon^{ab} \slnb_a \underline{\alpha}_b)^2 &\le& C  \sum_{i=1}^3 |\lie_{\Theta_i}\underline{\alpha}|^2 \; . \nonumber
	\end{eqnarray}
	We will workout the $\alpha$'s case only since the case of $\underline{\alpha}$ is completely analogous. We start by showing that
	\begin{equation}\label{sumliealpha}
	\sum_{i=1}^3 |\lie_{\Theta_i}\alpha|^2 = \slnb_a \alpha_b \slnb^a \alpha^b + |\alpha|^2 \; .
	\end{equation}
	We have
	\begin{eqnarray*}
		\sum_{i=1}^3 |\lie_{\Theta_i}\alpha|^2 &=& \sum_{i=1}^3 \left(\Theta_i^b \slnb_b \alpha_a + \alpha_b \slnb_a \Theta_i^b\right)\left(\Theta_i^c \slnb_c \alpha^a + \alpha_c \slnb^a \Theta_i^c\right) \\
		&=&  \slnb_b \alpha_a \slnb_c \alpha^a \sum_{i=1}^3\Theta_i^b \Theta_i^c + 2 \alpha_c \slnb_b \alpha_a \sum_{i=1}^3\Theta_i^b \slnb^a \Theta_i^c + \alpha_b \alpha_c \sum_{i=1}^3 \slnb_a \Theta_i^b \slnb^a \Theta_i^c \; .
	\end{eqnarray*}
	A lengthy but very straight forward computation\footnote{It is worth mentioning that in addition to checking these identities by hand, we had also carried out these calculations using Sage which proves to be a useful symbolic computational program for differential geometry, in both calculus and algebra.} shows that
	\begin{equation*}
	\sum_{i=1}^3\Theta_i^b \Theta_i^c=g_{\mathcal{S}^2}^{bc} \quad ; \quad  \sum_{i=1}^3\Theta_i^b \slnb^a \Theta_i^c=0 \quad ; \quad \sum_{i=1}^3 \slnb_a \Theta_i^b \slnb^a \Theta_i^c=g_{\mathcal{S}^2}^{bc} \; ,
	\end{equation*}
	and (\ref{sumliealpha}) follows immediately.
	
	Furthermore, if we denote the components of $\slnb \alpha$ in the $(\omega_1=\theta,\omega_2=\varphi)$ coordinates by
	\begin{equation*}
	\alpha_{\mathbf{a} ; \mathbf{b}}:= \dl_\mathbf{b} \alpha_\mathbf{a} + {}^{\mathcal{S}^2}\Gamma^\mathbf{c}_{\mathbf{ab}} \alpha_\mathbf{c} \; ,
	\end{equation*}
	then
	\begin{eqnarray*}
		(\slnb^a \alpha_a)^2 &=& \left(g_{\mathcal{S}^2}^{11}\alpha_{1;1} +  g_{\mathcal{S}^2}^{22}\alpha_{2;2}\right)^2 \; , \\
		(\epsilon^{ab}\slnb_a \alpha_b)^2&=&  g_{\mathcal{S}^2}^{22}\left(\alpha_{1;2} - \alpha_{2;1}\right)^2 \; .
	\end{eqnarray*}
	And so, using the identities $2(a^2+b^2)\ge(a\pm b)^2$, then the fact that $g_{\mathcal{S}^2}^{11}=1$, and finally (\ref{sumliealpha}), we have
	\begin{eqnarray*}
		(\slnb^a \alpha_a)^2 + (\epsilon^{ab}\slnb_a \alpha_b)^2 &\le& 2\left(\left(g_{\mathcal{S}^2}^{11}\alpha_{1;1}\right)^2 +  \left(g_{\mathcal{S}^2}^{22}\alpha_{2;2}\right)^2 +g_{\mathcal{S}^2}^{22} (\alpha_{1;2}^2 + \alpha_{2;1}^2) \right)  \\
		&=& 2 \slnb_a \alpha_b \slnb^a \alpha^b \\
		&\le&  2 \sum_{i=1}^3 |\lie_{\Theta_i}\alpha|^2 \; ,
	\end{eqnarray*}
	which proves (\ref{EPhi0andEF1}).
	
	Similar to what we did in the beginning of the proof to obtain (\ref{energylieOF}), one has
	\begin{eqnarray*}
		\sum_{X_1, \dots, X_5\in \mathbb{O}}E_T[\lie_{X_1} \dots \lie_{X_5}F](t)&=& \frac{1}{4}\sum_{X_1, \dots, X_5\in \mathbb{O}}\int_{\Sigma_t} 2V|\lie_{X_1}\dots \lie_{X_5} \Phi_0|^2 \\
		&& \hspace{2cm}+ |\lie_{X_1}\dots \lie_{X_5} \alpha|^2 +|\lie_{X_1}\dots \lie_{X_5} \underline{\alpha}|^2 \d r_* \d^2 \omega \; .
	\end{eqnarray*}
	By (\ref{rotationformofsldlt}), we see that
	\begin{equation*}
	\sldlt^2 \Phi_0 = \sum_{i,j=1}^3 \Theta_i^2 \Theta_j^2 \Phi_0 \; .
	\end{equation*}
	if we set
	\begin{equation*}
	{}^{ij}\Phi_0=\Theta_i^2 \Theta_j^2 \Phi_0 \; ,
	\end{equation*}
	then ${}^{ij}\Phi_0$ is the middle component of $\lie_{\Theta_i}^2 \lie_{\Theta_j}^2 F$. Using the identity $n(a_1^2+\dots +a_n^2)\le(a_1+\dots a_2)^2$ we have
	\begin{equation*}
	E\left[\sldlt^2\Phi_o\right](t)=E\left[\sum_{i,j=1}^3 {}^{ij}\Phi_0\right](t)\le C \sum_{i,j=1}^3 E\left[{}^{ij}\Phi_0\right](t) \; ,
	\end{equation*}
	then from (\ref{EPhi0andEF1}) applied to ${}^{ij}\Phi_0$ we get
	\begin{equation*}
	E\left[{}^{ij}\Phi_0\right](t) \le C \sum_{k=1}^3 E_T[\lie_{\Theta_k}\lie_{\Theta_i}^2 \lie_{\Theta_j}^2F](t) \; .
	\end{equation*}
	Summing over $i$ and $j$ and noting that
	\begin{equation*}
	\sum_{i,j,k=1}^3 E_T[\lie_{\Theta_k}\lie_{\Theta_i}^2 \lie_{\Theta_j}^2F](t) \le \sum_{X_1, \dots, X_5\in \mathbb{O}}E_T[\lie_{X_1} \dots \lie_{X_5}F](t) \; ,
	\end{equation*}
	we obtain (\ref{EPhi0andEF2}).
	
	The final inequality of the lemma is obtain using the same techniques. Indeed, the conformal energy of $\lie_{\Theta_i}F$ can be written as
	\begin{equation*}
	E_K[\lie_{\Theta_i}F](t)=\frac{1}{4}\int_{\Sigma_t} u_+^2 |\lie_{\Theta_i} \alpha|^2 + u_-^2 |\lie_{\Theta_i}\underline{\alpha}|^2 +(u_+^2 +u_-^2)V|\lie_{\Theta_i}\Phi_0|^2 \d r_* \d^2 \omega \; ,
	\end{equation*}
	and on the other side, using (\ref{conformaldensitypositive}) the conformal charge of $\Phi_0$ is
	\begin{equation*}
	E_{\mathcal{C}}[\Phi_0](t)=\frac{1}{4}\int_{\Sigma_t} \hf\left(u_+^2 |L\Phi_0|^2 + u_-^2 |N\Phi_0|^2\right) +(u_+^2 +u_-^2)V|\slnb\Phi_0|^2 \d r_* \d^2 \omega +   E[\Phi_0](t) \; ,
	\end{equation*}
	but since as we saw,
	\begin{eqnarray*}
		|L\Phi_0|^2 &=&   |(\bar{M}1 + \cot(\theta))\Phi_1|^2 =  (\slnb^a \alpha_a)^2 + (\epsilon^{ab} \slnb_a \alpha_b)^2 \\
		|N\Phi_0|^2 &=&  |(M + \cot(\theta))\Phi_{-1}|^2 = (\slnb^a \underline{\alpha}_a)^2 + (\epsilon^{ab} \slnb_a \underline{\alpha}_b)^2 \; ,
	\end{eqnarray*}
	thus applying (\ref{liegreaterdivpluscurl}) and (\ref{EPhi0andEF1}), we prove (\ref{EPhi0andEF3}) and the lemma.
\end{proof}

Those relations, along with all the bounds that we established, will allow us to complete the list of energy estimates required to prove the decay results. We shall use the following compact notation for summations of energies. For $Y=T$ or $K$ and $\mathbb{A}$ a set of (Killing) vector fields,
\begin{equation}\label{compactnotationofenergysum}
E_Y[\lie_{\mathbb{A}}^k F](t):=\sum_{X_1, \dots, X_k\in \mathbb{A}}\hspace{-0.5cm}E_Y[\lie_{X_1} \dots \lie_{X_k}F](t) \; .
\end{equation}

\begin{prop1}[Uniform bound on Conformal Energy]\label{uniformboundConfromalenergyProposition}
	Let $F$ be a Maxwell field and $n$ a non-negative integer, then
	\begin{eqnarray}
	\int\limits_{[0,+\infty[\times \R\times \mathcal{S}^2} t \chi_{trap} |\Phi_0|^2  \d t \d r_* \d^2 \omega &\le& C \left( E_K[\lie_{\mathbb{O}} F](0) + \sum_{k=0}^{5}E_T[\lie_{\mathbb{O}}^k F](0) \right)\; , \label{uniformboundontrappingtermMaxwell} \\
	E_K[\lie_{\mathbb{T}}^n F](t) &\le& C\left(\sum_{k=0}^{n+1}E_K[\lie_{\mathbb{T}}^k F](0) + \sum_{k=0}^{n+5} E_T[\lie_{\mathbb{T}}^k F](0)\right) \; . \label{uniformboundonconformalenergy}
	\end{eqnarray}
\end{prop1}

\begin{proof}
	Since we assume that $\Phi_0$ is of the form (\ref{nonstationarysolutions}), then using (\ref{angularbound})
	\begin{equation*}
	\int\limits_{[0,+\infty[\times \R\times \mathcal{S}^2} t \chi_{trap} |\Phi_0|^2  \d t \d r_* \d^2 \omega \le \int\limits_{[0,+\infty[\times \R\times \mathcal{S}^2} t \chi_{trap} |\slnb\Phi_0|^2  \d t \d r_* \d^2 \omega \; ,
	\end{equation*}
	which in turn by (\ref{unifomboundtrappingterm}) yields
	\begin{equation*}
	\int\limits_{[0,+\infty[\times \R\times \mathcal{S}^2} t \chi_{trap} |\Phi_0|^2  \d t \d r_* \d^2 \omega \le C(E_\mathcal{C}[\Phi_0](0)+E[\sldlt^2 \Phi_0](0)) \; ,
	\end{equation*}
	and (\ref{uniformboundontrappingtermMaxwell}) follows from (\ref{EPhi0andEF2}) and (\ref{EPhi0andEF3}).
	
	Applying the trapping lemma for $\lie_{X_1} \dots \lie_{X_n}F$ over $[0,t]$, then summing over $X_i \in \mathbb{T}$, (\ref{trappinginequality}) becomes
	\begin{equation*}
	E_K[\lie_{\mathbb{T}}^n F](t) - E_K[\lie_{\mathbb{T}}^n F](0) \le \sum_{X_1 \dots X_n \in \mathbb{T}}\int\limits_{[0,t]\times \R\times \mathcal{S}^2} t \chi_{trap} |\lie_{X_1} \dots \lie_{X_n}\Phi_0|^2  \d t \d r_* \d^2 \omega \; ,
	\end{equation*}
	which upon using (\ref{uniformboundontrappingtermMaxwell}) to bound the right hand side gives (\ref{uniformboundonconformalenergy}).
\end{proof}

\subsection{Decay Results}\label{sec:decayresults}

\subsubsection{Uniform Decay}
From the uniform bound on the integral of the trapping term, (\ref{uniformboundontrappingtermMaxwell}), the conformal energy is controlled on any achronal future oriented smooth hypersurface (i.e. a smooth hyper surface with future oriented casual normal) such that its union with $\Sigma_t=\{t\} \times \R \times \mathcal{S}^2$ is the boundary of an open submanifold of $\mathcal{N}$. To see why, let $S$ be such a hypersurface and $\mathcal{U}$ the  corresponding open subset of $[0,+\infty[_t \times \R_{r_*} \times \mathcal{S}^2_{\theta,\varphi}$, then by (\ref{divergencetheoremSEtensor}) and the last step in the proof of Lemma (\ref{trappinglemma}), the trapping lemma, we have
\begin{equation*}
\int_{\dl \mathcal{U}} \mathbf{T}_{ab}K^b \eta^a i_{\tau} \d^4x = \hf\int_{\mathcal{U}} {}^{(K)}\pi_{ab} \mathbf{T}^{ab} \d^4x \le \int\limits_{[0,+\infty[\times \R\times \mathcal{S}^2} t \chi_{trap} |\Phi_0|^2  \d t \d r_* \d^2 \omega \; .
\end{equation*}
So, by (\ref{uniformboundontrappingtermMaxwell}), we have a uniform bound for all such hypersurfaces $S$,
\begin{equation}\label{uniformboundonconfrmalenergyhypersurface}
E_K[F](S)=\int_S \mathbf{T}_{ab}K^b \eta^a i_{\tau} \d^4x \le C \left(\sum_{k=0}^{1} E_K[\lie_{\mathbb{O}}^k F](0) + \sum_{k=0}^{5}E_T[\lie_{\mathbb{O}}^k F](0) \right)\; .
\end{equation}

\begin{thm1}[Uniform Decay]\label{decayonhypersurfacesTHEOREM}
	Let $t_0 \ge0$ be a real parameter. Let $F$ with spin components $(\Phi_1,\Phi_0,\Phi_{-1})$, be a non-stationary finite energy solution of Maxwell's equations (\ref{Maxeq1}) and (\ref{Maxeq2}), that is, satisfying (\ref{nonstationarysolutions1}) and (\ref{nonstationarysolutions}). Let $S$ be any achronal future oriented smooth hypersurface, such that its union with $\Sigma_0=\{0\} \times \R \times \mathcal{S}^2$ is the boundary of an open submanifold of $\mathcal{N}$, and such that on $S$, $t\ge |r_*| + t_0$. If $F$ and its first five Lie derivatives with respect to $\mathbb{O}$ in (\ref{vectorfields}) have finite energies and finite conformal energies on $\Sigma_0$, then the energy of the Maxwell field on $S$, defined by (\ref{energyonasurface}) for $X=T=\dl_t$, decays like $t_0^{-2}$.
	
	In fact, there is a constant $C>0$ independent of $t_0,\; F,\; (t,r_*,\omega),$ and $S$, such that
	\begin{equation}\label{decayonhypersurfaces}
	E_T[F](S) \le {t_0}^{-2}C \left(\sum_{k=0}^{1} E_K[\lie_{\mathbb{O}}^k F](0) + \sum_{k=0}^{5}E_T[\lie_{\mathbb{O}}^k F](0) \right)\; .
	\end{equation}
\end{thm1}
\begin{proof}
	If we show that
	\begin{equation}\label{conformalenergydominatesenergyonsurface}
	E_T[F](S)\le t_0^{-2}E_K[F](S)
	\end{equation}
	then (\ref{uniformboundonconfrmalenergyhypersurface}) gives the desired decay. To prove (\ref{conformalenergydominatesenergyonsurface}), it is sufficient to show that
	\begin{equation*}
	\mathbf{T}_{ab}\eta_S^a(K^b-t_0^2 T^b) \ge 0 \; ,
	\end{equation*}
	which by the dominant energy condition (Corollary \ref{DEC}) and the remark after it, reduces to showing that $K^b-t_0^2 T^b$ is a future oriented causal vector field, that is, proving that
	\begin{equation*}
	t^2+r_*^2 - t_0^2\ge |2tr_*| \quad \textrm{i.e.} \quad (t-|r_*|)^2 \ge t_0^{2} \; ,
	\end{equation*}
	but this is true on $S$ by the hypothesis of the theorem.
\end{proof}

Anther way of looking at this results is by considering a collection of hypersurfaces $\{S_{t_0} ; t_0\ge 0\}$ indexed by $t_0$, and the theorem says that the energy flux across these surfaces decays as $t_0$ goes to infinity. If one takes the collection of parabolas $\{S_{t_0}: t=\sqrt{1+r_*^2} + t_0 ; t_0 \ge 0\}$, then we see that the energy decays as the surfaces approach the timelike infinity. This kind of decay is particularly useful in the construction of scattering theories on spacetimes (\cite{nicolas_conformal_2016}, \cite{mokdad_conformal_2016}).

\subsubsection{Pointwise Decay}

We shall divide the proof of the pointwise decay into lemmata. Let $F$ be a non-stationary finite energy solution of Maxwell's equations, and recall the norms defined in the begin of this section. We start by the following estimate.

\begin{lem1}
	Let $[{r_*}_1,{r_*}_2]$ be a compact interval of $r_*$. Then there is a constant $C_{({r_*}_1,{r_*}_2)}>0$ such that for all $t\ge0$,
	\begin{equation}\label{boundbyconformalenergy}
	\int\limits_{\{t\}\times ({r_*}_1,{r_*}_2)\times \mathcal{S}^2} t^2 |F|^2_{\mathbb{X}} \d r_*^2 \d \omega^2 \le C E_K[F](t)
	\end{equation}
\end{lem1}

\begin{proof}
	Let $a=2 \max\{|{r_*}_1|,|{r_*}_2|\}$. If $t \ge a$, then for all $r_* \in [{r_*}_1,{r_*}_2]$ we have $t\ge 2 |r_*|$ and so
	\begin{eqnarray*}
		t-r_* &\ge& \hf t \;  \Rightarrow \; (t - r_*)^2 \ge \frac{1}{4}t^2\; , \\
		t+r_* &\ge& \hf t \; \Rightarrow \; (t + r_*)^2 \ge \frac{1}{4}t^2  \; . \\
	\end{eqnarray*}
	Thus, for $t\ge a$, from the equivalence of norms discussed in the beginning of this section,
	\begin{eqnarray*}
		E_K[F](t) &=& \frac{1}{4}\int_{\Sigma_t} (t + r_*)^2|\Phi_1|^2 +  (t^2 + r_*^2)V|\Phi_0|^2 + (t-r_*)^2|\Phi_{-1}|^2 \d r_* \d^2 \omega \\
		&\ge& C_{({r_*}_1,{r_*}_2)} \int\limits_{\{t\}\times ({r_*}_1,{r_*}_2) \times \mathcal{S}^2}  t^2 \left( |\Phi_1|^2 +|\Phi_0|^2 +|\Phi_{-1}|^2 \right) \d r_* \d^2 \omega   \\
		&\ge&  C_{({r_*}_1,{r_*}_2)}  \int\limits_{\{t\}\times ({r_*}_1,{r_*}_2) \times \mathcal{S}^2} t^2 |F|^2_{\mathbb{X}} \d r_* \d^2 \omega \; .
	\end{eqnarray*}
	If for some $t' \in [0,a]$ we have $E_K[F](t')=0$, then since its integrant is a continuous function, it vanishes over $({r_*}_1,{r_*}_2) \times \mathcal{S}^2$, and so, $|\Phi_1(t')| = |\Phi_0(t')| = |\Phi_{-1}(t')|=0$. By the the uniqueness of the solution to the Cauchy problem and the linearity of Maxwell's equations, $F$ must be identically zero for all $t$, and the inequality trivially holds. Thus, for $F\ne 0$, $E_K[F](t)\ne 0$ for all $t$, and so the function
	\begin{equation*}
	(E_K[F](t))^{-1}\int\limits_{\{t\}\times ({r_*}_1,{r_*}_2) \times \mathcal{S}^2} t^2 |F|^2_{\mathbb{X}} \d r_* \d^2 \omega
	\end{equation*}
	is continuous and hence bounded over the interval $[0,a]$, which is what we want to prove.
\end{proof}

We will also need some Sobolev estimates which are direct consequences of classical results.

\begin{lem1}
	\begin{enumerate}
		\item Let $u$ be a smooth function on $\bar{U}$ with $U=(a,b)\subset \R$. Let $H^1(U)$ be a Sobolev space over $U$. Then,
		\begin{equation}\label{localsobolev1dim}
		\sup\limits_{\bar{U}} | u | \le C \| u \|_{H^1(U)}=C \left(\int\limits_{(a,b)} |u|^2 + |u'|^2 \d s\right)^{\hf}  \; .
		\end{equation}
		\item Let $u$ be a smooth function on the Sphere $\mathcal{S}^2$. Then,
		\begin{equation}\label{sphericalsobolev}
		\sup\limits_{\mathcal{S}^2} | u | \le C\left( \int_{\mathcal{S}^2} |u|^2 + |\slnb u|^2 + |\sldlt u|^2 \d^2 \omega \right)^{\hf}\; .
		\end{equation}
	\end{enumerate}
\end{lem1}
\begin{proof}
	Since $u$ is smooth, the first estimate is immediate from the classical Morrey's local inequality (see \cite{evans_partial_2010} for example). For the second estimate, we know from the Sobolev estimate on the Sphere that $H^2 (\mathcal{S}^2) \subset C^0(\mathcal{S}^2)$ with
	\begin{equation*}
	\sup\limits_{\mathcal{S}^2} | u | \le C\left( \int_{\mathcal{S}^2} |u|^2 + |\slnb^2 u|^2 \d^2 \omega \right)\; .
	\end{equation*}
	where $|\slnb^2 u|^2=|g^{ij}_{\mathcal{S}^2}g^{kl}_{\mathcal{S}^2} (\slnb^2 u)_{il} (\slnb^2 u)_{jk}|$. By the divergence theorem,
	\begin{equation*}
	\int_{\mathcal{S}^2} \sldlt(|\slnb u|^2) \d^2 \omega =0 \; .
	\end{equation*}
	Then using Bochner-Weitzenbock-Lichnerowicz formula,
	\begin{equation*}
	\hf \sldlt(|\slnb u|^2)= |\slnb^2 u|^2 + (\slnb u)\cdot (\slnb\sldlt u) + \textrm{Ric}_{\mathcal{S}^2}(\slnb u,\slnb u) \; ,
	\end{equation*}
	and applying an integration by parts then using the fact that $ \textrm{Ric}_{\mathcal{S}^2} \ge C g_{\mathcal{S}^2}$ for some $C > 0$, we get (\ref{sphericalsobolev}).
\end{proof}

We now state the pointwise decay and complete the proof.

\begin{thm1}[Pointwise Decay]\label{PointwisedecayTHEOREM}
	Let $[{r_*}_1, {r_*}_2]$ be a compact interval of $r_*$. Then there is a constant $C_{({r_*}_1, {r_*}_2)}>0$ such that if $F$ with spin components $(\Phi_1,\Phi_0,\Phi_{-1})$, is a non-stationary finite energy solution of Maxwell's equations (\ref{Maxeq1}) and (\ref{Maxeq2}), that is, satisfying (\ref{nonstationarysolutions1}) and (\ref{nonstationarysolutions}), then for all $t\ge0$, $r_* \in [{r_*}_1, {r_*}_2]$, $\omega=(\theta,\varphi) \in \mathcal{S}^2$,
	\begin{equation}\label{Pointwisedecay}
	|\Phi_1|+|\Phi_0| + |\Phi_{-1}| \le C_{({r_*}_1, {r_*}_2)} t^{-1} \left(\sum_{k=0}^{4}  E_K[\lie_{\mathbb{T}}^k F](0) + \sum_{k=0}^{8}E_T[\lie_{\mathbb{T}}^k F](0) \right) \; .
	\end{equation}
\end{thm1}
\begin{proof}
	In this proof, the constant $C$ may depend on the interval $[{r_*}_1, {r_*}_2]$. Since the covariant derivative and the Lie derivative of a smooth vector field with respect to another smooth vector field is again a smooth vector field, then it is a linear combination of the coordinate vector fields $\dl_t, \dlr, \dl_\theta ,$ and $\dl_\varphi$ with smooth coefficient functions. Thus, if we have $X$ and $Y$ in $\mathbb{X}$ from (\ref{vectorfields}), then $\lie_X Y$ and $\nabla_X Y$ will be independent of $t$, and over $[{r_*}_1, {r_*}_2]\times \mathcal{S}^2$ their coefficient functions in the the frame  $(\dl_t, \dlr, \dl_\theta,\dl_\varphi)$ will be bounded by a constant $C$ depending only on $[{r_*}_1, {r_*}_2]\times \mathcal{S}^2$. By the paragraph just before (\ref{normalizedTheta3}), we see that this is also true in the frame $\mathbb{X}$.
	
	To prove the decay, we need the control
	\begin{eqnarray}
	|\lie_{R}\lie_{T_1}\lie_{T_2}(F(X_1,X_2))|&\le& C |F|_{\mathbb{X},3,\mathbb{T}}\label{tobecontrolled} \; , \\
	|\lie_{R}\lie_{T_1}(F(X_1,X_2))|&\le& C |F|_{\mathbb{X},2,\mathbb{T}}\label{tobecontrolled2} \; , \\
	|\lie_{R}(F(X_1,X_2))|&\le& C |F|_{\mathbb{X},1,\mathbb{T}}\label{tobecontrolled1} \; ,
	\end{eqnarray}
	for all $T_1,T_2 \in \mathbb{T}$ and $X_1,X_2\in\mathbb{X}$.
	
	Let $X_1,X_2,X_3 \in \mathbb{X}$, then
	\begin{equation}\label{0reestimatedecay}
	X_1(F(X_2,X_3))=(\lie_{X_1} F)(X_2,X_3) + F(\lie_{X_1} X_2,X_3) + F(X_2,\lie_{X_1} X_3) \; ,
	\end{equation}
	and so,
	\begin{equation}\label{1stestimatedecay1}
	| X_1(F(X_2,X_3))| \le C \left( |(\lie_{X_1} F)(X_2,X_3)|  +  |F|_{\mathbb{X}}\right) \le C |F|_{\mathbb{X},1,\mathbb{X}} \; .
	\end{equation}
	If we apply $X_4\in\mathbb{X}$ to (\ref{0reestimatedecay}), then expand $X_4(\lie_{X_1} F)(X_2,X_3)$ as in (\ref{0reestimatedecay}) with $\lie_{X_1} F$ replacing $F$, then by (\ref{1stestimatedecay1}) we get
	\begin{equation*}
	|X_4 X_1(F(X_2,X_3))| \le C  |F|_{\mathbb{X},2,\mathbb{X}}\; .
	\end{equation*}
	Repeating this process $n$ times, we have the iterated estimate
	\begin{equation*}
	|\lie_{X_1} \dots \lie_{X_k} (F(X,Y))| \le C |F|_{\mathbb{X},k,\mathbb{X}} \; .
	\end{equation*}
	In particular, if we choose $k$ vectors $X_i \in \mathbb{T}$ and $X,Y \in \mathbb{X}$, we obtain
	\begin{equation}\label{2ndestimatedecay}
	|\lie_{X_1} \dots \lie_{X_k} (F(X,Y))| \le C |F|_{\mathbb{X},k,\mathbb{T}} \; .
	\end{equation}
	
	Similar to (\ref{0reestimatedecay}), we also have for $X_1,X_2,X_3 \in \mathbb{X}$,
	\begin{equation}\label{6thestimatedecay}
	X_1(F(X_2,X_3))=(\nabla_{X_1} F)(X_2,X_3) + F(\nabla_{X_1} X_2,X_3) + F(X_2,\nabla_{X_1} X_3) \; ,
	\end{equation}
	
	If we take $T_1,T_2\in \mathbb{T}$ then by applying $\lie_{T_1}\lie_{T_2}$ to (\ref{6thestimatedecay}) with $T_3 \in \mathbb{T}$ instead of $X_1\in\mathbb{X}$, keeping $X_2,X_3 \in \mathbb{X}$, we get
	\begin{eqnarray*}
		\lie_{T_1}\lie_{T_2}({T_3}(F(X_2,X_3)))&=&\lie_{T_1}\lie_{T_2}((\nabla_{T_3} F)(X_2,X_3)) \\
		&& \hspace{1cm}+ \lie_{T_1}\lie_{T_2}(F(\nabla_{T_3} X_2,X_3)) + \lie_{T_1}\lie_{T_2}(F(X_2,\nabla_{T_3} X_3)) \; .
	\end{eqnarray*}
	By (\ref{2ndestimatedecay}), the left hand side and the last two terms of the right hand side are controlled by $|F|_{\mathbb{X},3,\mathbb{T}}$ and  $|F|_{\mathbb{X},2,\mathbb{T}}$ respectively, and hence
	\begin{equation}\label{7thestimatedecay}
	|\lie_{T_1}\lie_{T_2}((\nabla_{T_3} F)(X_2,X_3))|\le C  |F|_{\mathbb{X},3,\mathbb{T}} \; .
	\end{equation}
	Similarly one has
	\begin{eqnarray}
	|\lie_{T_1}((\nabla_{T_3} F)(X_2,X_3))|&\le& C  |F|_{\mathbb{X},2,\mathbb{T}} \; . \nonumber \\
	|(\nabla_{T_3} F)(X_2,X_3)| &\le& C  |F|_{\mathbb{X},1,\mathbb{T}} \; . \label{7thestimatedecaylowerorder}
	\end{eqnarray}
	
	Let $T_1 ,\dots, T_4 \in \mathbb{T}$, then starting from (\ref{6thestimatedecay}) again, we apply two Lie derivatives and obtain
	\begin{eqnarray*}
		\lie_{T_1}\lie_{T_2}({R}(F(T_3,T_4)))&=&\lie_{T_1}\lie_{T_2}((\nabla_{R} F)(T_3,T_4)) \\
		&& \hspace{1cm}+ \lie_{T_1}\lie_{T_2}(F(\nabla_{R} T_3,T_4)) + \lie_{T_1}\lie_{T_2}(F(T_3,\nabla_{R} T_4)) \; .
	\end{eqnarray*}
	From (\ref{Maxeqabst2}),
	\begin{equation*}
	(\nabla_{R} F)(T_3,T_4)=-(\nabla_{T_3} F)( R,T_4)-(\nabla_{T_4} F)(T_3,R) \; .
	\end{equation*}
	and so we have,
	\begin{eqnarray*}
		|\lie_{T_1}\lie_{T_2}({R}(F(T_3,T_4)))|&\le& |\lie_{T_1}\lie_{T_2}((\nabla_{T_3} F)( R,T_4))| +|\lie_{T_1}\lie_{T_2}((\nabla_{T_4} F)(T_3,R))|+ \\
		&& \hspace{1cm}+ |\lie_{T_1}\lie_{T_2}(F(\nabla_{R} T_3,T_4))| + |\lie_{T_1}\lie_{T_2}(F(T_3,\nabla_{R} T_4))| \; ,
	\end{eqnarray*}
	and all the terms on the right hand side are controlled by (\ref{7thestimatedecay}) and (\ref{2ndestimatedecay}). Thus,
	\begin{equation*}
	|\lie_{T_1}\lie_{T_2}({R}(F(T_3,T_4)))| \le C  |F|_{\mathbb{X},3,\mathbb{T}}\; .
	\end{equation*}
	Similarly, we have,
	\begin{eqnarray*}
		|\lie_{T_1}({R}(F(T_3,T_4)))| &\le& C  |F|_{\mathbb{X},2,\mathbb{T}}\; , \\
		|{R}(F(T_3,T_4))| &\le& C  |F|_{\mathbb{X},1,\mathbb{T}}\; .
	\end{eqnarray*}
	
	Next, choose an a local orthonormal frame $\{rV,rW\}$ on $U\subset\mathcal{S}^2$, so that $$\{f^{-\hf}\dl_t, f^{-\hf}\dlr , V ,W\}$$ is a local orthonormal frame on the spacetime. In this frame, for $X\in\mathbb{X}$, (\ref{Maxeqabst1}) is
	\begin{equation}\label{orthonormalframceS2}
	f^{-1}\left((\nabla_{\dl_t}F)(\dl_t,X) - (\nabla_{R}F)(R,X)\right) -(\nabla_{V}F)(V,X)- (\nabla_{W}F)(W,X)=0\; .
	\end{equation}
	Shrinking $U$ if necessary, $V$ and $W$ can be locally written as linear combinations of vectors in $\mathbb{O}$ with smooth coefficients, and since $\mathcal{S}^2$ is locally compact, the coefficients and all their derivatives are bounded in some neighbourhood around each point. By covering $\mathcal{S}^2$ with such $U$'s, and as $\mathcal{S}^2$ is compact, a finite subcover of all these neighbourhoods, means that we have a uniform bound for all these coefficients defined on the subcover. Applying $\lie_{T_1}\lie_{T_2}$, for $T_1,T_2\in\mathbb{T}$, to (\ref{orthonormalframceS2}) with $V$ and $W$ replaced by there expressions in terms of $\Theta_i$'s with bounded coefficients, then using (\ref{7thestimatedecay}) and (\ref{7thestimatedecaylowerorder}), we get the control
	\begin{equation*}
	|\lie_{T_1}\lie_{T_2}((\nabla_{R} F)(R,X))|\le C  |F|_{\mathbb{X},3,\mathbb{T}} \; .
	\end{equation*}
	This control, by virtue of (\ref{6thestimatedecay}) gives the required control
	\begin{equation*}
	|\lie_{T_1}\lie_{T_2}(R (F(R,X)))|\le C  |F|_{\mathbb{X},3,\mathbb{T}} \; .
	\end{equation*}
	Similarly,
	\begin{eqnarray*}
		|\lie_{T_1}(R (F(R,X)))|&\le& C  |F|_{\mathbb{X},2,\mathbb{T}} \; , \\
		|R (F(R,X))|&\le& C  |F|_{\mathbb{X},1,\mathbb{T}} \; .
	\end{eqnarray*}
	This covers (\ref{tobecontrolled})-(\ref{tobecontrolled1}).
	

	Let $X_1,X_2\in\mathbb{X}$, from (\ref{sphericalsobolev}),
	\begin{equation*}
	|F(X_1,X_2)|^2 \le C\int_{\mathcal{S}^2}  |F(X_1,X_2)|^2 +  |\slnb(F(X_1,X_2))|^2 +  |\sldlt (F(X_1,X_2))|^2 \d^2 \omega \; .
	\end{equation*}
	Then using (\ref{rotationformofslnb}) and (\ref{rotationformofsldlt}),
	\begin{eqnarray*}
		|F(X_1,X_2)|^2 \!\!\!&\le& \!\!\!C\int_{\mathcal{S}^2}  |F(X_1,X_2)|^2 +  \sum_{i=1}^{3}\left(|\Theta_i(F(X_1,X_2))|^2 +  |\Theta_i^2 (F(X_1,X_2))|^2\right) \d^2 \omega \\
		\!\!\! &\le&\!\!\!C\int_{\mathcal{S}^2}  |F(X_1,X_2)|^2 + \!\!\!\! \sum_{T_1,T_,2\in\mathbb{T}}\!\!\!\!\left(|\lie_{T_1}(F(X_1,X_2))|^2 +  |\lie_{T_1}\lie_{T_2} (F(X_1,X_2))|^2\right) \d^2 \omega \; .
	\end{eqnarray*}
	By (\ref{localsobolev1dim}),
	\begin{eqnarray*}
		\hspace{-1cm}\int\limits_{\{t\}\times \{r_*\}\times \mathcal{S}^2}\hspace{-0.5cm} |\lie_{T_1}\lie_{T_2}(F(X_1,X_2))|^2 \d^2 \omega &\le& C\hspace{-0.5cm} \int\limits_{\{t\}\times ({r_*}_1, {r_*}_2)\times \mathcal{S}^2} \hspace{-0.5cm} |\lie_{T_1}\lie_{T_2}(F(X_1,X_2))|^2 \\
		&&\hspace{2cm}+ |\lie_{T_1}\lie_{T_2}(R(F(X_1,X_2)))|^2 \d r_* \d^2 \omega \; , \\
		\int\limits_{\{t\}\times \{r_*\}\times \mathcal{S}^2}\hspace{-0.5cm} |\lie_{T_1}(F(X_1,X_2))|^2 \d^2 \omega &\le& C\hspace{-0.5cm} \int\limits_{\{t\}\times ({r_*}_1, {r_*}_2)\times \mathcal{S}^2} \hspace{-0.5cm} |\lie_{T_1}(F(X_1,X_2))|^2 \\
		&&\hspace{2cm}+ |\lie_{T_1}(R(F(X_1,X_2)))|^2 \d r_* \d^2 \omega \; , \\
		\int\limits_{\{t\}\times \{r_*\}\times \mathcal{S}^2}\hspace{-0.5cm}
		|(F(X_1,X_2))|^2 \d^2 \omega &\le& C\hspace{-0.5cm} \int\limits_{\{t\}\times ({r_*}_1, {r_*}_2)\times \mathcal{S}^2} \hspace{-0.5cm} |(F(X_1,X_2))|^2 + |(R(F(X_1,X_2)))|^2 \d r_* \d^2 \omega \; . \\
	\end{eqnarray*}
	Summing over all $X_1,X_2\in\mathbb{X}$, then using (\ref{tobecontrolled}), (\ref{tobecontrolled2}), (\ref{tobecontrolled1}), and (\ref{2ndestimatedecay}), we have
	\begin{equation*}
	|F|^2_{\mathbb{X}} \le C \int\limits_{\{t\}\times ({r_*}_1, {r_*}_2)\times \mathcal{S}^2} |F|^2_{\mathbb{X},3,\mathbb{T}} \d r_* \d^2 \omega \; .
	\end{equation*}
	Finally, since for $T_1,\dots,T_k\in\mathbb{T}$, $\lie_{T_1}\dots\lie_{T_k}F$ is again a Maxwell field, then applying (\ref{boundbyconformalenergy}) for $\lie_{T_1}\dots\lie_{T_k}F$, we have
	\begin{equation*}
	\int\limits_{\{t\}\times ({r_*}_1, {r_*}_2)\times \mathcal{S}^2} t^2|\lie_{T_1}\dots\lie_{T_k}F|^2_{\mathbb{X}}  \d r_* \d^2 \omega \le C E_K[\lie_{T_1}\dots\lie_{T_k}F](t) \; ,
	\end{equation*}
	Taking summation over all $T_i\in\mathbb{T}$ and over $k$ from $0$ to $3$, we obtain,
	\begin{equation*}
	\int\limits_{\{t\}\times ({r_*}_1, {r_*}_2)\times \mathcal{S}^2} t^2 |F|^2_{\mathbb{X},3,\mathbb{T}} \le C \sum_{k=0}^{3} E_K[\lie_{\mathbb{T}}^k F](t) \; ,
	\end{equation*}
	and therefore,
	\begin{equation*}
	|F|^2_{\mathbb{X}} \le C t^{-2} E_K[\lie_{\mathbb{T}}^k F](t) \; ,
	\end{equation*}
	and the decay follows from (\ref{uniformboundonconformalenergy}).
\end{proof}

\subsection{Generic Spherically Symmetric Black Hole Spacetimes}\label{genericsphericalf}

We conclude these results by a remark on generic spherically symmetric black hole spacetimes, and specify what conditions on the metric are needed for our results to hold on a more general spacetime than RNdS.

Let the manifold $\mathcal{M}=\R_t \times ]0,+\infty[_r \times \mathcal{S}^2_{\theta,\varphi}$ be equipped with the metric defined in (\ref{RNdSmetric}), with $f$ a smooth function on $]0,+\infty[_r$. Consider the following conditions on $f$:
\begin{enumerate}
	\item There are three real numbers $r_i$, $0\le r_1 < r_2 < r_3 \le +\infty$, and these are the only possible zeros of  $f$.
	\item If $0< r_i <+\infty$ then $f(r_i)=0$ and $f'(r_i)\ne0$.
	\item $f(r)<0$ for $r\in ] r_1 , r_2 [$ and $f(r)>0$ for $r\in ] r_2 , r_3 [$
	\item If $r_3$ is infinite then $f(r)=1-Cr^{-1}+O(r^{-2})$ for $C>0$, and for $k=1,2,3$, $\dl_r^k f (r)=O(r^{-k-1})$, as $r\rightarrow +\infty$.
	\
\end{enumerate}
All these properties are satisfied by usual spherical black holes, like Schwarzschild, Reissner-Nordstrøm, Yasskin, asymptotically flat, or De Sitter, etc \dots

If $f$ indeed satisfies these conditions, then all of the above decay results hold equally true for such a generic spherically symmetric spacetimes. Our arguments make no use of the particular form of $f$ in RNdS, but rather these properties of $f$, and so the required arguments for a general spacetime are the same.

\section*{Acknowledgement}
The results of this paper and the mentioned conformal scattering results on RNDS spacetime \cite{mokdad_reissner-nordstrom-sitter_2017,mokdad_conformal_2016}, were obtained during my PhD thesis \cite{mokdad_maxwell_2016}. I would like to thank my thesis advisor Pr. Jean-Philippe Nicolas for his indispensable guidance during the thesis.

\printbibliography[heading=bibintoc]

\end{document}